# Some series and integrals involving the Riemann zeta function, binomial coefficients and the harmonic numbers

## Volume II(b)

Donal F. Connon

18 February 2008

**Abstract**

In this series of seven papers, predominantly by means of elementary analysis, we establish a number of identities related to the Riemann zeta function, including the following:

$$\gamma_p(u) = -\frac{1}{p+1}\sum_{n=0}^{\infty}\frac{1}{n+1}\sum_{k=0}^{n}\binom{n}{k}(-1)^k \log^{p+1}(u+k)$$

$$\int_1^x \gamma_p(u)\,du = \frac{(-1)^{p+1}}{p+1}[\varsigma^{(p+1)}(0,x) - \varsigma^{(p+1)}(0)]$$

$$\gamma_1\left(\frac{1}{4}\right) = \frac{1}{2}[2\gamma_1 - 15\log^2 2 - 6\gamma \log 2] - \frac{1}{2}\pi\left[\gamma + 4\log 2 + 3\log \pi - 4\log \Gamma\left(\frac{1}{4}\right)\right]$$

$$\int_1^\infty \left[\frac{1}{1-t} - \frac{1}{t\log(1/t)}\right]\frac{\log\log(1/t)}{t^u}\,dt = \gamma\gamma_0(u) + 2\gamma_1(u) + \gamma \log u + \log^2 u$$

$$\sum_{n=0}^{\infty}\frac{1}{n+1}\sum_{k=0}^{n}\binom{n}{k}(-1)^k k \log^2(k+1) = \gamma_1 - \log(2\pi) - \frac{1}{2}\gamma^2 + \frac{1}{24}\pi^2 + \frac{1}{2}\log^2(2\pi)$$

$$\sum_{n=1}^{\infty}t^n \sum_{k=1}^{n}\binom{n}{k}\frac{x^k}{(k+y)^s} = \frac{xt}{(1-t)\Gamma(s)}\int_0^\infty \frac{u^{s-1}e^{-(y+1)u}}{[1-(1+xe^{-u})t]}\,du = \frac{xt}{(1-t)^2}\Phi\left(\frac{xt}{(1-t)},s,y+1\right)$$

where $\gamma_p(u)$ is the Stieltjes constant and $\Phi(z,s,u)$ is the Hurwitz-Lerch zeta function.

Whilst this paper is mainly expository, some of the formulae reported in it are believed to be new, and the paper may also be of interest specifically due to the fact that most of the various identities have been derived by elementary methods.

**CONTENTS OF VOLUMES I TO VI:**  Volume/page

**SECTION:**

**ACKNOWLEDGEMENTS**

**REFERENCES**

# AN INTRODUCTION TO THE STIELTJES CONSTANTS

We recall Hasse's formula for the Hurwitz zeta function from (4.3.106a) that for $\operatorname{Re}(s) \neq 1$

$$(s-1)\varsigma(s,u) = \sum_{n=0}^{\infty}\frac{1}{n+1}\sum_{k=0}^{n}\binom{n}{k}\frac{(-1)^k}{(u+k)^{s-1}}$$

and note that for $s > 1$

$$(s-1)\sum_{n=0}^{\infty}\frac{1}{(n+u)^s} = \sum_{n=0}^{\infty}\frac{1}{n+1}\sum_{k=0}^{n}\binom{n}{k}\frac{(-1)^k}{(u+k)^{s-1}}$$

Hence, in passing, we see that $\sum_{n=0}^{\infty}\frac{1}{n+1}\sum_{k=0}^{n}\binom{n}{k}\frac{(-1)^k}{(u+k)^{s-1}}$ is positive for $s > 1$ and is monotonically decreasing with $u$ (which is not intuitively obvious).

We have the limit

(4.3.200) $$\lim_{s\to 1}(s-1)\varsigma(s,u) = \sum_{n=0}^{\infty}\frac{1}{n+1}\sum_{k=0}^{n}\binom{n}{k}(-1)^k = \sum_{n=0}^{\infty}\frac{1}{n+1}\delta_{n,0} = 1$$

This limit is independent of $u$ and we therefore have

$$\frac{\partial}{\partial u}\lim_{s\to 1}[(s-1)\varsigma(s,u)] = \lim_{s\to 1}\frac{\partial}{\partial u}[(s-1)\varsigma(s,u)] = 0$$

We also note that

(4.3.200a) $$\frac{\partial}{\partial u}\varsigma(s,u) = \frac{\partial}{\partial u}\sum_{n=0}^{\infty}\frac{1}{(n+u)^s} = -s\varsigma(s+1,u)$$

and hence as before we have

$$\lim_{s\to 1}\frac{\partial}{\partial u}[(s-1)\varsigma(s,u)] = -\lim_{s\to 1}[s(s-1)\varsigma(s+1,u)] = 0$$

As a result of (4.3.200) we may apply L'Hôpital's rule to the limit

$$\lim_{s\to 1}\left[\frac{(s-1)\varsigma(s,u)-1}{s-1}\right] = \lim_{s\to 1}\left[\frac{(s-1)\varsigma'(s,u)+\varsigma(s,u)}{1}\right]$$

From (4.3.107a) we see that

$$(4.3.201) \quad (s-1)\varsigma'(s,u)+\varsigma(s,u) = -\sum_{n=0}^{\infty}\frac{1}{n+1}\sum_{k=0}^{n}\binom{n}{k}(-1)^k \frac{\log(u+k)}{(u+k)^{s-1}}$$

and hence we have

$$\lim_{s\to 1}\left[(s-1)\varsigma'(s,u)+\varsigma(s,u)\right] = -\sum_{n=0}^{\infty}\frac{1}{n+1}\sum_{k=0}^{n}\binom{n}{k}(-1)^k \log(u+k)$$

Referring back to (4.3.74) in Volume II(a) we may then deduce that

$$(4.3.202) \quad \lim_{s\to 1}\left[(s-1)\varsigma'(s,u)+\varsigma(s,u)\right] = -\psi(u)$$

and therefore conclude that

$$(4.3.203) \quad \lim_{s\to 1}\left[\frac{(s-1)\varsigma(s,u)-1}{s-1}\right] = \lim_{s\to 1}\left[\varsigma(s,u)-\frac{1}{s-1}\right] = -\psi(u)$$

The above limit is derived in a more complex manner in [126, p.91] and [135, p.271] using Hermite's integral formula for the Hurwitz zeta function. A further derivation is contained in Volume VI.

With $u=1$ we immediately see that

$$(4.3.204) \quad \lim_{s\to 1}\left[\varsigma(s,1)-\frac{1}{s-1}\right] = \lim_{s\to 1}\left[\varsigma(s)-\frac{1}{s-1}\right] = -\psi(1) = \gamma$$

which is derived in an alternative way in (4.4.99m) in Volume III. The above limits then imply that

$$(4.3.204a) \quad \lim_{s\to 1}\left[\varsigma(s,u)-\varsigma(s)\right] = -\gamma - \psi(u)$$

In [135, p.266] we see that

$$(4.3.204b) \quad \lim_{s\to 1}\left[\frac{\varsigma(s,u)}{\Gamma(1-s)}\right] = -1$$

This is easily demonstrated as follows. We have

$$\frac{\varsigma(s,u)}{\Gamma(1-s)} = \frac{(1-s)\varsigma(s,u)}{(1-s)\Gamma(1-s)}$$

and hence we see that

$$\lim_{s\to 1}\left[\frac{(1-s)\varsigma(s,u)}{(1-s)\Gamma(1-s)}\right] = \lim_{s\to 1}\left[\frac{(1-s)\varsigma(s,u)}{\Gamma(2-s)}\right] = -1$$

L'Hôpital's rule then gives us

(4.3.205) $$\lim_{s\to 1}\left[\frac{\varsigma'(s,u)}{\Gamma'(1-s)}\right] = 1$$

Differentiating (4.3.202) with respect to $u$ gives us

$$\frac{d}{du}\lim_{s\to 1}\left[\varsigma(s,u) - \frac{1}{s-1}\right] = -\psi'(u)$$

and hence

$$\frac{d}{du}\lim_{s\to 1}\varsigma(s,u) = -\psi'(u)$$

We have by interchanging the derivative and the limit operations

$$\frac{d}{du}\lim_{s\to 1}\varsigma(s,u) = \lim_{s\to 1}\frac{\partial}{\partial u}\varsigma(s,u)$$

$$= \lim_{s\to 1}[-s\varsigma(s+1,u)] = -\varsigma(2,u)$$

giving the well-known result $\varsigma(2,u) = \psi'(u)$

This may also be derived by differentiating (4.3.200a).

From (4.3.203) we see that

$$\lim_{s\to 1}\left[\varsigma(s,u) - \frac{1}{s-1} + \psi(u)\right] = 0$$

and we have the Laurent expansion of the Hurwitz zeta function $\varsigma(s,u)$ about $s=1$

(4.3.206) $$\varsigma(s,u) = \frac{1}{s-1} + \sum_{p=0}^{\infty}\frac{(-1)^p}{p!}\gamma_p(u)(s-1)^p$$

where $\gamma_p(u)$ are known as the generalised Stieltjes constants. We have

(4.3.206a) $\gamma_0(u) = -\psi(u)$ and $\gamma_0(1) = -\psi(1) = \gamma$.

We have for $p \geq 0$

(4.3.207) $$\gamma_p(u) = (-1)^p \lim_{s \to 1}\left[\varsigma^{(p)}(s,u) - \frac{(-1)^p p!}{(s-1)^{p+1}}\right]$$

and, in particular, we have

(4.3.208) $$\gamma_0(u) = \lim_{s \to 1}\left[\varsigma(s,u) - \frac{1}{s-1}\right] = -\psi(u)$$

$$\gamma_1(u) = -\lim_{s \to 1}\left[\varsigma'(s,u) + \frac{1}{(s-1)^2}\right]$$

Applying L'Hôpital's rule again gives us

(4.3.208a) $$\lim_{s \to 1}\left[\frac{\varsigma(s,u) - \frac{1}{s-1} + \psi(u)}{s-1}\right] = \lim_{s \to 1}\left[\varsigma'(s,u) + \frac{1}{(s-1)^2}\right] = -\gamma_1(u)$$

Using (4.3.202) we see that

$$\lim_{s \to 1}[(s-1)\varsigma'(s,u) + \varsigma(s,u) + \psi(u)] = 0$$

and applying L'Hôpital's rule gives us

$$\lim_{s \to 1}\left[\frac{(s-1)\varsigma'(s,u) + \varsigma(s,u) + \psi(u)}{s-1}\right] = \lim_{s \to 1}\left[\frac{(s-1)\varsigma''(s,u) + 2\varsigma'(s,u)}{1}\right]$$

From (4.3.201) we see that

$$\frac{\partial}{\partial s}[(s-1)\varsigma'(s,u) + \varsigma(s,u)] = \sum_{n=0}^{\infty} \frac{1}{n+1} \sum_{k=0}^{n} \binom{n}{k}(-1)^k \frac{\log^2(k+u)}{(u+k)^{s-1}}$$

and thus we have

(4.3.209) $$\lim_{s \to 1}[(s-1)\varsigma''(s,u) + 2\varsigma'(s,u)] = \sum_{n=0}^{\infty} \frac{1}{n+1} \sum_{k=0}^{n} \binom{n}{k}(-1)^k \log^2(k+u)$$

Accordingly we see that

$$\lim_{s \to 1}\left[\frac{(s-1)\varsigma'(s,u)+\varsigma(s,u)+\psi(u)}{s-1}\right]=\sum_{n=0}^{\infty}\frac{1}{n+1}\sum_{k=0}^{n}\binom{n}{k}(-1)^k \log^2(k+u)$$

and the limit may be written as

$$\lim_{s \to 1}\left[\varsigma'(s,u)+\frac{1}{(s-1)^2}+\frac{\varsigma(s,u)+\psi(u)}{s-1}-\frac{1}{(s-1)^2}\right]$$

$$=\lim_{s \to 1}\left[\varsigma'(s,u)+\frac{1}{(s-1)^2}\right]+\lim_{s \to 1}\left[\frac{\varsigma(s,u)+\psi(u)}{s-1}-\frac{1}{(s-1)^2}\right]$$

$$=-2\gamma_1(u)$$

where we have used (4.3.207) and (4.3.208a).

Hence we obtain

(4.3.210) $$\gamma_1(u)=-\frac{1}{2}\sum_{n=0}^{\infty}\frac{1}{n+1}\sum_{k=0}^{n}\binom{n}{k}(-1)^k \log^2(k+u)$$

The corresponding Stieltjes constants $\gamma_p = \gamma_p(1)$ for the Riemann zeta function $\varsigma(s)$ are given by

$$\gamma_p = (-1)^p \lim_{s \to 1}\left[\varsigma^{(p)}(s)-\frac{(-1)^p p!}{(s-1)^{p+1}}\right]$$

and hence we have

(4.3.211) $$\gamma_0(1)=\gamma_0=\gamma=-\sum_{n=0}^{\infty}\frac{1}{n+1}\sum_{k=0}^{n}\binom{n}{k}(-1)^k \log(k+1)$$

(where we have employed the Ser/Sondow formula (4.4.92a) for $\gamma$ noting that no contribution comes from the term in (4.4.92a) with $n = 0$). Furthermore we see that

(4.3.212) $$\gamma_1(1)=\gamma_1=-\frac{1}{2}\sum_{n=0}^{\infty}\frac{1}{n+1}\sum_{k=0}^{n}\binom{n}{k}(-1)^k \log^2(k+1)$$

We note from [The On-Line Encyclopedia of Integer Sequences](#) that

$$\gamma_1 = -0.728158458...$$

On suspects from (4.3.212) that a pattern is developing and this is investigated below in a more generalised setting.

Suppose that we have the following functional equation for suitably analytic functions

$$(s-a)f(s,u) = g(s,u)$$

and suppose that $\lim_{s \to a} g(s,u) = g(a,u)$ exists. We then have

$$\lim_{s \to a}\left[(s-a)f(s,u) - g(a,u)\right] = 0$$

Applying L'Hôpital's rule results in

$$\lim_{s \to a}\left[\frac{(s-a)f(s,u) - g(a,u)}{s-a}\right] = \lim_{s \to a}\frac{\partial}{\partial s}\left[(s-a)f(s,u) - g(a,u)\right]$$

$$= \lim_{s \to a}\frac{\partial}{\partial s} g(s,u) = g'(a,u)$$

Therefore we have

$$\lim_{s \to a}\left[f(s,u) - \frac{g(a,u)}{s-a}\right] = \lim_{s \to a}\frac{\partial}{\partial s} g(s,u) = g'(a,u)$$

and this may be written as

$$\lim_{s \to a}\left[f(s,u) - \frac{g(a,u)}{s-a} - g'(a,u)\right] = \lim_{s \to a}\frac{\partial}{\partial s} g(s,u) - g'(a,u) = 0$$

A further application of L'Hôpital's rule produces

$$\lim_{s \to a}\left[\frac{f(s,u) - \frac{g(a,u)}{s-a} - g'(a,u)}{s-a}\right] = \lim_{s \to a}\left[f'(s,u) + \frac{g(a,u)}{(s-a)^2}\right]$$

$$= \frac{\partial}{\partial s}\left[\lim_{s \to a}\frac{\partial}{\partial s} g(s,u) - g'(a,u)\right] = \lim_{s \to a}\frac{\partial^2}{\partial s^2} g(s,u)$$

Therefore we may deduce that

$$\lim_{s \to a}\left[ f'(s,u) + \frac{g(a,u)}{(s-a)^2} \right] = \lim_{s \to a} \frac{\partial^2}{\partial s^2} g(s,u)$$

and the process may be continued indefinitely to result in

$$\lim_{s \to a}\left[ f^{(p)}(s,u) - \frac{(-1)^p p!\, g(a,u)}{(s-a)^{p+1}} \right] = \lim_{s \to a} \frac{\partial^{p+1}}{\partial s^{p+1}} g(s,u)$$

as may be readily proved by mathematical induction.

We also have where $f(s,u)$ is simply of the form $f(s)$

(4.3.213) $$\lim_{s \to a}\left[ f^{(p)}(s) - \frac{(-1)^p p!\, g(a)}{(s-a)^{p+1}} \right] = \lim_{s \to a} \frac{d^{p+1}}{ds^{p+1}} g(s)$$

In the case of the Hurwitz zeta function we let

$$f(s,u) = \varsigma(s,u)$$

and, with reference to (4.3.106a), we see that

$$g(s,u) = \sum_{n=0}^{\infty} \frac{1}{n+1} \sum_{k=0}^{n} \binom{n}{k} \frac{(-1)^k}{(u+k)^{s-1}}$$

$$g(1,u) = \sum_{n=0}^{\infty} \frac{1}{n+1} \sum_{k=0}^{n} \binom{n}{k} (-1)^k = 1$$

$$\frac{\partial^{p+1}}{\partial s^{p+1}} g(s,u) = (-1)^{p+1} \sum_{n=0}^{\infty} \frac{1}{n+1} \sum_{k=0}^{n} \binom{n}{k} \frac{(-1)^k \log^{p+1}(u+k)}{(u+k)^{s-1}}$$

$$\lim_{s \to 1}\left[ \varsigma^{(p)}(s,u) - \frac{(-1)^p p!}{(s-1)^{p+1}} \right] = (-1)^{p+1} \sum_{n=0}^{\infty} \frac{1}{n+1} \sum_{k=0}^{n} \binom{n}{k} (-1)^k \log^{p+1}(1+k)$$

From (4.3.207) we had

$$\gamma_p(u) = (-1)^p \lim_{s \to 1}\left[ \varsigma^{(p)}(s,u) - \frac{(-1)^p p!}{(s-1)^{p+1}} \right]$$

and therefore we may determine the following expression for the Stieltjes constants

$$\gamma_p(u) = -\frac{1}{p+1}\sum_{n=0}^{\infty}\frac{1}{n+1}\sum_{k=0}^{n}\binom{n}{k}(-1)^k \log^{p+1}(u+k)$$

$$\gamma_p = \gamma_p(1) = -\frac{1}{p+1}\sum_{n=0}^{\infty}\frac{1}{n+1}\sum_{k=0}^{n}\binom{n}{k}(-1)^k \log^{p+1}(1+k)$$

$$\gamma_0(u) = -\sum_{n=0}^{\infty}\frac{1}{n+1}\sum_{k=0}^{n}\binom{n}{k}(-1)^k \log(u+k) = -\psi(u)$$

An alternative derivation of the series for the generalised Stieltjes constants is shown below. We have

$$\varsigma(s,u) = \frac{1}{s-1} + \sum_{p=0}^{\infty}\frac{(-1)^p}{p!}\gamma_p(u)(s-1)^p$$

and therefore

$$(s-1)\varsigma(s,u) = 1 + \sum_{p=0}^{\infty}\frac{(-1)^p}{p!}\gamma_p(u)(s-1)^{p+1}$$

Differentiation with respect to $s$ gives us

$$\frac{\partial}{\partial s}[(s-1)\varsigma(s,u)] = \sum_{p=0}^{\infty}\frac{(-1)^p}{p!}\gamma_p(u)(p+1)(s-1)^p = \gamma_0(u) + \sum_{p=1}^{\infty}\frac{(-1)^p}{p!}\gamma_p(u)(p+1)(s-1)^p$$

and as $s \to 1$ we have using (4.3.202)

$$-\psi(u) = \gamma_0(u)$$

A further differentiation gives us

$$\frac{\partial^2}{\partial s^2}[(s-1)\varsigma(s,u)] = \sum_{p=0}^{\infty}\frac{(-1)^p}{p!}\gamma_p(u)(p+1)p(s-1)^{p-1}$$

$$= -2\gamma_1(u) + \sum_{p=2}^{\infty}\frac{(-1)^p}{p!}\gamma_p(u)(p+1)p(s-1)^{p-1}$$

As $s \to 1$ we have

$$\lim_{s\to 1}[(s-1)\varsigma''(s,u) + 2\varsigma'(s,u)] = \sum_{n=0}^{\infty}\frac{1}{n+1}\sum_{k=0}^{n}\binom{n}{k}(-1)^k \log^2(k+u)$$

and therefore

$$\sum_{n=0}^{\infty}\frac{1}{n+1}\sum_{k=0}^{n}\binom{n}{k}(-1)^k \log^2(k+u) = -2\gamma_1(u)$$

More generally we have

$$\frac{\partial^m}{\partial s^m}[(s-1)\varsigma(s,u)] = \sum_{p=0}^{\infty}\frac{(-1)^p}{p!}\gamma_p(u)(p+1)p(p-1)...(p+2-m)(s-1)^{p+1-m}$$

$$= \sum_{p=m-1}^{\infty}\frac{(-1)^p}{p!}\gamma_p(u)(p+1)p(p-1)...(p+2-m)(s-1)^{p+1-m}$$

$$= m(-1)^{m-1}\gamma_{m-1}(u) + \sum_{p=m}^{\infty}\frac{(-1)^p}{p!}\gamma_p(u)(p+1)p(p-1)...(p+2-m)(s-1)^{p+1-m}$$

Therefore we see that

$$\lim_{s\to 1}\frac{\partial^m}{\partial s^m}[(s-1)\varsigma(s,u)] = m(-1)^{m-1}\gamma_{m-1}(u)$$

The Hasse identity is

$$(s-1)\varsigma(s,u) = \sum_{n=0}^{\infty}\frac{1}{n+1}\sum_{k=0}^{n}\binom{n}{k}\frac{(-1)^k}{(u+k)^{s-1}}$$

and therefore we have an alternative derivative expression

$$\lim_{s\to 1}\frac{\partial^m}{\partial s^m}[(s-1)\varsigma(s,u)] = (-1)^m \sum_{n=0}^{\infty}\frac{1}{n+1}\sum_{k=0}^{n}\binom{n}{k}(-1)^k \log^m(u+k)$$

This then gives us for $m \geq 0$

(4.3.214) $$\gamma_m(u) = -\frac{1}{m+1}\sum_{n=0}^{\infty}\frac{1}{n+1}\sum_{k=0}^{n}\binom{n}{k}(-1)^k \log^{m+1}(u+k)$$

In particular we have

(4.3.215) $$\gamma_0(u) = -\sum_{n=0}^{\infty}\frac{1}{n+1}\sum_{k=0}^{n}\binom{n}{k}(-1)^k \log(u+k) = -\psi(u)$$

The identity (4.3.215) was noted by Coffey in [45e] (and this concurs with (4.4.99aix) of Volume III where the summation started at $n=1$)

$$\sum_{n=1}^{\infty} \frac{1}{n+1} \sum_{k=0}^{n} \binom{n}{k} (-1)^k \log(u+k) = \psi(u) - \log u$$

$\square$

We have the following summation

$$\sum_{p=0}^{\infty} \frac{(-1)^p}{p!} \gamma_p(u) = \sum_{p=0}^{\infty} \frac{(-1)^p}{p!} \left\{ -\frac{1}{p+1} \sum_{n=0}^{\infty} \frac{1}{n+1} \sum_{k=0}^{n} \binom{n}{k} (-1)^k \log^{p+1}(u+k) \right\}$$

$$= \sum_{n=0}^{\infty} \frac{1}{n+1} \sum_{k=0}^{n} \binom{n}{k} (-1)^k \sum_{p=0}^{\infty} \frac{(-1)^{p+1} \log^{p+1}(u+k)}{(p+1)!}$$

$$= \sum_{n=0}^{\infty} \frac{1}{n+1} \sum_{k=0}^{n} \binom{n}{k} (-1)^k \left( \exp[-\log(u+k)] - 1 \right)$$

$$= \sum_{n=0}^{\infty} \frac{1}{n+1} \sum_{k=0}^{n} \binom{n}{k} (-1)^k \left( \frac{1}{u+k} - 1 \right)$$

$$= -1 + \sum_{n=0}^{\infty} \frac{1}{n+1} \sum_{k=0}^{n} \binom{n}{k} \frac{(-1)^k}{u+k}$$

Therefore we obtain

(4.3.216) $$\sum_{p=0}^{\infty} \frac{(-1)^p}{p!} \gamma_p(u) = \psi'(u) - 1$$

This equation was also derived by Coffey [45c] in 2005 using a different method. The result is in fact most easily seen by letting $s=2$ in the Laurent expansion for $\varsigma(s,u)$

$$\varsigma(s,u) = \frac{1}{s-1} + \sum_{p=0}^{\infty} \frac{(-1)^p}{p!} \gamma_p(u)(s-1)^p$$

This gives us

$$\varsigma(2,u) = 1 + \sum_{p=0}^{\infty} \frac{(-1)^p}{p!} \gamma_p(u)$$

For $s > 1$ we may use the following form of the Hurwitz zeta function

$$\varsigma(2,u) = \sum_{n=0}^{\infty} \frac{1}{(n+u)^2} = \psi'(u)$$

and the result follows immediately.

Using the same method as above, it is easily seen that

$$\sum_{p=0}^{\infty} \frac{\gamma_p(u)}{p!} = \frac{3}{2} - u$$

and this may also be confirmed by letting $s = 0$ in the Laurent expansion for $\varsigma(s,u)$. We then see that

(4.3.217) $$\sum_{p=0}^{\infty} \frac{\gamma_p}{p!} = \frac{1}{2}$$

We may slightly generalise (4.3.216) as follows by including a factor of $t^p$.

$$\sum_{p=0}^{\infty} \frac{(-1)^p}{p!} t^p \gamma_p(u) = \sum_{p=0}^{\infty} \frac{(-1)^p}{p!} t^p \left\{ -\frac{1}{p+1} \sum_{n=0}^{\infty} \frac{1}{n+1} \sum_{k=0}^{n} \binom{n}{k} (-1)^k \log^{p+1}(u+k) \right\}$$

$$= \frac{1}{t} \sum_{n=0}^{\infty} \frac{1}{n+1} \sum_{k=0}^{n} \binom{n}{k} (-1)^k \sum_{p=0}^{\infty} \frac{(-1)^{p+1} t^{p+1} \log^{p+1}(u+k)}{(p+1)!}$$

$$= \frac{1}{t} \sum_{n=0}^{\infty} \frac{1}{n+1} \sum_{k=0}^{n} \binom{n}{k} (-1)^k \left( \exp[-t \log(u+k)] - 1 \right)$$

$$= \frac{1}{t} \sum_{n=0}^{\infty} \frac{1}{n+1} \sum_{k=0}^{n} \binom{n}{k} (-1)^k \left( \frac{1}{[u+k]^t} - 1 \right)$$

$$= -\frac{1}{t} + \frac{1}{t} \sum_{n=0}^{\infty} \frac{1}{n+1} \sum_{k=0}^{n} \binom{n}{k} \frac{(-1)^k}{(u+k)^t}$$

Therefore we obtain (which is simply a defining expression for the Stieltjes constants)

(4.3.217a) $$\sum_{p=0}^{\infty} \frac{(-1)^p}{p!} t^p \gamma_p(u) = \varsigma(t+1,u) - \frac{1}{t}$$

and with $t = 1$ this reverts to (4.3.216)

$$\sum_{p=0}^{\infty}\frac{(-1)^p}{p!}\gamma_p(u)=\varsigma(2,u)-1=\psi'(u)-1$$

## A POSSIBLE CONNECTION WITH THE FRESNEL INTEGRAL

In 1955 Briggs [35b] showed that

(4.3.217b) $$\gamma_p=2\sum_{k=1}^{\infty}\int_1^{\infty}\cos 2\pi kx\,\frac{\log^p x}{x}dx$$

and we have the summation

$$\sum_{p=0}^{\infty}\frac{(-1)^p}{p!}t^p\gamma_p=2\sum_{p=0}^{\infty}\frac{(-1)^p}{p!}t^p\sum_{k=1}^{\infty}\int_1^{\infty}\cos 2\pi kx\,\frac{\log^p x}{x}dx$$

$$=2\sum_{k=1}^{\infty}\int_1^{\infty}\frac{\cos 2\pi kx}{x}\sum_{p=0}^{\infty}\frac{(-1)^p t^p \log^p x}{p!}dx$$

$$=2\sum_{k=1}^{\infty}\int_1^{\infty}\frac{\cos 2\pi kx\exp[-t\log x]}{x}dx$$

$$=2\sum_{k=1}^{\infty}\int_1^{\infty}\frac{\cos 2\pi kx}{x^{t+1}}dx$$

Then using (4.3.217a) we see that

$$\varsigma(t+1)-\frac{1}{t}=2\sum_{k=1}^{\infty}\int_1^{\infty}\frac{\cos 2\pi kx}{x^{t+1}}dx$$

$$=2\sum_{k=1}^{\infty}\int_0^{\infty}\frac{\cos 2\pi kx}{x^{t+1}}dx-2\sum_{k=1}^{\infty}\int_0^1\frac{\cos 2\pi kx}{x^{t+1}}dx$$

$$=2(2\pi)^t\sum_{k=1}^{\infty}\int_0^{\infty}\frac{\cos ku}{u^{t+1}}du-2\sum_{k=1}^{\infty}\int_0^1\frac{\cos 2\pi kx}{x^{t+1}}dx$$

In a manner similar to (F.15) in Volume VI it may be shown that for $0<q<1$

$$\int_0^\infty \frac{\cos u}{u^q}\,du = \frac{\pi}{\Gamma(q)\cos(q\pi/2)}$$

and hence we have for $-1 < t < 0$

$$\int_0^\infty \frac{\cos u}{u^{t+1}}\,du = -\frac{\pi}{\Gamma(t+1)\sin(t\pi/2)}$$

$$\int_0^\infty \frac{\cos ky}{y^{t+1}}\,dy = -\frac{\pi k^t}{\Gamma(t+1)\sin(t\pi/2)}$$

This then gives us

$$\varsigma(t+1) - \frac{1}{t} = -\frac{(2\pi)^{t+1}}{\Gamma(t+1)\sin(t\pi/2)}\sum_{k=1}^\infty k^t - 2\sum_{k=1}^\infty \int_0^1 \frac{\cos 2\pi kx}{x^{t+1}}\,dx$$

Designating $t = -s$ we obtain for $0 < s < 1$

$$\varsigma(1-s) + \frac{1}{s} = \frac{(2\pi)^{1-s}\varsigma(s)}{\Gamma(1-s)\sin(s\pi/2)} - 2\sum_{k=1}^\infty \int_0^1 \frac{\cos 2\pi kx}{x^{1-s}}\,dx$$

and using the Riemann functional equation this becomes

(4.3.217c) $$\varsigma(1-s) - \frac{1}{s} = 2\sum_{k=1}^\infty \int_0^1 \frac{\cos 2\pi kx}{x^{1-s}}\,dx$$

We note that

$$\lim_{s\to 0}\left[\varsigma(1-s) - \frac{1}{s}\right] = 2\sum_{k=1}^\infty \int_0^1 \frac{\cos 2\pi kx}{x}\,dx$$

and hence we have

$$\gamma = 2\sum_{k=1}^\infty \int_0^1 \frac{\cos 2\pi kx}{x}\,dx$$

in accordance with (4.3.217b).

Having regard to (4.3.217b) it is interesting to note that in their paper, An integral representation of the generalized Euler-Mascheroni constants, Ainsworth and Howell [6y] used a contour integral of the form $\int_C \cot \pi z \frac{\log^n}{z} dz$ in order to evaluate the Stieltjes constants: this is reminiscent of the integral identity (6.5a)

$$\frac{1}{2}\int_a^b p(x)\cot(\alpha x/2)\,dx = \sum_{n=0}^{\infty}\int_a^b p(x)\sin\alpha nx\,dx$$

With $s = 1/2$ we get

$$\varsigma\left(\frac{1}{2}\right) + 2 = 2\varsigma\left(\frac{1}{2}\right) - 2\sum_{k=1}^{\infty}\int_0^1 \frac{\cos 2\pi kx}{\sqrt{x}}\,dx$$

and we see with a simple substitution that

$$\int_0^1 \frac{\cos 2\pi kx}{\sqrt{x}}\,dx = \frac{\sqrt{2\pi}}{\pi\sqrt{k}}\int_0^{\sqrt{2\pi k}} \cos(t^2)\,dt$$

The Fresnel integral is defined by G&R [74, p.880]

$$C(x) = \frac{2}{\sqrt{2\pi}}\int_0^x \cos(t^2)\,dt$$

and we then see that

$$\int_0^1 \frac{\cos 2\pi kx}{\sqrt{x}}\,dx = \frac{C(\sqrt{2\pi k})}{\sqrt{k}}$$

Hence we obtain

(4.3.217d) $\quad \frac{3}{2}\varsigma\left(\frac{1}{2}\right) - 2 = 2\sum_{k=1}^{\infty}\frac{C(\sqrt{2\pi k})}{\sqrt{k}}$

It therefore appears that a double series could be established for $\varsigma\left(\frac{1}{2}\right)$ using the known expansion for the Fresnel integral

$$C(x) = \sum_{k=0}^{\infty}\frac{(\pi/2)^{2k}(-1)^k x^{4k+1}}{(4k+1)(2k)!}$$

Letting $p = 0$ in (4.3.217b) gives us

$$\gamma = \gamma_0 = 2\sum_{k=1}^{\infty} \int_1^{\infty} \frac{\cos 2\pi kx}{x} dx = 2\sum_{k=1}^{\infty} [Ci(\infty) - Ci(2\pi k)]$$

$$= -2\sum_{k=1}^{\infty} Ci(2\pi k)$$

where $Ci(x)$ is the cosine integral defined in G&R [74, p.878] and [1, p.231] as

$$Ci(x) = \gamma + \log x + \int_0^x \frac{\cos t - 1}{t} dt = \gamma + \log x + \sum_{n=1}^{\infty} \frac{(-1)^n x^{2n}}{2n(2n)!}$$

(see also (6.94g) et seq. in Volume V for a more detailed exposition) and we can see from (6.94ga) that $Ci(\infty) = 0$.

We will see in (6.117r) in Volume V that

$$\gamma = \frac{1}{2} - 2\sum_{k=1}^{\infty} Ci(2k\pi)$$

and we therefore have an unexplained difference of $1/2$ which remains to be investigated. Using the Wolfram Integrator, we find that integrals of the form $\int \cos 2\pi kx \frac{\log^p x}{x} dx$ may be evaluated in terms of generalised hypergeometric series.

We briefly reconsider the equation

$$\sum_{p=0}^{\infty} \frac{(-1)^p}{p!} t^p \gamma_p = 2\sum_{p=0}^{\infty} \frac{(-1)^p}{p!} t^p \sum_{k=1}^{\infty} \int_1^{\infty} \cos 2\pi kx \frac{\log^p x}{x} dx$$

$$= 2\sum_{p=0}^{\infty} \frac{(-1)^p}{p!} t^p \sum_{k=1}^{\infty} \int_0^{\infty} \cos 2\pi kx \frac{\log^p x}{x} dx - 2\sum_{p=0}^{\infty} \frac{(-1)^p}{p!} t^p \sum_{k=1}^{\infty} \int_0^1 \cos 2\pi kx \frac{\log^p x}{x} dx$$

We now have regard to the second summation. Dilcher [54a] proved by induction that if $f_p(x) = \frac{\log^p x}{x}$ then

$$f_p^{(n)}(x) = \frac{p!}{x^{n+1}} \sum_{j=0}^{p} s(n+1, p+1-j) \frac{\log^j x}{j!}$$

where $s(n+1, p+1)$ is a Stirling number of the first kind. We then have

$$f_p^{(n)}(1) = p!\,s(n+1, p+1)$$

which results in the Taylor series expansion

$$f_p(x) = \frac{\log^p x}{x} = \sum_{n=0}^{\infty} \frac{f_p^{(n)}(1)}{n!}(x-1)^n = p!\sum_{n=0}^{\infty} \frac{s(n+1, p+1)}{n!}(x-1)^n$$

If $x \to 1+z$ this may be written as

$$\log^p(1+z) = p!(1+z)\sum_{n=0}^{\infty} \frac{s(n+1, p+1)}{n!} z^n = p!\sum_{n=0}^{\infty} \frac{s(n+1, p+1)}{n!} z^n + p!\sum_{n=0}^{\infty} \frac{s(n+1, p+1)}{n!} z^{n+1}$$

We have [126, p.56] $s(n+1, p+1) = s(n, p) - ns(n, p+1)$ and thus

$$\log^p(1+z) = p!\sum_{n=0}^{\infty} \frac{s(n, p) - ns(n, p+1)}{n!} z^n + p!\sum_{n=0}^{\infty} \frac{s(n+1, p+1)}{n!} z^{n+1}$$

$$= p!\sum_{n=0}^{\infty} \frac{s(n, p)}{n!} z^n - p!\sum_{n=0}^{\infty} \frac{ns(n, p+1)}{n!} z^n + p!\sum_{n=0}^{\infty} \frac{s(n+1, p+1)}{n!} z^{n+1}$$

$$= p!\sum_{n=0}^{\infty} \frac{s(n, p)}{n!} z^n - p!\sum_{n=0}^{\infty} \frac{s(n, p+1)}{(n-1)!} z^n + p!\sum_{n=0}^{\infty} \frac{s(n+1, p+1)}{n!} z^{n+1}$$

$$= p!\sum_{n=0}^{\infty} \frac{s(n, p)}{n!} z^n - p!\sum_{m=1}^{\infty} \frac{s(m+1, p+1)}{m!} z^{m+1} + p!\sum_{n=0}^{\infty} \frac{s(n+1, p+1)}{n!} z^{n+1}$$

and we therefore have obtained the expansion (3.105) from Volume I for the Stirling number of the first kind which is valid for $|z| < 1$

$$\log^p(1+z) = p!\sum_{n=0}^{\infty} \frac{s(n, p)}{n!} z^n \quad \text{or} \quad \log^p x = p!\sum_{n=0}^{\infty} \frac{s(n, p)}{n!}(x-1)^n$$

We then have

$$\int_0^1 \cos 2\pi kx \frac{\log^p x}{x} dx = p!\int_0^1 \frac{\cos 2\pi kx}{x} \sum_{n=0}^{\infty} \frac{s(n, p)}{n!}(x-1)^n dx$$

or alternatively without $x$ in the denominator

$$\int_0^1 \cos 2\pi kx \frac{\log^p x}{x} dx = p! \int_0^1 \cos 2\pi kx \sum_{n=0}^{\infty} \frac{s(n+1, p+1)}{n!} (x-1)^n dx$$

We have

$$\int_0^1 (x-1)^n \cos 2\pi kx \, dx = (-1)^n \int_0^1 (1-x)^n \cos 2\pi kx \, dx$$

and with the substitution $y = 1 - x$ this becomes

$$= (-1)^n \int_0^1 y^n \cos 2\pi ky \, dy$$

Unfortunately, any further analysis in this direction will have to wait for another day.

## EVALUATION OF SOME STIELTJES CONSTANTS

Let us consider the function $Z(s,u)$ defined by

$$Z(s,u) = \sum_{n=0}^{\infty} \frac{1}{n+1} \sum_{k=0}^{n} \binom{n}{k} (-1)^k (k+u)^s$$

$$= \sum_{n=0}^{\infty} \frac{1}{n+1} \sum_{k=0}^{n} \binom{n}{k} (-1)^k \exp[s \log(k+u)]$$

$$= \sum_{n=0}^{\infty} \frac{1}{n+1} \sum_{k=0}^{n} \binom{n}{k} (-1)^k \sum_{p=0}^{\infty} \frac{s^p \log^p(k+u)}{p!}$$

$$= \sum_{p=0}^{\infty} \frac{s^p}{p!} \sum_{n=0}^{\infty} \frac{1}{n+1} \sum_{k=0}^{n} \binom{n}{k} (-1)^k \log^p(k+u)$$

and using

$$\gamma_p(u) = -\frac{1}{p+1} \sum_{n=0}^{\infty} \frac{1}{n+1} \sum_{k=0}^{n} \binom{n}{k} (-1)^k \log^{p+1}(u+k)$$

this becomes

$$= 1 - \sum_{p=0}^{\infty} \frac{s^p}{p!} \gamma_p(u)$$

Then, using (4.3.217a), we have

$$Z(s,u) = \sum_{n=0}^{\infty} \frac{1}{n+1} \sum_{k=0}^{n} \binom{n}{k} (-1)^k (k+u)^s = -s\varsigma(1-s,u)$$

The fact that $-s\varsigma(1-s,u) = \sum_{n=0}^{\infty} \frac{1}{n+1} \sum_{k=0}^{n} \binom{n}{k} (-1)^k (k+u)^s$ may in fact be immediately be seen from the Hasse formula (3.12a) from Volume I.

Differentiating $Z(s,u)$ gives us

$$\frac{\partial}{\partial s} Z(s,u) = \sum_{n=0}^{\infty} \frac{1}{n+1} \sum_{k=0}^{n} \binom{n}{k} (-1)^k (k+u)^s \log(k+u)$$

$$= s\varsigma'(1-s,u) - \varsigma(1-s,u)$$

Letting $s = 2m+1$ we obtain another derivation of (4.3.133b)

(4.3.217e) $\quad \varsigma(2m+1) = (-1)^m \dfrac{2(2\pi)^{2m}}{(2m+1)!} \sum_{n=0}^{\infty} \dfrac{1}{n+1} \sum_{k=0}^{n} \binom{n}{k} (-1)^k (k+1)^{2m+1} \log(k+1)$

and with $s = 2m$ we get

(4.3.217f) $\quad 2m\varsigma'(1-2m,u) - \varsigma(1-2m,u) = \sum_{n=0}^{\infty} \dfrac{1}{n+1} \sum_{k=0}^{n} \binom{n}{k} (-1)^k (k+1)^{2m} \log(k+1)$

Also, using Lerch's identity (4.3.116), with $s = 1$ we obtain

$$\log \Gamma(u) - \frac{1}{2}\log(2\pi) - \varsigma(0,u) = \sum_{n=0}^{\infty} \frac{1}{n+1} \sum_{k=0}^{n} \binom{n}{k} (-1)^k (k+u) \log(k+u)$$

which we have seen before in (4.3.119). Note that we previously we used (4.3.119) to derive Lerch's identity.

With regard to the definition of $Z(s,u)$ we may also note from (A.25) that

$$B_N(u) = \sum_{n=0}^{N} \frac{1}{n+1} \sum_{k=0}^{n} \binom{n}{k} (-1)^k (k+u)^N$$

$$= \sum_{n=0}^{\infty} \frac{1}{n+1} \sum_{k=0}^{n} \binom{n}{k} (-1)^k (k+u)^N$$

It may also be useful to consider a function of the type

$$Z(r,s,u,v) = \sum_{n=0}^{\infty} \frac{1}{n+1} \sum_{k=0}^{n} \binom{n}{k} (-1)^k (k+u)^r (k+v)^s$$

We also note for reference that Kluyver showed in 1922 that

(4.3.217g) $$\gamma_p = \frac{\log^p 2}{p+1} \sum_{k=1}^{\infty} \frac{(-1)^k}{k} B_{p+1}\left(\frac{\log k}{\log 2}\right)$$

$\square$

We may also consider a more complex series

$$\sum_{n=1}^{\infty} z^n \sum_{k=0}^{\infty} \frac{(-1)^k}{k!} n^k \gamma_k(u) = \sum_{n=1}^{\infty} z^n \sum_{k=0}^{\infty} \frac{(-1)^k}{k!} \left\{ -\frac{1}{k+1} \sum_{m=0}^{\infty} \frac{1}{m+1} \sum_{r=0}^{m} \binom{m}{r} (-1)^r n^k \log^{k+1}(u+r) \right\}$$

$$= \sum_{n=1}^{\infty} \frac{z^n}{n} \sum_{m=0}^{\infty} \frac{1}{m+1} \sum_{r=0}^{m} \binom{m}{r} (-1)^r \sum_{k=0}^{\infty} \frac{(-1)^k n^{k+1} \log^{k+1}(u+r)}{(k+1)!}$$

$$= \sum_{n=1}^{\infty} \frac{z^n}{n} \sum_{m=0}^{\infty} \frac{1}{m+1} \sum_{r=0}^{m} \binom{m}{r} (-1)^r \exp\{[-n\log(u+r)] - 1\}$$

$$= \sum_{n=1}^{\infty} \frac{z^n}{n} \sum_{m=0}^{\infty} \frac{1}{m+1} \sum_{r=0}^{m} \binom{m}{r} (-1)^r \left( \frac{1}{(u+r)^n} - 1 \right)$$

The Hasse identity tells us that

$$\sum_{m=0}^{\infty} \frac{1}{m+1} \sum_{r=0}^{m} \binom{m}{r} \frac{(-1)^r}{(u+r)^n} = n\varsigma(n+1,u)$$

and we therefore obtain

(4.3.217h) $$\sum_{n=1}^{\infty} z^n \sum_{k=0}^{\infty} \frac{(-1)^k}{k!} n^k \gamma_k(u) = \sum_{n=1}^{\infty} \varsigma(n+1,u) z^n - \sum_{n=1}^{\infty} \frac{z^n}{n}$$

With $u = 1$ we get

$$\sum_{n=1}^{\infty} z^n \sum_{k=0}^{\infty} \frac{(-1)^k}{k!} n^k \gamma_k = \sum_{n=1}^{\infty} \varsigma(n+1) z^n - \sum_{n=1}^{\infty} \frac{z^n}{n}$$

which may be written for $|z|<1$ as

(4.3.218) $$\sum_{n=1}^{\infty} z^n \sum_{k=0}^{\infty} \frac{(-1)^k}{k!} n^k \gamma_k = \log(1-z) - \gamma - \psi(1-z)$$

This result was obtained in a different way by Lopez [101ab] in 1998.

Coffey [45e] employed a different method in 2006 to show that for integers $j \geq 0$ and $|z|<1$

$$\sum_{n=1}^{\infty} z^n \sum_{k=0}^{\infty} \frac{(-1)^k}{(k-j)!} n^{k-j} \gamma_k(u) = \sum_{n=1}^{\infty} \varsigma^{(j)}(n+1,u) z^n - (-1)^j j! Li_{j+1}(z)$$

and with $j=0$ this becomes

$$\sum_{n=1}^{\infty} z^n \sum_{k=0}^{\infty} \frac{(-1)^k}{k!} n^k \gamma_k(u) = \sum_{n=1}^{\infty} \varsigma(n+1,u) z^n - \log(1-z)$$

$$= \psi(u) - \psi(u-z) - \log(1-z)$$

Letting $z \to -z$ we obtain

(4.3.219) $$\sum_{n=1}^{\infty} (-1)^n z^n \sum_{k=0}^{\infty} \frac{(-1)^k}{k!} n^k \gamma_k(u) = -\sum_{n=2}^{\infty} (-1)^n \varsigma(n,u) z^{n-1} + \log(1+z)$$

and integrating this with respect to $z$ produces

$$\sum_{n=1}^{\infty} \frac{(-1)^n}{n+1} x^{n+1} \sum_{k=0}^{\infty} \frac{(-1)^k}{k!} n^k \gamma_k(u) = -\sum_{n=2}^{\infty} (-1)^n \varsigma(n,u) \frac{x^n}{n} + (1+x)\log(1+x) - x$$

We have the well-known result [126, p.159] for $|x| < |u|$

(4.3.220) $$\sum_{n=2}^{\infty} (-1)^n \varsigma(n,u) \frac{x^n}{n} = \log \Gamma(u+x) - \log \Gamma(u) - x\psi(u)$$

and hence we have

$$\sum_{n=1}^{\infty}\frac{(-1)^{n}}{n+1}x^{n+1}\sum_{k=0}^{\infty}\frac{(-1)^{k}}{k!}n^{k}\gamma_{k}(u)=-\log\Gamma(u+x)+\log\Gamma(u)+x\psi(u)+(1+x)\log(1+x)-x$$

Differentiating (4.3.220) results in

$$\sum_{n=2}^{\infty}(-1)^{n}\varsigma(n,u)x^{n-1}=\psi(u+x)-\psi(u)$$

and we obtain another proof of Coffey's formula (4.3.219) above

$$\sum_{n=1}^{\infty}(-1)^{n}z^{n}\sum_{k=0}^{\infty}\frac{(-1)^{k}}{k!}n^{k}\gamma_{k}(u)=-\psi(u+z)+\psi(u)+\log(1+z)$$

Multiplying (4.3.219) by $z$ and then integrating gives us

$$\sum_{n=1}^{\infty}(-1)^{n}\frac{z^{n+2}}{n+2}\sum_{k=0}^{\infty}\frac{(-1)^{k}}{k!}n^{k}\gamma_{k}(u)=-\sum_{n=2}^{\infty}(-1)^{n}\varsigma(n,u)\frac{z^{n+1}}{n+1}+\frac{1}{4}[2(z^{2}-1)\log(1+z)-z^{2}+2z]$$

and we know from [126, p.211] that for $|z|<|u|$ (see also Coffey [45j, p.13])

$$\sum_{n=2}^{\infty}(-1)^{n}\varsigma(n,u)\frac{z^{n+1}}{n+1}=[2u-1-\log(2\pi)]\frac{z}{2}+[1-\psi(u)]\frac{z^{2}}{2}+(1-u)\log\Gamma(u+z)$$

$$+\log G(u+z)+(u-1)\log\Gamma(u)-\log G(u)$$

Hence we obtain

(4.3.221)

$$\sum_{n=1}^{\infty}(-1)^{n}\frac{z^{n+2}}{n+2}\sum_{k=0}^{\infty}\frac{(-1)^{k}}{k!}n^{k}\gamma_{k}(u)=-[2u-1-\log(2\pi)]\frac{z}{2}-[1-\psi(u)]\frac{z^{2}}{2}-(1-u)\log\Gamma(u+z)$$

$$-\log G(u+z)-(u-1)\log\Gamma(u)+\log G(u)+\frac{1}{4}[2(z^{2}-1)\log(1+z)-z^{2}+2z]$$

$\square$

We may also consider the following summation

$$\sum_{p=0}^{\infty}\frac{(-1)^{p}}{p!}(s-1)^{p}\gamma_{p+r}(u)=\sum_{p=0}^{\infty}\frac{(-1)^{p}(s-1)^{p}}{p!}\left\{-\frac{1}{p+r+1}\sum_{n=0}^{\infty}\frac{1}{n+1}\sum_{k=0}^{n}\binom{n}{k}(-1)^{k}\log^{p+r+1}(u+k)\right\}$$

$$= \frac{(-1)^r}{(s-1)^{r+1}} \sum_{n=0}^{\infty} \frac{1}{n+1} \sum_{k=0}^{n} \binom{n}{k}(-1)^k \sum_{p=0}^{\infty}(p+r)(p+r-1)...(p+1)\frac{(-1)^{p+r+1}(s-1)^{p+r+1}\log^{p+r+1}(u+k)}{(p+r+1)!}$$

With $r=1$ we get

$$\sum_{p=0}^{\infty}\frac{(-1)^p}{p!}(s-1)^p \gamma_{p+1}(u) = \sum_{p=0}^{\infty}\frac{(-1)^p(s-1)^p}{p!}\left\{-\frac{1}{p+2}\sum_{n=0}^{\infty}\frac{1}{n+1}\sum_{k=0}^{n}\binom{n}{k}(-1)^k \log^{p+2}(u+k)\right\}$$

$$= -\frac{1}{(s-1)^2}\sum_{n=0}^{\infty}\frac{1}{n+1}\sum_{k=0}^{n}\binom{n}{k}(-1)^k \sum_{p=0}^{\infty}(p+1)\frac{(-1)^{p+2}(s-1)^{p+2}\log^{p+2}(u+k)}{(p+2)!}$$

We see that

$$\sum_{p=0}^{\infty}\frac{x^{p+2}}{(p+2)!} = e^x - 1 - x$$

and also that

$$x^2 \frac{d}{dx}\left(\frac{e^x-1-x}{x}\right) = \sum_{p=0}^{\infty}(p+1)\frac{x^{p+2}}{(p+2)!}$$

We therefore have

$$1+(x-1)e^x = \sum_{p=0}^{\infty}(p+1)\frac{x^{p+2}}{(p+2)!}$$

and hence we get

$$\sum_{p=0}^{\infty}\frac{(-1)^p}{p!}(s-1)^p \gamma_{p+1}(u) = -\frac{1}{(s-1)^2}\sum_{n=0}^{\infty}\frac{1}{n+1}\sum_{k=0}^{n}\binom{n}{k}(-1)^k\left[1-(s-1)\frac{\log(u+k)}{(u+k)^{s-1}} - \frac{1}{(u+k)^{s-1}}\right]$$

$$= -\frac{1}{(s-1)^2} + \frac{1}{s-1}\sum_{n=0}^{\infty}\frac{1}{n+1}\sum_{k=0}^{n}\binom{n}{k}(-1)^k \frac{\log(u+k)}{(u+k)^{s-1}} + \frac{1}{(s-1)^2}\sum_{n=0}^{\infty}\frac{1}{n+1}\sum_{k=0}^{n}\binom{n}{k}(-1)^k \frac{1}{(u+k)^{s-1}}$$

$$= -\frac{1}{(s-1)^2} - \varsigma'(s,u)$$

Therefore we obtain

$$(4.3.222) \qquad -\varsigma'(s,u) = \frac{1}{(s-1)^2} + \sum_{p=0}^{\infty} \frac{(-1)^p}{p!}(s-1)^p \gamma_{p+1}(u)$$

This then implies that $\lim_{s \to 1}[(s-1)^2 \varsigma'(s,u)] = -1$ which concurs with (4.3.203).

More generally, with a little more algebra, we may show that

$$(4.3.222a) \qquad (-1)^r \varsigma^{(r)}(s,u) = \frac{r!}{(s-1)^{r+1}} + \sum_{p=0}^{\infty} \frac{(-1)^p}{p!}(s-1)^p \gamma_{p+r}(u)$$

In 1995, Choudhury [45ad] mentioned that, for $r \geq 1$ and $\mathrm{Re}(s) > 1$, Ramanujan determined that

$$(4.3.222b) \qquad (-1)^r \varsigma^{(r)}(s) = \sum_{p=1}^{\infty} \frac{\log^r p}{p^s} = \frac{r!}{(s-1)^{r+1}} + \sum_{p=0}^{\infty} \frac{(-1)^p}{p!}(s-1)^p \gamma_{p+r}$$

(see also Ramanujan's Notebooks [21], Part I, p.224). This is a particular case of (4.3.222a) where $u = 1$, and it should be noted that (4.3.222a) is valid for all $s \neq 1$.

Letting $s = 0$ in (4.3.222) we see that

$$-\varsigma'(0,u) = 1 + \sum_{p=0}^{\infty} \frac{\gamma_{p+1}(u)}{p!}$$

and applying Lerch's identity (4.3.116) this may be written as

$$(4.3.223) \qquad \frac{1}{2}\log(2\pi) - \log \Gamma(u) = 1 + \sum_{p=0}^{\infty} \frac{\gamma_{p+1}(u)}{p!}$$

which may be compared with (4.3.216)

$$\psi'(u) = 1 + \sum_{p=0}^{\infty} \frac{(-1)^p}{p!} \gamma_p(u)$$

Differentiating (4.3.223) we have

$$\psi(u) = -\sum_{p=0}^{\infty} \frac{\gamma'_{p+1}(u)}{p!}$$

$$= \sum_{p=0}^{\infty} \frac{1}{p!} \left\{ \sum_{n=0}^{\infty} \frac{1}{n+1} \sum_{k=0}^{n} \binom{n}{k}(-1)^k \frac{\log^{p+1}(u+k)}{u+k} \right\}$$

$$= \sum_{p=0}^{\infty} \frac{\log^p(u+k)}{p!} \left\{ \sum_{n=0}^{\infty} \frac{1}{n+1} \sum_{k=0}^{n} \binom{n}{k} (-1)^k \frac{\log(u+k)}{u+k} \right\}$$

$$= \sum_{n=0}^{\infty} \frac{1}{n+1} \sum_{k=0}^{n} \binom{n}{k} (-1)^k \log(u+k)$$

which we saw earlier in (4.3.74).

A further differentiation of (4.3.223) gives us

$$\psi'(u) = -\sum_{p=0}^{\infty} \frac{\gamma''_{p+1}(u)}{p!}$$

$$= \sum_{p=0}^{\infty} \frac{1}{p!} \left\{ (p+1) \sum_{n=0}^{\infty} \frac{1}{n+1} \sum_{k=0}^{n} \binom{n}{k} (-1)^k \frac{\log^p(u+k)}{(u+k)^2} - \sum_{n=0}^{\infty} \frac{1}{n+1} \sum_{k=0}^{n} \binom{n}{k} (-1)^k \frac{\log^{p+1}(u+k)}{(u+k)^2} \right\}$$

$$= \sum_{n=0}^{\infty} \frac{1}{n+1} \sum_{k=0}^{n} \binom{n}{k} (-1)^k \frac{1}{(u+k)^2} \sum_{p=0}^{\infty} (p+1) \frac{\log^p(u+k)}{p!}$$

$$- \sum_{n=0}^{\infty} \frac{1}{n+1} \sum_{k=0}^{n} \binom{n}{k} (-1)^k \frac{\log(u+k)}{(u+k)^2} \sum_{p=0}^{\infty} \frac{\log^p(u+k)}{p!}$$

$$= \sum_{n=0}^{\infty} \frac{1}{n+1} \sum_{k=0}^{n} \binom{n}{k} (-1)^k \frac{1}{(u+k)^2} [1+\log(u+k)](u+k)$$

$$- \sum_{n=0}^{\infty} \frac{1}{n+1} \sum_{k=0}^{n} \binom{n}{k} (-1)^k \frac{\log(u+k)}{(u+k)^2} (u+k)$$

$$= \sum_{n=0}^{\infty} \frac{1}{n+1} \sum_{k=0}^{n} \binom{n}{k} (-1)^k \frac{1}{u+k}$$

$= \psi'(u)$ (which is where we started from).

Letting $s = -1$ in (4.3.222) gives us

$$-\varsigma'(-1,u) = \frac{1}{4} + \sum_{p=0}^{\infty} \frac{2^p}{p!} \gamma_{p+1}(u)$$

Differentiation gives us

$$-\frac{d}{du}\varsigma'(-1,u) = \sum_{p=0}^{\infty}\frac{2^p}{p!}\gamma'_{p+1}(u)$$

and using

$$\gamma'_{p+1}(u) = -\sum_{n=0}^{\infty}\frac{1}{n+1}\sum_{k=0}^{n}\binom{n}{k}(-1)^k\frac{\log^{p+1}(u+k)}{u+k}$$

we have

$$\frac{d}{du}\varsigma'(-1,u) = \sum_{p=0}^{\infty}\frac{2^p}{p!}\sum_{n=0}^{\infty}\frac{1}{n+1}\sum_{k=0}^{n}\binom{n}{k}(-1)^k\frac{\log^{p+1}(u+k)}{u+k}$$

$$= \sum_{p=0}^{\infty}\frac{2^p\log^p(u+k)}{p!}\sum_{n=0}^{\infty}\frac{1}{n+1}\sum_{k=0}^{n}\binom{n}{k}(-1)^k\frac{\log(u+k)}{u+k}$$

$$= \sum_{n=0}^{\infty}\frac{1}{n+1}\sum_{k=0}^{n}\binom{n}{k}(-1)^k\frac{\log(u+k)}{u+k}\sum_{p=0}^{\infty}\frac{2^p\log^p(u+k)}{p!}$$

$$= \sum_{n=0}^{\infty}\frac{1}{n+1}\sum_{k=0}^{n}\binom{n}{k}(-1)^k\frac{\log(u+k)}{u+k}\exp[2\log(u+k)]$$

$$= \sum_{n=0}^{\infty}\frac{1}{n+1}\sum_{k=0}^{n}\binom{n}{k}(-1)^k(u+k)\log(u+k)$$

Then referring to (4.3.108) we see that

$$\frac{d}{du}\varsigma'(-1,u) = \varsigma'(0,u) - \varsigma(0,u)$$

We also see that

$$1 + \sum_{p=0}^{\infty}\frac{(-1)^p}{p!}\gamma_p(u) = -\sum_{p=0}^{\infty}\frac{\gamma''_{p+1}(u)}{p!}$$

We have from (4.3.222a)

$$(-1)^r\varsigma^{(r)}(0,u) = (-1)^{r+1}r! + \sum_{p=0}^{\infty}\frac{\gamma_{p+r}(u)}{p!}$$

and hence with $u = 1$

$$(-1)^r \varsigma^{(r)}(0,1) = (-1)^r \varsigma^{(r)}(0) = (-1)^{r+1} r! + \sum_{p=0}^{\infty} \frac{\gamma_{p+r}}{p!}$$

From (4.3.222b) we see that

$$(-1)^r \left[ \varsigma^{(r)}(s,u) - \frac{(-1)^r r!}{(s-1)^{r+1}} \right] = \sum_{p=0}^{\infty} \frac{(-1)^p}{p!} (s-1)^p \gamma_{p+r}(u)$$

and hence as we have seen before

$$(-1)^r \lim_{s \to 1} \left[ \varsigma^{(r)}(s,u) - \frac{(-1)^r r!}{(s-1)^{r+1}} \right] = \gamma_r(u)$$

It seems that a lot of the above analysis was superfluous: having regard to

$$\varsigma(s,u) = \frac{1}{s-1} + \sum_{p=0}^{\infty} \frac{(-1)^p}{p!} \gamma_p(u)(s-1)^p$$

we immediately see that

$$\varsigma^{(r)}(s,u) = \frac{(-1)^r r!}{(s-1)^{r+1}} + \sum_{p=0}^{\infty} \frac{(-1)^p}{p!} p(p-1)\ldots(p-r+1)\gamma_p(u)(s-1)^{p-r}$$

which is in fact exactly equivalent to (4.3.222a)

$$\varsigma^{(r)}(s,u) = \frac{(-1)^r r!}{(s-1)^{r+1}} + (-1)^r \sum_{p=0}^{\infty} \frac{(-1)^p}{p!} (s-1)^p \gamma_{p+r}(u)$$

We have

$$\frac{\partial}{\partial u} \varsigma(s,u) = \sum_{p=0}^{\infty} \frac{(-1)^p}{p!} \gamma'_p(u)(s-1)^p = -s\varsigma(s+1,u)$$

and hence we obtain

$$\sum_{p=0}^{\infty} \frac{(-1)^p}{p!} \gamma'_p(u)(s-1)^p = -1 - \sum_{p=0}^{\infty} \frac{(-1)^p}{p!} \gamma_p(u) s^{p+1}$$

The left-hand side is equal to

$$-\sum_{p=0}^{\infty} \frac{(-1)^p}{p!} \sum_{n=0}^{\infty} \frac{1}{n+1} \sum_{k=0}^{n} \binom{n}{k} (-1)^k \frac{\log^p(u+k)}{u+k} (s-1)^p$$

$$= -\sum_{n=0}^{\infty} \frac{1}{n+1} \sum_{k=0}^{n} \binom{n}{k} \frac{(-1)^k}{u+k} \sum_{p=0}^{\infty} \frac{(-1)^p (s-1)^p \log^p(u+k)}{p!}$$

$$= -\sum_{n=0}^{\infty} \frac{1}{n+1} \sum_{k=0}^{n} \binom{n}{k} \frac{(-1)^k}{u+k} \exp[-(s-1)\log(u+k)]$$

$$= -\sum_{n=0}^{\infty} \frac{1}{n+1} \sum_{k=0}^{n} \binom{n}{k} \frac{(-1)^k}{(u+k)^s}$$

$$= -s\varsigma(s+1, u)$$

$\square$

We see that

$$\frac{d}{du}\gamma_0(u) = -\sum_{n=0}^{\infty} \frac{1}{n+1} \sum_{k=0}^{n} \binom{n}{k} \frac{(-1)^k}{u+k}$$

$$= -\varsigma(2, u) = -\sum_{n=0}^{\infty} \frac{1}{(n+u)^2}$$

where we have employed the Hasse formula (4.3.106a). Integration then results in the familiar formula for the digamma function [126, p.14]

(4.3.223a) $\quad \gamma_0(x) - \gamma_0(1) = \sum_{n=0}^{\infty} \left[\frac{1}{n+x} - \frac{1}{n+1}\right] = -(x-1)\sum_{n=0}^{\infty} \frac{1}{(n+1)(n+x)} = -\psi(x) - \gamma$

Similarly we find that

$$\frac{d}{du}\gamma_1(u) = -\sum_{n=0}^{\infty} \frac{1}{n+1} \sum_{k=0}^{n} \binom{n}{k} \frac{(-1)^k \log(u+k)}{u+k}$$

The Hasse formula (4.3.106a) gives us

$$\frac{\partial}{\partial s}[(s-1)\varsigma(s,u)] = (s-1)\varsigma'(s,u) + \varsigma(s,u) = -\sum_{n=0}^{\infty} \frac{1}{n+1} \sum_{k=0}^{n} \binom{n}{k} \frac{(-1)^k \log(u+k)}{(u+k)^{s-1}}$$

and with $s = 2$ we have

$$\varsigma'(2,u) + \varsigma(2,u) = -\sum_{n=0}^{\infty} \frac{1}{n+1} \sum_{k=0}^{n} \binom{n}{k} \frac{(-1)^k \log(u+k)}{u+k}$$

Therefore we get

(4.3.223b)

$$\frac{d}{du}\gamma_1(u) = \varsigma'(2,u) + \varsigma(2,u) = -\sum_{n=0}^{\infty}\frac{\log(n+u)}{(n+u)^2} + \sum_{n=0}^{\infty}\frac{1}{(n+u)^2} = \sum_{n=0}^{\infty}\frac{1-\log(n+u)}{(n+u)^2}$$

We may easily evaluate the following integral

$$\int_1^x \frac{1-\log(n+u)}{(n+u)^2}du = \frac{\log(n+u)}{n+u}\bigg|_1^x$$

and hence we have upon integrating (see also (4.3.228b))

(4.3.223c) $\qquad \gamma_1(x) - \gamma_1(1) = \sum_{n=0}^{\infty}\left[\frac{\log(n+x)}{n+x} - \frac{\log(n+1)}{n+1}\right]$

As mentioned previously, the Stieltjes constants $\gamma_n$ are the coefficients in the Laurent expansion of $\varsigma(s)$ about $s=1$.

$$\varsigma(s) = \frac{1}{s-1} + \sum_{n=0}^{\infty}\frac{(-1)^n}{n!}\gamma_n(s-1)^n \quad \text{or} \quad \varsigma(s+1) = \frac{1}{s} + \sum_{n=0}^{\infty}\frac{(-1)^n}{n!}\gamma_n s^n$$

Since $\lim_{s\to 1}\left[\varsigma(s) - \frac{1}{s-1}\right] = \gamma$ (see also (4.4.99m)), it is clear that $\gamma_0 = \gamma$. It may be shown, as in [81, p.4], that

(4.3.224) $\qquad \gamma_n = \lim_{N\to\infty}\left[\sum_{k=1}^{N}\frac{\log^n k}{k} - \frac{\log^{n+1} N}{n+1}\right] = \lim_{N\to\infty}\left[\sum_{k=1}^{N}\frac{\log^n k}{k} - \int_1^N \frac{\log^n t}{t}dt\right]$

For example, Bohman and Fröberg [24a] noted that

$$(s-1)\varsigma(s) = \sum_{k=1}^{\infty}\frac{(s-1)}{k^s} \quad \text{and} \quad \sum_{k=1}^{\infty}\left[\frac{1}{k^{s-1}} - \frac{1}{(k+1)^{s-1}}\right] = 1$$

and, assuming that $s$ is real and greater than 1, the above two equations may be subtracted to give

(4.3.225) $\qquad (s-1)\varsigma(s) = 1 + \sum_{k=1}^{\infty}\left[\frac{1}{(k+1)^{s-1}} - \frac{1}{k^{s-1}} + \frac{(s-1)}{k^s}\right]$

It will be noted from the above that

$$\lim_{s \to 1}(s-1)\varsigma(s) = 1$$

which we shall also see in (4.4.100x).

Equation (4.3.225) may be written as

$$(s-1)\varsigma(s) = 1 + \sum_{k=1}^{\infty}\left[\exp(-(s-1)\log(k+1)) - \exp(-(s-1)\log k) + (s-1)k^{-1}\exp(-(s-1)\log k)\right]$$

$$= 1 + \sum_{k=1}^{\infty}\left[\sum_{n=0}^{\infty}\frac{(-1)^n(s-1)^n}{n!}\left(\log^n(k+1) - \log^n k\right) + \frac{s-1}{k}\sum_{n=0}^{\infty}\frac{(-1)^n(s-1)^n}{n!}\log^n k\right]$$

Dividing by $(s-1)$ we get (4.3.224) where

(4.3.225a) $\quad \gamma_n = \sum_{k=1}^{\infty}\left[\frac{\log^n k}{k} - \frac{\log^n(k+1) - \log^n k}{n+1}\right] = \sum_{k=1}^{\infty}\left[\frac{\log^n k}{k} - \int_{k}^{k+1}\frac{\log^n t}{t}dt\right]$

Osler [105(v)] has provided an interesting derivation of the above formula using partial sums of the Riemann zeta function.

Differentiating (4.3.225) we get

$$(s-1)\varsigma'(s) + \varsigma(s) = \sum_{k=1}^{\infty}\left[-\frac{\log(k+1)}{(k+1)^{s-1}} + \frac{\log k}{k^{s-1}} - \frac{(s-1)\log k}{k^s} + \frac{1}{k^s}\right]$$

and hence we have

$$(s-1)\varsigma'(s) = \sum_{k=1}^{\infty}\left[-\frac{\log(k+1)}{(k+1)^{s-1}} + \frac{\log k}{k^{s-1}} - \frac{(s-1)\log k}{k^s}\right]$$

which implies that (as is easily seen by inspection)

$$\sum_{k=1}^{\infty}\frac{\log(k+1)}{(k+1)^{s-1}} = \sum_{k=1}^{\infty}\frac{\log k}{k^{s-1}}$$

Ivić [81, p.41] reports that the coefficients $\gamma_p(x)$ for $x \in (0,1]$ may be expressed as

(4.3.226) $\quad \gamma_p(x) = \lim_{N \to \infty}\left(\sum_{n=0}^{N}\frac{\log^p(n+x)}{n+x} - \frac{\log^{p+1}(N+x)}{p+1}\right)$

and in particular we have

(4.3.226i) $$\gamma_0(x) = \lim_{N \to \infty} \left( \sum_{n=0}^{N} \frac{1}{n+x} - \log(N+x) \right)$$

(4.3.226ii) $$\gamma_1(x) = \lim_{N \to \infty} \left( \sum_{n=0}^{N} \frac{\log(n+x)}{n+x} - \frac{1}{2} \log^2(N+x) \right)$$

We also see that

$$\gamma_p(1) = \lim_{N \to \infty} \left( \sum_{n=0}^{N} \frac{\log^p(n+1)}{n+1} - \frac{\log^{p+1}(N+1)}{p+1} \right)$$

$$= \lim_{N \to \infty} \left( \sum_{n=1}^{N} \frac{\log^p(n+1)}{n+1} - \frac{\log^{p+1}(N+1)}{p+1} \right)$$

$$= \lim_{N \to \infty} \left( \sum_{k=1}^{N} \frac{\log^p k}{k} + \frac{\log^p(N+1)}{N+1} - \frac{\log^{p+1}(N+1)}{p+1} \right)$$

and since $\lim_{N \to \infty} \frac{\log^p(N+1)}{N+1} = 0$ this becomes

$$= \lim_{N \to \infty} \left( \sum_{k=1}^{N} \frac{\log^p k}{k} - \frac{\log^{p+1}(N+1)}{p+1} \right)$$

$$= \lim_{N \to \infty} \left( \sum_{k=1}^{N} \frac{\log^p k}{k} - \frac{\log^{p+1} N}{p+1} + \frac{\log^{p+1} N}{p+1} - \frac{\log^{p+1}(N+1)}{p+1} \right)$$

Therefore, assuming that $\lim_{N \to \infty} \left[ \log^{p+1} N - \log^{p+1}(N+1) \right] = 0$, we have another proof of the well-known result (see for example the 1955 paper, "The power series coefficients of $\varsigma(s)$", by Briggs and Chowla [35a])

$$\gamma_p = \gamma_p(1) = \lim_{N \to \infty} \left( \sum_{k=1}^{N} \frac{\log^p k}{k} - \frac{\log^{p+1} N}{p+1} \right) = \lim_{N \to \infty} \left( \sum_{k=1}^{N} \frac{\log^p k}{k} - \int_{1}^{N} \frac{\log^p x}{x} dx \right)$$

As shown by Dilcher in [54a] we also have for $k \geq 0$

(4.3.226iii) $$(-1)^{k+1} \varsigma_a^{(k)}(1) = \sum_{n=1}^{\infty} (-1)^n \frac{\log^k n}{n} = \sum_{j=0}^{k-1} \binom{k}{j} (-1)^j \gamma_j \log^{k-j} 2 - \frac{\log^{k+1} 2}{k+1}$$

of which (C.61) in Volume VI is a particular case. It should be noted that I have inserted a factor of $(-1)^j$ which appears to be required in (4.3.226iii). We then see that (see also (4.3.237) and (4.4.44j))

(4.3.226iv) $$\varsigma_a^{(2)}(1) = \frac{1}{3}\log^3 2 - \gamma \log^2 2 + 2\gamma_1 \log 2$$

See also the paper by Collins [46aa] entitled "The role of Bell polynomials in integration" where the integrals $\int_0^\infty \frac{\log^n x}{e^x + 1} dx$ are evaluated in terms of the Stieltjes constants. It may of course be noted that

$$\int_0^\infty \frac{\log^n x}{e^x+1} dx = \frac{d^n}{ds^n} \int_0^\infty \frac{x^{s-1}}{e^x+1} dx \bigg|_{s=1} = \frac{d^n}{ds^n}[\Gamma(s)\varsigma_a(s)]\bigg|_{s=1}$$

From (4.3.226ii) we have

$$\gamma_1(x) - \gamma_1(1) = \lim_{N\to\infty}\left( \sum_{n=0}^N \left[\frac{\log(n+x)}{n+x} - \frac{\log(n+1)}{n+1}\right] - \frac{1}{2}\log^2(N+x) + \frac{1}{2}\log^2(N+1) \right)$$

We see that

$$\log^2(N+1) - \log^2(N+x) = [\log(N+1) + \log(N+x)][\log(N+1) - \log(N+x)]$$

$$= \log[(N+1)(N+x)]\log\frac{(N+1)}{(N+x)}$$

$$= \log\left[N^2\left(1+\frac{1}{N}\right)\left(1+\frac{x}{N}\right)\right]\log\frac{\left(1+\frac{1}{N}\right)}{\left(1+\frac{x}{N}\right)}$$

$$= \log\left[N^2\left(1+\frac{1}{N}\right)\left(1+\frac{x}{N}\right)\right]\log\frac{\left(1+\frac{1}{N}\right)}{\left(1+\frac{x}{N}\right)} = 2\log N \log\frac{\left(1+\frac{1}{N}\right)}{\left(1+\frac{x}{N}\right)} + \log\left(1+\frac{1}{N}\right)\log\frac{\left(1+\frac{1}{N}\right)}{\left(1+\frac{x}{N}\right)}$$

$$+\log\left(1+\frac{x}{N}\right)\log\frac{\left(1+\frac{1}{N}\right)}{\left(1+\frac{x}{N}\right)}$$

Therefore, using L'Hôpital's rule we may prove that

$$\lim_{N\to\infty}[\log^2(N+1)-\log^2(N+x)]=0$$

and this limit is almost intuitively obvious.

Hence we see that (as originally obtained in (4.3.223b))

(4.3.226v) $$\gamma_1(x)-\gamma_1(1)=\sum_{n=0}^{\infty}\left[\frac{\log(n+x)}{n+x}-\frac{\log(n+1)}{n+1}\right]$$

More generally, with the identity

$$a^{p+1}-b^{p+1}=(a-b)(a^p+ba^{p-1}+b^2a^{p-2}+\ldots+b^p)$$

we may also prove that

$$\lim_{N\to\infty}[\log^{p+1}(N+1)-\log^{p+1}(N+x)]=0$$

The following alternative proof is quite succinct. Define $f(x)=\log^{p+1}(N+x)$ and, by the mean value theorem of calculus, we have

$$\log^{p+1}(N+x)-\log^{p+1}(N+1)=(x-1)(p+1)\frac{\log^p(N+\alpha)}{N+\alpha}$$

where $1<\alpha<x$. Successive applications of L'Hôpital's rule easily show that

$$\lim_{N\to\infty}\frac{\log^p(N+\alpha)}{N+\alpha}=0$$

Hence we have

(4.3.227) $$\gamma_p(x)-\gamma_p(1)=\sum_{n=0}^{\infty}\left[\frac{\log^p(n+x)}{n+x}-\frac{\log^p(n+1)}{n+1}\right]$$

and in particular we have

(4.3.228) $$\gamma_0(x) - \gamma_0(1) = \sum_{n=0}^{\infty}\left[\frac{1}{n+x} - \frac{1}{n+1}\right] = -(x-1)\sum_{n=0}^{\infty}\frac{1}{(n+1)(n+x)} = -\gamma - \psi(x)$$

We therefore have

$$\gamma_0'(x) = -\sum_{n=0}^{\infty}\frac{1}{(n+x)^2} = -\psi'(x)$$

and hence upon integration we see that

$$\gamma_0(x) = -\psi(x) + c$$

We already know that the integration constant $c = 0$ but we may also determine this in the following way:

As $x \to \infty$ we see that

$$\sum_{j=0}^{k}\binom{k}{j}(-1)^j \log(x+j) = \log\left[\prod_{j=0}^{k}(x+j)^{(-1)^j\binom{k}{j}}\right]$$

$$\sim \log\left[x^{\sum_{j=0}^{k}\binom{k}{j}(-1)^j}\right] = \delta_{k,0}\log x$$

Therefore we have

$$\sum_{k=0}^{\infty}\frac{1}{k+1}\sum_{j=0}^{k}\binom{k}{j}(-1)^j \log(x+j) \to \log x \text{ as } x \to \infty$$

and hence $-\gamma_0(x) \to \log x$ as $x \to \infty$. We have

$$\lim_{x\to\infty}[\gamma_0(x) + \log x] = 0$$

and therefore

$$\lim_{x\to\infty}[-\psi(x) + c + \log x] = 0$$

which implies that $c = 0$ since we already know from (E.66a) in Volume VI that $\lim_{x\to\infty}[-\psi(x) + \log x] = 0$.

We easily see that

(4.3.228a)
$$\gamma_p(x) - \gamma_p(1) = \sum_{n=0}^{\infty}\left[\frac{\log^p(n+x)}{n+x} - \frac{\log^p(n+1)}{n+1}\right] = (-1)^p \lim_{s \to 1} \frac{\partial^p}{\partial s^p}[\varsigma(s,x) - \varsigma(s,1)]$$

$$= (-1)^p \lim_{s \to 1} \frac{\partial^p}{\partial s^p}[\varsigma(s,x) - \varsigma(s)]$$

From this we easily see that

(4.3.228b) $\quad \gamma_p(x) - \gamma_p(y) = \lim_{s \to 1}(-1)^p \frac{\partial^p}{\partial s^p}[\varsigma(s,x) - \varsigma(s,y)]$

which we shall derive in a slightly different manner in (4.3.233e). It may also be readily derived from (4.3.207). This rather innocuous little equation actually turns out to be quite useful later in (4.3.233f). With $p = 0$ and $y = 1$ we obtain

(4.3.228c) $\quad \gamma_0(x) - \gamma = \lim_{s \to 1}[\varsigma(s,x) - \varsigma(s)]$

which we saw in (4.3.204a). We see that

$$\lim_{s \to 1}[\varsigma(s,x) - \varsigma(s)] = \lim_{s \to 1}\left(\frac{1}{x} + \sum_{n=1}^{\infty}\left[\frac{1}{(n+x)^s} - \frac{1}{n^s}\right]\right)$$

$$= \frac{1}{x} + \sum_{n=1}^{\infty}\left[\frac{1}{n+x} - \frac{1}{n}\right]$$

and, using (E.14) in Volume VI, we have

$$\frac{1}{x} + \sum_{n=1}^{\infty}\left[\frac{1}{n+x} - \frac{1}{n}\right] = -\gamma - \psi(x)$$

which gives us another derivation of (4.3.228c).

From (4.3.228b) we also see that

(4.3.228d) $\quad \varsigma'\left(1,\frac{1}{2}\right) - \varsigma'(1,1) = -\gamma_1\left(\frac{1}{2}\right) + \gamma_1 = \log^2 2 + 2\gamma$

as reported by Adamchik in [2a]. See also (4.3.233b) et seq.

Differentiation of (4.3.228a) with respect to $x$ gives us

$$\frac{d}{dx}\gamma_p(x) = (-1)^p \frac{\partial}{\partial x} \lim_{s \to 1} \frac{\partial^p}{\partial s^p}[\varsigma(s,x) - \varsigma(s)]$$

and, since the partial derivatives commute in the region where $\varsigma(s,x)$ is analytic, we have

$$= (-1)^p \lim_{s \to 1} \frac{\partial^p}{\partial s^p} \frac{\partial}{\partial x}[\varsigma(s,x) - \varsigma(s)]$$

$$= (-1)^{p+1} \lim_{s \to 1} \frac{\partial^p}{\partial s^p}[s\varsigma(s+1,x)]$$

since $\frac{\partial}{\partial x}\varsigma(s,x) = -s\varsigma(s+1,x)$. With $p=1$ we obtain

$$\frac{d}{dx}\gamma_1(x) = \varsigma'(2,x) + \varsigma(2,x)$$

which we have previously seen in (4.3.223b).

This may also be obtained directly as follows. Referring to (4.3.227) we see that

$$\frac{d}{dx}\gamma_p(x) = p\sum_{n=0}^{\infty} \frac{\log^{p-1}(n+x)}{(n+x)^2} - \sum_{n=0}^{\infty} \frac{\log^p(n+x)}{(n+x)^2}$$

and in particular

$$\frac{d}{dx}\gamma_1(x) = \sum_{n=0}^{\infty} \frac{1}{(n+x)^2} - \sum_{n=0}^{\infty} \frac{\log(n+x)}{(n+x)^2}$$

$$= \varsigma(2,x) + \varsigma'(2,x)$$

We now recall (4.3.112b)

$$(s-1)\varsigma'(s,x) + \varsigma(s,x) = -\sum_{n=0}^{\infty} \frac{1}{n+1} \sum_{k=0}^{n} \binom{n}{k}(-1)^k \frac{\log(x+k)}{(x+k)^{s-1}}$$

and note that

$$\varsigma'(2,x) + \varsigma(2,x) = -\sum_{n=0}^{\infty} \frac{1}{n+1} \sum_{k=0}^{n} \binom{n}{k}(-1)^k \frac{\log(x+k)}{x+k}$$

Therefore we get as expected

$$\frac{d}{dx}\gamma_1(x) = -\sum_{n=0}^{\infty}\frac{1}{n+1}\sum_{k=0}^{n}\binom{n}{k}(-1)^k\frac{\log(x+k)}{x+k}$$

Integration results in

$$\gamma_1(x)-\gamma_1(1) = -\frac{1}{2}\sum_{n=0}^{\infty}\frac{1}{n+1}\sum_{k=0}^{n}\binom{n}{k}(-1)^k\log^2(x+k)+\frac{1}{2}\sum_{n=0}^{\infty}\frac{1}{n+1}\sum_{k=0}^{n}\binom{n}{k}(-1)^k\log^2(1+k)$$

and by symmetry we see that

$$\gamma_1(x) = -\frac{1}{2}\sum_{n=0}^{\infty}\frac{1}{n+1}\sum_{k=0}^{n}\binom{n}{k}(-1)^k\log^2(x+k)+c$$

It remains for us to determine the integration constant. Employing the same method as before, as $x \to \infty$ we see that

$$\sum_{j=0}^{k}\binom{k}{j}(-1)^j\log^2(x+j) = \sum_{j=0}^{k}\binom{k}{j}(-1)^j\log(x+j)\log(x+j)$$

$$= \sum_{j=0}^{k}\log(x+j)^{(-1)^j\binom{k}{j}}\log(x+j)$$

$$= \delta_{k,0}\log^2 x$$

Therefore we have

$$-\frac{1}{2}\sum_{k=0}^{\infty}\frac{1}{k+1}\sum_{j=0}^{k}\binom{k}{j}(-1)^j\log^2(x+j) \to -\frac{1}{2}\log^2 x \text{ as } x\to\infty$$

It is seen that

$$\gamma_1(x) = -\frac{1}{2}\log^2 x - \frac{1}{2}\sum_{n=1}^{\infty}\frac{1}{n+1}\sum_{k=0}^{n}\binom{n}{k}(-1)^k\log^2(x+k)+c$$

However, where we go from here I do not know! When $x=1$ we do however find that $c=0$ provided we already know that

$$\gamma_1(u) = -\frac{1}{2}\sum_{n=0}^{\infty}\frac{1}{n+1}\sum_{k=0}^{n}\binom{n}{k}(-1)^k \log^2(u+k)$$

□

We see that starting the summation of (4.3.227) at $n=1$ gives us

(4.3.228e) $\quad \gamma_p(x) - \gamma_p(1) - \dfrac{\log^p x}{x} = \sum_{n=1}^{\infty}\left[\dfrac{\log^p(n+x)}{n+x} - \dfrac{\log^p(n+1)}{n+1}\right]$

and thus, in view of the telescoping of the above series, as $x \to 0$ we have

$$\lim_{x\to 0}[\gamma_p(x) - \gamma_p(1) - \frac{\log^p x}{x}] = \lim_{x\to 0}\sum_{n=1}^{\infty}\left[\frac{\log^p(n+x)}{n+x} - \frac{\log^p(n+1)}{n+1}\right]$$

$$= \sum_{n=1}^{\infty}\left[\frac{\log^p n}{n} - \frac{\log^p(n+1)}{n+1}\right] = 0$$

We then have

(4.3.228f) $\quad \lim_{x\to 0}[\gamma_p(x) - \dfrac{\log^p x}{x}] = \gamma_p$

Therefore we get

$$\lim_{x\to 0}\left[\frac{\log^p x}{x} + \frac{1}{p+1}\sum_{n=0}^{\infty}\frac{1}{n+1}\sum_{k=0}^{n}\binom{n}{k}(-1)^k \log^{p+1}\frac{x+k}{1+k}\right] = 0$$

Using L'Hôpital's rule we see that

$$\lim_{x\to 0}\left[\frac{\dfrac{\log^p x}{x} + \dfrac{1}{p+1}\sum_{n=0}^{\infty}\dfrac{1}{n+1}\sum_{k=0}^{n}\binom{n}{k}(-1)^k \log^{p+1}\dfrac{x+k}{1+k}}{x}\right]$$

$$= \lim_{x\to 0}\left[\frac{p\log^{p-1} x - \log^p x}{x^2} + \sum_{n=0}^{\infty}\frac{1}{n+1}\sum_{k=0}^{n}\binom{n}{k}(-1)^k \frac{\log^p(x+k)}{x+k}\right]$$

$$= \lim_{x\to 0}\left[\frac{p\log^{p-1} x - \log^p x}{x^2} - \gamma'_p(x)\right]$$

By differentiating (4.3.228e) we get

$$\gamma'_p(x) - \frac{p\log^{p-1} x - \log^p x}{x^2} = -\sum_{n=1}^{\infty}\frac{\log^p(n+x)}{(n+x)^2} + p\sum_{n=1}^{\infty}\frac{\log^{p-1}(n+x)}{(n+x)^2}$$

and hence

$$\lim_{x\to 0}\left[\gamma'_p(x) - \frac{p\log^{p-1} x - \log^p x}{x^2}\right] = -\sum_{n=1}^{\infty}\frac{\log^p n}{n^2} + p\sum_{n=1}^{\infty}\frac{\log^{p-1} n}{n^2}$$

This gives us

(4.3.228f) $\quad\displaystyle\lim_{x\to 0}\left[\gamma'_p(x) - \frac{p\log^{p-1} x - \log^p x}{x^2}\right] = (-1)^{p+1}\varsigma^{(p)}(2) + (-1)^{p+1}p\,\varsigma^{(p-1)}(2)$

We also have

$$\frac{d}{dx}\gamma_2(x) = 2\sum_{n=0}^{\infty}\frac{\log(n+x)}{(n+x)^2} - \sum_{n=0}^{\infty}\frac{\log^2(n+x)}{(n+x)^2}$$

$$= -2\varsigma'(2,x) - \varsigma''(2,x)$$

We have

$$(s-1)\varsigma''(s,x) + 2\varsigma'(s,x) = \sum_{n=0}^{\infty}\frac{1}{n+1}\sum_{k=0}^{n}\binom{n}{k}(-1)^k\frac{\log^2(x+k)}{(x+k)^{s-1}}$$

and therefore we see that

$$\varsigma''(2,x) + 2\varsigma'(2,x) = \sum_{n=0}^{\infty}\frac{1}{n+1}\sum_{k=0}^{n}\binom{n}{k}(-1)^k\frac{\log^2(x+k)}{x+k}$$

Hence we get

$$\frac{d}{dx}\gamma_2(x) = -\sum_{n=0}^{\infty}\frac{1}{n+1}\sum_{k=0}^{n}\binom{n}{k}(-1)^k\frac{\log^2(x+k)}{x+k}$$

and upon integration we obtain

$$\gamma_2(x) - \gamma_2(1) = -\frac{1}{3}\sum_{n=0}^{\infty}\frac{1}{n+1}\sum_{k=0}^{n}\binom{n}{k}(-1)^k\log^3(x+k) + \frac{1}{3}\sum_{n=0}^{\infty}\frac{1}{n+1}\sum_{k=0}^{n}\binom{n}{k}(-1)^k\log^3(1+k)$$

and thus we have

$$\gamma_2(x) = -\frac{1}{3}\sum_{n=0}^{\infty}\frac{1}{n+1}\sum_{k=0}^{n}\binom{n}{k}(-1)^k \log^3(x+k) + c$$

Using (4.3.227) we see that

$$\sum_{p=0}^{\infty}\frac{(-1)^p}{p!}[\gamma_p(x) - \gamma_p(1)] = \sum_{p=0}^{\infty}\frac{(-1)^p}{p!}\sum_{n=0}^{\infty}\left[\frac{\log^p(n+x)}{n+x} - \frac{\log^p(n+1)}{n+1}\right]$$

$$= \sum_{n=0}^{\infty}\frac{1}{n+x}\sum_{p=0}^{\infty}\frac{(-1)^p \log^p(n+x)}{p!} - \sum_{n=0}^{\infty}\frac{1}{n+1}\sum_{p=0}^{\infty}\frac{(-1)^p \log^p(n+1)}{p!}$$

$$= \sum_{n=0}^{\infty}\frac{1}{n+x}\exp[-\log(n+x)] - \sum_{n=0}^{\infty}\frac{1}{n+1}\exp[-\log(n+1)]$$

$$= \sum_{n=0}^{\infty}\frac{1}{(n+x)^2} - \sum_{n=0}^{\infty}\frac{1}{(n+1)^2}$$

$$= \psi'(x) - \psi'(1)$$

This concurs with (4.3.216) above.

We see from (4.3.215) that

$$\int_1^x \gamma_0(u)\,du = -\sum_{n=0}^{\infty}\frac{1}{n+1}\sum_{k=0}^{n}\binom{n}{k}(-1)^k \int_1^x \log(u+k)\,du$$

$$= -\sum_{n=0}^{\infty}\frac{1}{n+1}\sum_{k=0}^{n}\binom{n}{k}(-1)^k [(u+k)\log(u+k) - u]\Big|_1^x$$

$$= x - 1 - \sum_{n=0}^{\infty}\frac{1}{n+1}\sum_{k=0}^{n}\binom{n}{k}(-1)^k (x+k)\log(x+k) + \sum_{n=0}^{\infty}\frac{1}{n+1}\sum_{k=0}^{n}\binom{n}{k}(-1)^k (1+k)\log(1+k)$$

and using (4.3.74a) we obtain for $x > 0$

(4.3.229) $$\int_1^x \gamma_0(u)\,du = -\log\Gamma(x) = -[\varsigma'(0,x) - \varsigma'(0)]$$

as expected since $\gamma_0(u) = -\psi(u)$. Similarly, we have

$$\int_1^x \gamma_1(u)\,du = -\frac{1}{2}\sum_{n=0}^{\infty}\frac{1}{n+1}\sum_{k=0}^{n}\binom{n}{k}(-1)^k\int_1^x \log^2(u+k)\,du$$

$$= -\frac{1}{2}\sum_{n=0}^{\infty}\frac{1}{n+1}\sum_{k=0}^{n}\binom{n}{k}(-1)^k(u+k)[\log^2(u+k)-2\log(u+k)+2]\Big|_1^x$$

We recall (4.3.107a)

$$(s-1)\varsigma'(s,u)+\varsigma(s,u) = -\sum_{n=0}^{\infty}\frac{1}{n+1}\sum_{k=0}^{n}\binom{n}{k}(-1)^k\frac{\log(u+k)}{(u+k)^{s-1}}$$

and differentiation results in

$$(s-1)\varsigma''(s,u)+2\varsigma'(s,u) = \sum_{n=0}^{\infty}\frac{1}{n+1}\sum_{k=0}^{n}\binom{n}{k}(-1)^k\frac{\log^2(u+k)}{(u+k)^{s-1}}$$

$$(s-1)\varsigma^{(3)}(s,u)+3\varsigma^{(2)}(s,u) = -\sum_{n=0}^{\infty}\frac{1}{n+1}\sum_{k=0}^{n}\binom{n}{k}(-1)^k\frac{\log^3(u+k)}{(u+k)^{s-1}}$$

With $s=0$ we get

$$-\varsigma'(0,u)+\varsigma(0,u) = -\sum_{n=0}^{\infty}\frac{1}{n+1}\sum_{k=0}^{n}\binom{n}{k}(-1)^k(u+k)\log(u+k)$$

$$-\varsigma''(0,u)+2\varsigma'(0,u) = \sum_{n=0}^{\infty}\frac{1}{n+1}\sum_{k=0}^{n}\binom{n}{k}(-1)^k(u+k)\log^2(u+k)$$

$$-\varsigma^{(3)}(0,u)+3\varsigma^{(2)}(0,u) = -\sum_{n=0}^{\infty}\frac{1}{n+1}\sum_{k=0}^{n}\binom{n}{k}(-1)^k(u+k)\log^3(u+k)$$

Using the Leibniz rule for differentiation we note that

(4.3.229a) $$\frac{\partial^n}{\partial s^n}[(s-1)\varsigma(s,u)] = (s-1)\frac{\partial^n}{\partial s^n}\varsigma(s,u)+n\frac{\partial^{n-1}}{\partial s^{n-1}}\varsigma(s,u)$$

and using the above we obtain

(4.3.230)

$$-\varsigma^{(m)}(0,u) + m\varsigma^{(m-1)}(0,u) = (-1)^m \sum_{n=0}^{\infty} \frac{1}{n+1} \sum_{k=0}^{n} \binom{n}{k} (-1)^k (u+k) \log^m(u+k)$$

With $u = 1$ we get

(4.3.230i)

$$-\varsigma^{(m)}(0) + m\varsigma^{(m-1)}(0) = (-1)^m \sum_{n=0}^{\infty} \frac{1}{n+1} \sum_{k=0}^{n} \binom{n}{k} (-1)^k (1+k) \log^m(1+k)$$

and therefore we see that

(4.3.230ii)

$$-\varsigma^{(m)}(0) + m\varsigma^{(m-1)}(0) + (-1)^m m\gamma_{m-1} = (-1)^m \sum_{n=0}^{\infty} \frac{1}{n+1} \sum_{k=0}^{n} \binom{n}{k} (-1)^k k \log^m(1+k)$$

Hence, using Apostol's paper [14aa], we may express $\sum_{n=0}^{\infty} \frac{1}{n+1} \sum_{k=0}^{n} \binom{n}{k} (-1)^k k \log^m(1+k)$ in terms of Stieltjes constants and other known constants.

For example we have with $m = 1$

$$-\varsigma^{(1)}(0) + \varsigma^{(0)}(0) - \gamma = -\sum_{n=0}^{\infty} \frac{1}{n+1} \sum_{k=0}^{n} \binom{n}{k} (-1)^k k \log(1+k)$$

which gives us

(4.3.230iii) $\sum_{n=0}^{\infty} \frac{1}{n+1} \sum_{k=0}^{n} \binom{n}{k} (-1)^k k \log(1+k) = -\frac{1}{2} \log(2\pi) - \frac{1}{2} + \gamma$

When $m = 2$ we get

$$-\varsigma^{(2)}(0) + 2\varsigma^{(1)}(0) + 2\gamma_1 = \sum_{n=0}^{\infty} \frac{1}{n+1} \sum_{k=0}^{n} \binom{n}{k} (-1)^k k \log^2(1+k)$$

and hence we have

(4.3.230iv)

$$\frac{1}{2} \log^2(2\pi) + \frac{1}{4} \varsigma(2) - \frac{1}{2} \gamma^2 + \gamma_1 - \log(2\pi) = \sum_{n=0}^{\infty} \frac{1}{n+1} \sum_{k=0}^{n} \binom{n}{k} (-1)^k k \log^2(1+k)$$

In passing, we note that

$$\varsigma'(s,u) = -\frac{1}{s-1}\varsigma(s,u) - \frac{1}{s-1}\sum_{n=0}^{\infty}\frac{1}{n+1}\sum_{k=0}^{n}\binom{n}{k}(-1)^k \frac{\log(u+k)}{(u+k)^{s-1}}$$

$$\varsigma''(s,u) = \frac{2}{(s-1)^2}\varsigma(s,u) + \frac{2}{(s-1)^2}\sum_{n=0}^{\infty}\frac{1}{n+1}\sum_{k=0}^{n}\binom{n}{k}(-1)^k \frac{\log(u+k)}{(u+k)^{s-1}}$$

$$+ \frac{1}{s-1}\sum_{n=0}^{\infty}\frac{1}{n+1}\sum_{k=0}^{n}\binom{n}{k}(-1)^k \frac{\log^2(u+k)}{(u+k)^{s-1}}$$

$$\varsigma^{(3)}(s,u) = -\frac{3!}{(s-1)^3}\varsigma(s,u) - \frac{3!}{(s-1)^3}\sum_{n=0}^{\infty}\frac{1}{n+1}\sum_{k=0}^{n}\binom{n}{k}(-1)^k \frac{\log(u+k)}{(u+k)^{s-1}}$$

$$- \frac{3}{(s-1)^2}\sum_{n=0}^{\infty}\frac{1}{n+1}\sum_{k=0}^{n}\binom{n}{k}(-1)^k \frac{\log^2(u+k)}{(u+k)^{s-1}}$$

$$- \frac{1}{s-1}\sum_{n=0}^{\infty}\frac{1}{n+1}\sum_{k=0}^{n}\binom{n}{k}(-1)^k \frac{\log^3(u+k)}{(u+k)^{s-1}}$$

We then obtain

$$\int_{1}^{x}\gamma_1(u)\,du = \frac{1}{2}\varsigma''(0,u) - \varsigma'(0,u) + \varsigma'(0,u) - \varsigma(0,u) - u + \frac{1}{2}\Big|_{1}^{x}$$

Hence we get

(4.3.231) $$\int_{1}^{x}\gamma_1(u)\,du = \frac{1}{2}[\varsigma''(0,x) - \varsigma''(0)]$$

Moving on we have

$$\int_{1}^{x}\gamma_2(u)\,du = -\frac{1}{3}\sum_{n=0}^{\infty}\frac{1}{n+1}\sum_{k=0}^{n}\binom{n}{k}(-1)^k \int_{1}^{x}\log^3(u+k)\,du$$

$$\int_{1}^{x}\log^3(u+k)\,du = (u+k)[\log^3(u+k) - 3\log^2(u+k) + 6\log(u+k) - 6]\Big|_{1}^{x}$$

Therefore we obtain

$$\int_1^x \gamma_2(u)\,du = \left(-\frac{1}{3}\varsigma^{(3)}(0,u) + \varsigma^{(2)}(0,u)\right) + \left(-\varsigma^{(2)}(0,u) + 2\varsigma^{(1)}(0,u)\right) + \left(-2\varsigma^{(1)}(0,u) + 2\varsigma(0,u)\right)$$

$$+ 2\left(u - \frac{1}{2}\right)\Big|_1^x$$

and thus we have

(4.3.231i)  $$\int_1^x \gamma_2(u)\,du = -\frac{1}{3}[\varsigma^{(3)}(0,x) - \varsigma^{(3)}(0)]$$

More generally we have

$$\int_1^x \gamma_p(u)\,du = -\frac{1}{p+1}\sum_{n=0}^{\infty}\frac{1}{n+1}\sum_{k=0}^{n}\binom{n}{k}(-1)^k \int_1^x \log^{p+1}(u+k)\,du$$

With the substitution $y = \log(u+k)$ we get

$$\int \log^{p+1}(u+k)\,du = \int y^{p+1}e^y\,dy = (-1)^{p+1}(p+1)!\,e^y \sum_{m=0}^{p+1}(-1)^m \frac{y^m}{m!}$$

and therefore

$$\int \log^{p+1}(u+k)\,du = (-1)^{p+1}(p+1)!\sum_{m=0}^{p+1}\frac{(-1)^m}{m!}(u+k)\log^m(u+k)$$

Hence we obtain

$$\int_1^x \gamma_p(u)\,du = (-1)^{p+1}p!\sum_{m=0}^{p+1}\frac{(-1)^m}{m!}\sum_{n=0}^{\infty}\frac{1}{n+1}\sum_{k=0}^{n}\binom{n}{k}(-1)^k(1+k)\log^m(x+k)$$

$$-(-1)^{p+1}p!\sum_{m=0}^{p+1}\frac{(-1)^m}{m!}\sum_{n=0}^{\infty}\frac{1}{n+1}\sum_{k=0}^{n}\binom{n}{k}(-1)^k(1+k)\log^m(1+k)$$

and using (4.3.107a) this may be written as

$$\int_1^x \gamma_p(u)\,du = p!\sum_{m=1}^{p+1}\frac{(-1)^m}{m!}[\varsigma^{(m)}(0,x) - m\varsigma^{(m-1)}(0,x)]$$

$$-p!\sum_{m=1}^{p+1}\frac{(-1)^m}{m!}[\varsigma^{(m)}(0,1)-m\varsigma^{(m-1)}(0,1)]$$

Hence we obtain the general result

(4.3.231ii) $$\int_1^x \gamma_p(u)\,du = \frac{(-1)^{p+1}}{p+1}[\varsigma^{(p+1)}(0,x)-\varsigma^{(p+1)}(0)]$$

We saw above that

$$\int_1^x \gamma_1(u)\,du = \frac{1}{2}[\varsigma''(0,x)-\varsigma''(0)]$$

and differentiation gives us

$$\gamma_1(x) = \frac{1}{2}\frac{d}{dx}\varsigma''(0,x) = \frac{1}{2}\frac{\partial}{\partial x}\frac{\partial^2}{\partial s^2}\varsigma(s,x)\bigg|_{s=0}$$

We have

$$\frac{\partial}{\partial x}\frac{\partial^2}{\partial s^2}\varsigma(s,x) = -s\varsigma^{(2)}(s+1,x)-2\varsigma^{(1)}(s+1,x)$$

and hence we get

$$\gamma_1(x) = -\frac{1}{2}\lim_{s\to 0}[s\varsigma^{(2)}(s+1,x)+2\varsigma^{(1)}(s+1,x)]$$

We have therefore come around full circle since

$$\lim_{s\to 0}[s\varsigma^{(2)}(s+1,x)+2\varsigma^{(1)}(s+1,x)] = \lim_{s\to 1}[(s-1)\varsigma^{(2)}(s,x)+2\varsigma^{(1)}(s,x)]$$

and thus, as we have already seen, $\gamma_1(x) = -\frac{1}{2}\sum_{n=0}^{\infty}\frac{1}{n+1}\sum_{k=0}^{n}\binom{n}{k}(-1)^k \log^2(u+k)$.

By considering $\int_0^1 \varsigma(s,u)\,du = 0$ for $\mathrm{Re}(s)<1$ Coffey [45c] has shown that

(4.3.231iii) $$\sum_{p=0}^{\infty}\frac{(-1)^p}{p!}(s-1)^p \int_0^1 \gamma_p(u)\,du = \frac{1}{1-s}$$

but it should be noted that my expression for $\gamma_p(u)$ is not defined at $x=0$.

We see from (4.3.214) that

$$\sum_{p=0}^{\infty}\frac{(-1)^p}{p!}(s-1)^p\gamma_p(u) = \sum_{p=0}^{\infty}\frac{(-1)^p(s-1)^p}{p!}\left\{-\frac{1}{p+1}\sum_{n=0}^{\infty}\frac{1}{n+1}\sum_{k=0}^{n}\binom{n}{k}(-1)^k\log^{p+1}(u+k)\right\}$$

$$=\frac{1}{s-1}\sum_{n=0}^{\infty}\frac{1}{n+1}\sum_{k=0}^{n}\binom{n}{k}(-1)^k\sum_{p=0}^{\infty}\frac{(-1)^{p+1}(s-1)^{p+1}\log^{p+1}(u+k)}{(p+1)!}$$

$$=\frac{1}{s-1}\sum_{n=0}^{\infty}\frac{1}{n+1}\sum_{k=0}^{n}\binom{n}{k}(-1)^k\left(\exp[-(s-1)\log(u+k)]-1\right)$$

$$=\frac{1}{s-1}\sum_{n=0}^{\infty}\frac{1}{n+1}\sum_{k=0}^{n}\binom{n}{k}(-1)^k\left[\frac{1}{(u+k)^{s-1}}-1\right]$$

$$=\varsigma(s,u)-\frac{1}{s-1}$$

We have thus rediscovered that

$$\varsigma(s,u)=\frac{1}{s-1}+\sum_{p=0}^{\infty}\frac{(-1)^p}{p!}(s-1)^p\gamma_p(u)$$

$\square$

Coffey [45c] has shown that for integers $q\geq 2$

(4.3.232) $$\sum_{r=1}^{q-1}\gamma_p\left(\frac{r}{q}\right)=-\gamma_p+q(-1)^p\frac{\log^{p+1}q}{p+1}+q\sum_{j=0}^{p}\binom{p}{j}(-1)^j\gamma_{p-j}\log^j q$$

and, taking the simplest case $q=2$, we get

(4.3.232a) $$\gamma_p\left(\frac{1}{2}\right)=-\gamma_p+2(-1)^p\frac{\log^{p+1}2}{p+1}+2\sum_{j=0}^{p}\binom{p}{j}(-1)^j\gamma_{p-j}\log^j 2$$

With $p=0$ we get

$$\gamma_0\left(\frac{1}{2}\right)=-\gamma_0+2\log 2+2\sum_{j=0}^{0}\binom{1}{j}(-1)^j\gamma_{0-j}\log^j 2$$

$$= \gamma + 2\log 2 = -\psi(1/2)$$

in agreement with (4.3.215).

With $p=1$ we see that

$$\gamma_1\left(\frac{1}{2}\right) = -\gamma_1 - \log^2 2 + 2\sum_{j=0}^{1} \binom{1}{j}(-1)^j \gamma_{1-j} \log^j 2$$

$$= -\gamma_1 - \log^2 2 + 2\gamma_1 - 2\gamma \log 2$$

and thus we have

(4.3.233) $\quad \gamma_1\left(\frac{1}{2}\right) = \gamma_1 - \log^2 2 - 2\gamma \log 2$

A different proof of this is given in (4.3.263).

Using (4.3.233) in conjunction with (4.3.226iv)

$$\gamma_1(x) - \gamma_1(1) = \sum_{n=0}^{\infty}\left[\frac{\log(n+x)}{n+x} - \frac{\log(n+1)}{n+1}\right]$$

we obtain with $x=1/2$

$$\log^2 2 + 2\gamma \log 2 = \sum_{n=0}^{\infty}\left[\frac{\log(n+1)}{n+1} - 2\frac{\log(2n+1)}{2n+1} + 2\frac{\log 2}{2n+1}\right]$$

and this corrects a misprint in Dilcher's paper [54a].

We also have

(4.3.233a) $\quad \gamma_2\left(\frac{1}{2}\right) = -\gamma_2 + \frac{2}{3}\log^3 2 + 2\sum_{j=0}^{1}\binom{1}{j}(-1)^j \gamma_{2-j} \log^j 2$

$$= \gamma_2 - 2\gamma_1 \log 2 + \frac{2}{3}\log^3 2$$

□

Integrating (4.3.226iv) gives us

$$\int_1^u \gamma_1(x)dx - (u-1)\gamma_1 = \sum_{n=0}^{\infty}\left[\frac{1}{2}\log^2(n+u) - \frac{1}{2}\log^2(n+1) - (u-1)\frac{\log(n+1)}{n+1}\right]$$

$$= \sum_{n=0}^{\infty}\left[\frac{1}{2}\log^2(n+u)-\frac{1}{2}\log^2(n+1)+(u-1)\left(\frac{\log(n+u)}{n+u}-\frac{\log(n+1)}{n+1}\right)-(u-1)\frac{\log(n+u)}{n+u}\right]$$

and using (4.3.226iv) this becomes

$$= \sum_{n=0}^{\infty}\left[\frac{1}{2}\log^2(n+u)-\frac{1}{2}\log^2(n+1)-(u-1)\frac{\log(n+u)}{n+u}\right]+(u-1)[\gamma_1(u)-\gamma_1]$$

Then referring to (4.3.231) we may write this as

$$\frac{1}{2}\varsigma''(0,u)-\frac{1}{2}\varsigma''(0)-(u-1)\gamma_1 =$$

$$\sum_{n=0}^{\infty}\left[\frac{1}{2}\log^2(n+u)-\frac{1}{2}\log^2(n+1)-(u-1)\frac{\log(n+u)}{n+u}\right]+(u-1)[\gamma_1(u)-\gamma_1]$$

We have

$$\varsigma\left(s,\frac{1}{2}\right)=(2^s-1)\varsigma(s)$$

$$\varsigma'\left(s,\frac{1}{2}\right)=(2^s-1)\varsigma'(s)+2^s\varsigma(s)\log 2$$

$$\varsigma''\left(s,\frac{1}{2}\right)=(2^s-1)\varsigma''(s)+2.2^s\varsigma'(s)\log 2+2^s\varsigma(s)\log^2 2$$

and with $s=0$ we have

$$\varsigma''\left(0,\frac{1}{2}\right)=2\varsigma'(0)\log 2+\varsigma(0)\log^2 2$$

$$=-\log(2\pi)\log 2-\frac{1}{2}\log^2 2$$

Then with $u=1/2$ and using (4.3.233) we obtain

$$-\frac{1}{2}\log(2\pi)\log 2-\frac{1}{4}\log^2 2-\frac{1}{2}\varsigma''(0)+\frac{1}{2}\gamma_1 =$$

$$\sum_{n=0}^{\infty}\left[\frac{1}{2}\log^2\left(n+\frac{1}{2}\right)-\frac{1}{2}\log^2(n+1)+\frac{1}{2n+1}\log\left(n+\frac{1}{2}\right)\right]+\frac{1}{2}\log^2 2+\gamma\log 2$$

We note from (4.3.247) that

$$\varsigma''(0) = \gamma_1 + \frac{1}{2}\gamma^2 - \frac{1}{24}\pi^2 - \frac{1}{2}\log^2(2\pi)$$

and we obtain (on the assumption that integrating (4.3.226iv) term by term is valid)

$$\frac{1}{4}\log(2\pi)\log(\pi/2) - \frac{3}{4}\log^2 2 - \gamma\log 2 - \frac{1}{2}\gamma^2 + \frac{1}{48}\pi^2 =$$

$$\sum_{n=0}^{\infty}\left[\frac{1}{2}\log^2\left(n+\frac{1}{2}\right)-\frac{1}{2}\log^2(n+1)+\frac{1}{2n+1}\log\left(n+\frac{1}{2}\right)\right]$$

□

Adamchik [2a], using Rademacher's formula [14, p.261] for the Hurwitz zeta function, showed that

(4.3.233b)

$$\varsigma'\left(1,\frac{p}{q}\right) - \varsigma'\left(1,1-\frac{p}{q}\right) = \pi[\log(2\pi q)+\gamma]\cot\left(\frac{\pi p}{q}\right) - 2\pi\sum_{j=1}^{q-1}\log\Gamma\left(\frac{j}{q}\right)\sin\left(\frac{2\pi jp}{q}\right)$$

where $p$ and $q$ are positive integers and $p < q$. In particular we have

(4.3.233c) $\quad \varsigma'\left(1,\frac{1}{4}\right) - \varsigma'\left(1,\frac{3}{4}\right) = \pi\left[\gamma + 4\log 2 + 3\log\pi - 4\log\Gamma\left(\frac{1}{4}\right)\right]$

(4.3.233d) $\quad \varsigma'\left(1,\frac{1}{3}\right) - \varsigma'\left(1,\frac{2}{3}\right) = \frac{\pi}{2\sqrt{3}}\left[2\gamma - \log 3 + 8\log(2\pi) - 12\log\Gamma\left(\frac{1}{3}\right)\right]$

We have from (4.3.206) that

$$\varsigma(s,u) = \frac{1}{s-1} + \sum_{p=0}^{\infty}\frac{(-1)^p}{p!}\gamma_p(u)(s-1)^p$$

and therefore

$$\varsigma(s,u) - \varsigma(s,1-u) = \sum_{p=0}^{\infty} \frac{(-1)^p}{p!} [\gamma_p(u) - \gamma_p(1-u)](s-1)^p$$

Differentiation with respect to $s$ results in

(4.3.233di) $$\varsigma'(s,u) - \varsigma'(s,1-u) = \sum_{p=0}^{\infty} \frac{(-1)^p}{p!} p[\gamma_p(u) - \gamma_p(1-u)](s-1)^{p-1}$$

and in the limit as $s \to 1$ we have

(4.3.233e) $$\varsigma'(1,u) - \varsigma'(1,1-u) = -[\gamma_1(u) - \gamma_1(1-u)]$$

In fact, as seen in (4.3.228b), more generally we have

$$\gamma_p(x) - \gamma_p(y) = \lim_{s \to 1} (-1)^p \frac{\partial^p}{\partial s^p} [\varsigma(s,x) - \varsigma(s,y)]$$

Hence we obtain

(4.3.233f) $$\gamma_1\left(\frac{1}{4}\right) - \gamma_1\left(\frac{3}{4}\right) = -\pi\left[\gamma + 4\log 2 + 3\log \pi - 4\log \Gamma\left(\frac{1}{4}\right)\right]$$

and equivalently

(4.3.233g) $$\varsigma'\left(1,\frac{3}{4}\right) - \varsigma'\left(1,\frac{1}{4}\right) = \pi\left[\gamma + 4\log 2 + 3\log \pi - 4\log \Gamma\left(\frac{1}{4}\right)\right]$$

As mentioned above in (4.3.232), Coffey [45c] has shown that for integers $q \geq 2$

$$\sum_{r=1}^{q-1} \gamma_p\left(\frac{r}{q}\right) = -\gamma_p + q(-1)^p \frac{\log^{p+1} q}{p+1} + q \sum_{j=0}^{p} \binom{p}{j} (-1)^j \gamma_{p-j} \log^j q$$

and from (4.3.233) we saw that

$$\gamma_1\left(\frac{1}{2}\right) = \gamma_1 - \log^2 2 - 2\gamma \log 2$$

We see that

$$\sum_{r=1}^{3} \gamma_1\left(\frac{r}{4}\right) = -\gamma_1 - 4\log^2 4 + 4\sum_{j=0}^{1} \binom{1}{j} (-1)^j \gamma_{1-j} \log^j 4$$

$$= 3\gamma_1 - 16\log^2 2 - 8\gamma \log 2$$

and we also have

$$\sum_{r=1}^{3} \gamma_1\left(\frac{r}{4}\right) = \gamma_1\left(\frac{1}{4}\right) + \gamma_1\left(\frac{1}{2}\right) + \gamma_1\left(\frac{3}{4}\right)$$

$$= \gamma_1\left(\frac{1}{4}\right) + \gamma_1\left(\frac{3}{4}\right) + \gamma_1 - \log^2 2 - 2\gamma \log 2$$

Therefore we get

$$\gamma_1\left(\frac{1}{4}\right) + \gamma_1\left(\frac{3}{4}\right) = 2\gamma_1 - 15\log^2 2 - 6\gamma \log 2$$

and using (4.3.233a) we easily see that

(4.3.233h)
$$\gamma_1\left(\frac{1}{4}\right) = \frac{1}{2}[2\gamma_1 - 15\log^2 2 - 6\gamma \log 2] - \frac{1}{2}\pi\left[\gamma + 4\log 2 + 3\log \pi - 4\log \Gamma\left(\frac{1}{4}\right)\right]$$

$$\gamma_1\left(\frac{3}{4}\right) = \frac{1}{2}[2\gamma_1 - 15\log^2 2 - 6\gamma \log 2] + \frac{1}{2}\pi\left[\gamma + 4\log 2 + 3\log \pi - 4\log \Gamma\left(\frac{1}{4}\right)\right]$$

Assuming that each of the terms $\pi, \gamma, \log 2, \log \pi$ and $\log \Gamma\left(\frac{1}{4}\right)$ has weight equal to one and $\gamma_1$ has weight 2, then the above represents a homogenous equation of weight 2.

An alternative proof is contained in (4.3.268). Since $\psi\left(\frac{1}{4}\right) = -\gamma - \frac{1}{2}\pi - 3\log 2$ we may write this as

$$\gamma_1\left(\frac{1}{4}\right) = \gamma\psi\left(\frac{1}{4}\right) + \gamma^2 + \frac{1}{2}[2\gamma_1 - 15\log^2 2] - \frac{1}{2}\pi\left[4\log 2 + 3\log \pi - 4\log \Gamma\left(\frac{1}{4}\right)\right]$$

or

$$\gamma_1\left(\frac{1}{4}\right) = -\gamma\gamma\left(\frac{1}{4}\right) + \gamma^2 + \gamma_1 - \frac{15}{2}\log^2 2 - \frac{1}{2}\pi\left[4\log 2 + 3\log \pi - 4\log \Gamma\left(\frac{1}{4}\right)\right]$$

Similarly we have

$$\gamma_1\left(\frac{1}{3}\right) - \gamma_1\left(\frac{2}{3}\right) = -\left[\varsigma'\left(1,\frac{1}{3}\right) - \varsigma'\left(1,\frac{2}{3}\right)\right]$$

and using (4.3.233d) this becomes

$$= -\frac{\pi}{2\sqrt{3}}\left[2\gamma - \log 3 + 8\log(2\pi) - 12\log\Gamma\left(\frac{1}{3}\right)\right]$$

Equation (4.3.232) gives us

$$\gamma_1\left(\frac{1}{3}\right) + \gamma_1\left(\frac{2}{3}\right) = 2\gamma_1 - 3\gamma\log 3 - \frac{3}{2}\log^2 3$$

and this then gives us

(4.3.233i) $\gamma_1\left(\frac{1}{3}\right) = \gamma_1 - \frac{3}{2}\gamma\log 3 - \frac{3}{4}\log^2 3 - \frac{\pi}{4\sqrt{3}}\left[2\gamma - \log 3 + 8\log(2\pi) - 12\log\Gamma\left(\frac{1}{3}\right)\right]$

with a similar expression for $\gamma_1\left(\frac{2}{3}\right)$.

In [54a] Dilcher defined a generalised polygamma function by

$$\psi_p(x) = -\gamma_p - \frac{\log^p x}{x} - \sum_{n=1}^{\infty}\left[\frac{\log^p(n+x)}{n+x} - \frac{\log^p n}{n}\right]$$

$$= -\gamma_p - \frac{\log^p x}{x} - \sum_{n=1}^{\infty}\left[\frac{\log^p(n+x)}{n+x} - \frac{\log^p(n+1)}{n+1} - \frac{\log^p n}{n} + \frac{\log^p(n+1)}{n+1}\right]$$

$$= -\gamma_p - \frac{\log^p x}{x} - \sum_{n=1}^{\infty}\left[\frac{\log^p(n+x)}{n+x} - \frac{\log^p(n+1)}{n+1}\right]$$

Hence, comparing this with (4.3.228e), we see that

(4.3.233j)  $\psi_p(x) = -\gamma_p(x)$

This relationship was also noted by Coffey in [45c].

Dilcher [54a] also showed that

$$\psi_1\left(\frac{1}{3}\right) = -\gamma_1 + \frac{1}{2}\left[3\log 3 + \frac{\pi}{\sqrt{3}}\right]\gamma + \frac{3}{4}\log^3 3 + \pi\sqrt{3}\left[\frac{2}{3}\log(2\pi) - \frac{1}{12}\log 3 - \log\Gamma\left(\frac{1}{3}\right)\right]$$

and a little algebra shows that this is the same as (4.3.233i) above. It was also shown by Dilcher [54a] that

$$\psi_k\left(\frac{1}{2}\right) = -\gamma_k\left(\frac{1}{2}\right) = (-1)^{k+1} 2\frac{\log^{k+1} 2}{k+1} - \gamma_k - 2\sum_{j=0}^{k}\binom{k}{j}(-1)^j \gamma_{k-j} \log^j 2$$

which is the same as Coffey's formula (4.3.232a).

With $s = 2$ in (4.3.233di) we obtain

(4.3.233k) $$\varsigma'(2,u) - \varsigma'(2,1-u) = \sum_{p=0}^{\infty}\frac{(-1)^p}{p!}p[\gamma_p(u) - \gamma_p(1-u)]$$

We see that

$$\sum_{p=0}^{\infty}\frac{(-1)^p p\gamma_p(u)}{p!} = \sum_{p=0}^{\infty}\frac{(-1)^p p}{p!}\left\{-\frac{1}{p+1}\sum_{n=0}^{\infty}\frac{1}{n+1}\sum_{k=0}^{n}\binom{n}{k}(-1)^k \log^{p+1}(u+k)\right\}$$

$$= \sum_{n=0}^{\infty}\frac{1}{n+1}\sum_{k=0}^{n}\binom{n}{k}(-1)^k \sum_{p=0}^{\infty}\frac{p(-1)^{p+1}\log^{p+1}(u+k)}{p!(p+1)}$$

We see that

$$\sum_{p=0}^{\infty}\frac{(-1)^{p+1}x^{p+1}}{(p+1)!} = e^{-x} - 1$$

and also that

$$\sum_{p=0}^{\infty}\frac{(-1)^{p+1}px^{p+1}}{(p+1)!} = x^2\frac{d}{dx}\left(\frac{e^{-x}-1}{x}\right)$$

We therefore have

$$\sum_{p=0}^{\infty}\frac{(-1)^{p+1}px^{p+1}}{(p+1)!} = 1 - (x+1)e^{-x}$$

and hence we get

$$\sum_{p=0}^{\infty}\frac{(-1)^{p+1}p\gamma_p(u)}{p!} = \sum_{n=0}^{\infty}\frac{1}{n+1}\sum_{k=0}^{n}\binom{n}{k}(-1)^k\left\{1-\frac{[1+\log(u+k)]}{u+k}\right\}$$

$$= \sum_{n=0}^{\infty}\frac{1}{n+1}\sum_{k=0}^{n}\binom{n}{k}(-1)^k - \sum_{n=0}^{\infty}\frac{1}{n+1}\sum_{k=0}^{n}\binom{n}{k}\frac{(-1)^k}{u+k} - \sum_{n=0}^{\infty}\frac{1}{n+1}\sum_{k=0}^{n}\binom{n}{k}(-1)^k\frac{\log(u+k)}{u+k}$$

$$= 1 - \varsigma(2,u) - \sum_{n=0}^{\infty}\frac{1}{n+1}\sum_{k=0}^{n}\binom{n}{k}(-1)^k\frac{\log(u+k)}{u+k}$$

Using (4.3.107a)

$$(s-1)\varsigma'(s,u) + \varsigma(s,u) = -\sum_{n=0}^{\infty}\frac{1}{n+1}\sum_{k=0}^{n}\binom{n}{k}(-1)^k\frac{\log(u+k)}{(u+k)^{s-1}}$$

we see that

$$\varsigma'(2,u) + \varsigma(2,u) = -\sum_{n=0}^{\infty}\frac{1}{n+1}\sum_{k=0}^{n}\binom{n}{k}(-1)^k\frac{\log(u+k)}{u+k}$$

This gives us

$$\sum_{p=0}^{\infty}\frac{(-1)^{p+1}p\gamma_p(u)}{p!} = 1 + \varsigma'(2,u)$$

and this verifies (4.3.233k).

Similarly, with $s = 3$ in (4.3.233di) we obtain

$$\varsigma'(3,u) - \varsigma'(3,1-u) = \sum_{p=0}^{\infty}\frac{(-1)^p 2^p}{p!}p[\gamma_p(u) - \gamma_p(1-u)]$$

Differentiating (4.3.217a)

$$\sum_{p=0}^{\infty}\frac{(-1)^p}{p!}t^p\gamma_p(u) = \varsigma(t+1,u) - \frac{1}{t}$$

we see that

$$\sum_{p=0}^{\infty}\frac{(-1)^p}{p!}pt^{p-1}\gamma_p(u) = \frac{\partial}{\partial t}\varsigma(t+1,u) + \frac{1}{t^2}$$

and with $t = 1$ we obtain the same result as before.

□

We recall from (4.3.210) that

$$\gamma_1(u) = -\frac{1}{2}\sum_{n=0}^{\infty}\frac{1}{n+1}\sum_{k=0}^{n}\binom{n}{k}(-1)^k \log^2(k+u)$$

and therefore

$$\gamma_1\left(\frac{1}{2}\right) = -\frac{1}{2}\sum_{n=0}^{\infty}\frac{1}{n+1}\sum_{k=0}^{n}\binom{n}{k}(-1)^k \log^2 \frac{(2k+1)}{2}$$

$$= -\frac{1}{2}\sum_{n=0}^{\infty}\frac{1}{n+1}\sum_{k=0}^{n}\binom{n}{k}(-1)^k [\log^2(2k+1) - 2\log 2 \log(2k+1) + \log^2 2]$$

$$= -\frac{1}{2}\sum_{n=0}^{\infty}\frac{1}{n+1}\sum_{k=0}^{n}\binom{n}{k}(-1)^k \log^2(2k+1) + \log 2 \sum_{n=0}^{\infty}\frac{1}{n+1}\sum_{k=0}^{n}\binom{n}{k}(-1)^k \log(2k+1) - \frac{1}{2}\log^2 2$$

From (4.3.74a) we have previously seen that

$$\sum_{n=0}^{\infty}\frac{1}{n+1}\sum_{k=0}^{n}\binom{n}{k}(-1)^k \log(2k+1) = -\gamma - \log 2$$

and accordingly we obtain

$$\gamma_1\left(\frac{1}{2}\right) = -\frac{1}{2}\sum_{n=0}^{\infty}\frac{1}{n+1}\sum_{k=0}^{n}\binom{n}{k}(-1)^k \log^2(2k+1) - \gamma \log 2 - \frac{3}{2}\log^2 2$$

Hence we have using (4.3.233)

(4.3.233l) $$\frac{1}{2}\sum_{n=0}^{\infty}\frac{1}{n+1}\sum_{k=0}^{n}\binom{n}{k}(-1)^k \log^2(2k+1) = \gamma_1 + \gamma \log 2 - \frac{1}{2}\log^2 2$$

□

Using (4.3.212) we see that

$$\int_0^u \gamma_1(x)\,dx = -\frac{1}{2}\sum_{n=0}^{\infty}\frac{1}{n+1}\sum_{k=0}^{n}\binom{n}{k}(-1)^k \int_0^u \log^2(x+k)\,dx$$

$$\int_0^u \log^2(x+k)\,dx = (x+k)\log^2(x+k) - 2(x+k)\log(x+k) + 2(x+k)\Big|_0^u$$

$$= (u+k)\log^2(u+k) - 2(u+k)\log(u+k) + 2(u+k)$$

$$-k\log^2(k) + 2k\log(k) - 2k$$

and we shall abbreviate this as $I(u) - I(0)$. Therefore we have

$$\int_0^u \gamma_1(x)\,dx = -\frac{1}{2}\sum_{n=0}^{\infty}\frac{1}{n+1}\sum_{k=0}^{n}\binom{n}{k}(-1)^k[I(u) - I(0)]$$

where

$$I(u) = (u+k)\{[\log(u+k) - 1]^2 + 1\}$$

With reference to

$$(s-1)\varsigma(s,u) = \sum_{n=0}^{\infty}\frac{1}{n+1}\sum_{k=0}^{n}\binom{n}{k}\frac{(-1)^k}{(u+k)^{s-1}}$$

we see that

$$\sum_{n=0}^{\infty}\frac{1}{n+1}\sum_{k=0}^{n}\binom{n}{k}(-1)^k(u+k)\log(u+k) = -\lim_{s\to 0}\frac{\partial}{\partial s}\sum_{n=0}^{\infty}\frac{1}{n+1}\sum_{k=0}^{n}\binom{n}{k}\frac{(-1)^k}{(u+k)^{s-1}}$$

$$= -\lim_{s\to 0}\frac{\partial}{\partial s}[(s-1)\varsigma(s,u)]$$

$$= \varsigma'(0,u) - \varsigma(0,u)$$

which we have in fact already seen from (4.3.108). Then using Lerch's identity (4.3.116) we see that

(4.3.233m) $\displaystyle\sum_{n=0}^{\infty}\frac{1}{n+1}\sum_{k=0}^{n}\binom{n}{k}(-1)^k(u+k)\log(u+k) = \log\Gamma(u) - \frac{1}{2}\log(2\pi) - \varsigma(0,u)$

We note that it is not possible to substitute $u = 0$ in the above equation. With $u = 1$ we get

$$\sum_{n=0}^{\infty}\frac{1}{n+1}\sum_{k=0}^{n}\binom{n}{k}(-1)^k(k+1)\log(k+1)=\frac{1}{2}-\frac{1}{2}\log(2\pi)$$

which we have seen before in (4.3.115).

$\square$

We see from (4.3.110) that

$$\varsigma(0,u)=-\sum_{n=0}^{\infty}\frac{1}{n+1}\sum_{k=0}^{n}\binom{n}{k}(-1)^k(u+k)$$

and therefore using (4.3.108) we have

$$\varsigma'(0,u)=\sum_{n=0}^{\infty}\frac{1}{n+1}\sum_{k=0}^{n}\binom{n}{k}(-1)^k(u+k)[\log(u+k)-1]$$

Furthermore we have

$$(4.3.233\text{n})\quad \sum_{n=0}^{\infty}\frac{1}{n+1}\sum_{k=0}^{n}\binom{n}{k}(-1)^k(u+k)\log^2(u+k)=\lim_{s\to 0}\frac{\partial^2}{\partial s^2}\sum_{n=0}^{\infty}\frac{1}{n+1}\sum_{k=0}^{n}\binom{n}{k}\frac{(-1)^k}{(u+k)^{s-1}}$$

$$=\lim_{s\to 0}\frac{\partial^2}{\partial s^2}[(s-1)\varsigma(s,u)]$$

$$=-\varsigma''(0,u)+2\varsigma'(0,u)$$

$$=-\varsigma''(0,u)+2\log\Gamma(u)-\log(2\pi)$$

We now differentiate (4.3.233n) to obtain

$$2\sum_{n=0}^{\infty}\frac{1}{n+1}\sum_{k=0}^{n}\binom{n}{k}(-1)^k\log(u+k)+\sum_{n=0}^{\infty}\frac{1}{n+1}\sum_{k=0}^{n}\binom{n}{k}(-1)^k\log^2(u+k)=-\frac{d}{du}\varsigma''(0,u)+2\psi(u)$$

We therefore see that

$$\frac{d}{du}\varsigma''(0,u)=-\sum_{n=0}^{\infty}\frac{1}{n+1}\sum_{k=0}^{n}\binom{n}{k}(-1)^k\log^2(u+k)=2\gamma_1(u)$$

It is readily seen that

$$-\frac{\partial}{\partial u}\frac{\partial^2}{\partial s^2}\varsigma(s,u) = -\frac{\partial^2}{\partial s^2}\frac{\partial}{\partial u}\varsigma(s,u)$$

$$= -\frac{\partial^2}{\partial s^2}[-s\varsigma(s+1,u)]$$

$$= s\varsigma''(s+1,u) + 2\varsigma'(s+1,u)$$

We also have

$$\sum_{n=0}^{\infty}\frac{1}{n+1}\sum_{k=0}^{n}\binom{n}{k}(-1)^k(u+k)\log^2(u+k) =$$

$$-\varsigma''(0,u) + 2\sum_{n=0}^{\infty}\frac{1}{n+1}\sum_{k=0}^{n}\binom{n}{k}(-1)^k(u+k)\log(u+k) + 2\varsigma(0,u)$$

and therefore we get

$$\varsigma''(0,u) = -\sum_{n=0}^{\infty}\frac{1}{n+1}\sum_{k=0}^{n}\binom{n}{k}(-1)^k(u+k)[\log(u+k)-1]^2$$

We also note for reference from Apostol's paper [14aa] that

(4.3.234) $$\varsigma''(0) = \gamma_1 + \frac{1}{2}\gamma^2 - \frac{1}{24}\pi^2 - \frac{1}{2}\log^2(2\pi)$$

and we may write this as

(4.3.234i)

$$\varsigma''(0) = -\frac{1}{2}\sum_{n=0}^{\infty}\frac{1}{n+1}\sum_{k=0}^{n}\binom{n}{k}(-1)^k \log^2(1+k) + \frac{1}{2}\left[\sum_{n=0}^{\infty}\frac{1}{n+1}\sum_{k=0}^{n}\binom{n}{k}(-1)^k \log(1+k)\right]^2$$

$$-\frac{1}{4}\varsigma(2) - 2\left[\sum_{n=0}^{\infty}\frac{1}{n+1}\sum_{k=0}^{n}\binom{n}{k}(-1)^k(1+k)[\log(1+k)+1]\right]^2$$

Upon differentiating (4.3.107a) we obtain

$$(s-1)\varsigma''(s) + 2\varsigma'(s) = \sum_{n=0}^{\infty}\frac{1}{n+1}\sum_{k=0}^{n}\binom{n}{k}\frac{(-1)^k \log^2(k+1)}{(k+1)^{s-1}}$$

and as $s \to 0$ we get

$$2\varsigma'(0) - \varsigma''(0) = \sum_{n=0}^{\infty} \frac{1}{n+1} \sum_{k=0}^{n} \binom{n}{k} (-1)^k (k+1) \log^2(k+1)$$

$$= \sum_{n=0}^{\infty} \frac{1}{n+1} \sum_{k=0}^{n} \binom{n}{k} (-1)^k k \log^2(k+1) + \sum_{n=0}^{\infty} \frac{1}{n+1} \sum_{k=0}^{n} \binom{n}{k} (-1)^k \log^2(k+1)$$

and we therefore see that

$$\sum_{n=0}^{\infty} \frac{1}{n+1} \sum_{k=0}^{n} \binom{n}{k} (-1)^k k \log^2(k+1) = 2\varsigma'(0) - \varsigma''(0) - \sum_{n=0}^{\infty} \frac{1}{n+1} \sum_{k=0}^{n} \binom{n}{k} (-1)^k \log^2(k+1)$$

$$= 2\varsigma'(0) - \varsigma''(0) + 2\gamma_1$$

and using (4.3.234) we then obtain

$$\sum_{n=0}^{\infty} \frac{1}{n+1} \sum_{k=0}^{n} \binom{n}{k} (-1)^k k \log^2(k+1) = \gamma_1 - \log(2\pi) - \frac{1}{2}\gamma^2 + \frac{1}{24}\pi^2 + \frac{1}{2}\log^2(2\pi)$$

In Coffey's recent paper [45c], "New results on the Stieltjes constants: Asymptotic and exact evaluation", it was shown that the Stieltjes constant could be written in the form

(4.3.235) $$\gamma_0 = \gamma = \frac{1}{2}\log 2 - \frac{1}{\log 2} \sum_{n=0}^{\infty} \frac{1}{2^{n+1}} \sum_{k=1}^{n} \binom{n}{k} (-1)^k \frac{\log(k+1)}{k+1}$$

Since $\log 1 = 0$ we may also write this as

$$\gamma_0 = \gamma = \frac{1}{2}\log 2 - \frac{1}{\log 2} \sum_{n=0}^{\infty} \frac{1}{2^{n+1}} \sum_{k=0}^{n} \binom{n}{k} (-1)^k \frac{\log(k+1)}{k+1}$$

and using (c.61) in Volume VI

$$\sum_{n=0}^{\infty} \frac{1}{2^{n+1}} \sum_{k=0}^{n} \binom{n}{k} \frac{(-1)^k \log(k+1)}{k+1} = -\varsigma'_a(1) = -\log 2 \left[ \gamma - \frac{\log 2}{2} \right]$$

it is seen that Coffey's result becomes evident. Coffey also reported that

(4.3.236)

$$\gamma_1 = \frac{1}{12}\log^2 2 - \frac{1}{2}\sum_{n=0}^{\infty}\frac{1}{2^{n+1}}\sum_{k=1}^{n}\binom{n}{k}(-1)^k\frac{\log(k+1)}{k+1} + \frac{1}{2\log 2}\sum_{n=0}^{\infty}\frac{1}{2^{n+1}}\sum_{k=1}^{n}\binom{n}{k}(-1)^k\frac{\log^2(k+1)}{k+1}$$

and, using (C.61), this may be written as

$$\gamma_1 = \frac{1}{12}\log^2 2 + \frac{1}{2}\log 2\left[\gamma - \frac{\log 2}{2}\right] + \frac{1}{2\log 2}\varsigma_a''(1)$$

which, in turn, is equivalent to (see also Dilcher's result (4.3.226iv))

(4.3.237) $$\varsigma_a''(1) = 2\gamma_1 \log 2 - \gamma \log^2 2 + \frac{1}{3}\log^2 2$$

The above is in fact a particular case of the general formula reported by Briggs and Chowla [35a] in 1955 where they showed that (see also Dilcher's result (4.3.226iv))

(4.3.238) $$\varsigma_a^{(k)}(1) = k!\sum_{r=1}^{k+1}\frac{(-1)^{r+1}\log^r 2}{r!}A_{k-r}$$

where $A_n = \frac{(-1)^n}{n!}\gamma_n$ and $A_{-1} = 1$.

In [45e] Coffey extended his results to show that

$$\gamma_m = -m!\sum_{l=1}^{m}\frac{B_{m-l+1}}{(m-l+1)!}\frac{\log^{m-l} 2}{l!}\sum_{n=1}^{\infty}\frac{1}{2^{n+1}}\sum_{k=1}^{n}\binom{n}{k}(-1)^k\frac{\log^l(k+1)}{k+1}$$

(4.3.239)

$$-\frac{1}{(m+1)\log 2}\sum_{n=1}^{\infty}\frac{1}{2^{n+1}}\sum_{k=1}^{n}\binom{n}{k}(-1)^k\frac{\log^{m+1}(k+1)}{k+1} - \frac{B_{m+1}}{m+1}\log^{m+1} 2$$

or alternatively

(4.3.239i)

$$\gamma_m = -m!\sum_{l=1}^{m}\frac{B_{m-l+1}}{(m-l+1)!}\frac{\log^{m-l} 2}{l!}(-1)^l\varsigma_a^{(l)}(1) - \frac{1}{(m+1)\log 2}(-1)^{m+1}\varsigma_a^{(m+1)}(1) - \frac{B_{m+1}}{m+1}\log^{m+1} 2$$

□

Differentiating the Hasse formula (3.12) we obtain

$$\varsigma'(s) = -\frac{1}{s-1}\sum_{n=0}^{\infty}\frac{1}{n+1}\sum_{k=0}^{n}\binom{n}{k}\frac{(-1)^k \log(k+1)}{(k+1)^{s-1}} - \frac{1}{(s-1)^2}\sum_{n=0}^{\infty}\frac{1}{n+1}\sum_{k=0}^{n}\binom{n}{k}\frac{(-1)^k}{(k+1)^{s-1}}$$

(4.3.239ii) $$= -\frac{1}{s-1}\sum_{n=0}^{\infty}\frac{1}{n+1}\sum_{k=0}^{n}\binom{n}{k}\frac{(-1)^k \log(k+1)}{(k+1)^{s-1}} - \frac{1}{s-1}\varsigma(s)$$

and with $s = 0$ we have using (F.6)

$$\varsigma'(0) = -\frac{1}{2}\log(2\pi) = \sum_{n=0}^{\infty}\frac{1}{n+1}\sum_{k=0}^{n}\binom{n}{k}(-1)^k (k+1)\log(k+1) - \frac{1}{2}$$

Hence we see that

$$\sum_{n=0}^{\infty}\frac{1}{n+1}\sum_{k=0}^{n}\binom{n}{k}(-1)^k (k+1)\log(k+1) = \frac{1}{2} - \frac{1}{2}\log(2\pi)$$

and therefore we see that

(4.3.240) $$\sum_{n=0}^{\infty}\frac{1}{n+1}\sum_{k=0}^{n}\binom{n}{k}(-1)^k k \log(k+1) = \gamma + \frac{1}{2} - \frac{1}{2}\log(2\pi)$$

With $s = 2$ we have

$$\varsigma'(2) = -\sum_{n=0}^{\infty}\frac{1}{n+1}\sum_{k=0}^{n}\binom{n}{k}\frac{(-1)^k \log(k+1)}{k+1} - \varsigma(2)$$

and using

$$\gamma_p(u) = -\frac{1}{p+1}\sum_{n=0}^{\infty}\frac{1}{n+1}\sum_{k=0}^{n}\binom{n}{k}(-1)^k \log^{p+1}(u+k)$$

we see that

(4.3.241) $$\gamma'_p(u) = -\sum_{n=0}^{\infty}\frac{1}{n+1}\sum_{k=0}^{n}\binom{n}{k}(-1)^k \frac{\log^p(u+k)}{u+k}$$

and

(4.3.242) $$\gamma'_p(1) = -\sum_{n=0}^{\infty}\frac{1}{n+1}\sum_{k=0}^{n}\binom{n}{k}(-1)^k \frac{\log^p(1+k)}{1+k}$$

Hence we have shown that

(4.3.243) $\quad \varsigma'(2) = \gamma_1'(1) - \varsigma(2)$

which we have already seen in a more general form in (4.3.223b). Using (F.7)

$$\varsigma'(-1) = \frac{1}{12}(1 - \gamma - \log 2\pi) + \frac{1}{2\pi^2}\varsigma'(2)$$

we see that

$$\varsigma'(2) = 2\pi^2 \varsigma'(-1) - \varsigma(2)(1 - \gamma - \log 2\pi)$$

and hence we have

(4.3.244) $\quad \gamma_1'(1) = 2\pi^2 \varsigma'(-1) + \varsigma(2)(\gamma + \log 2\pi)$

With $s = -2m$ in (4.3.239ii) we have

$$\varsigma'(-2m) = \frac{1}{2m+1}\sum_{n=0}^{\infty}\frac{1}{n+1}\sum_{k=0}^{n}\binom{n}{k}(-1)^k (k+1)^{2m+1}\log(k+1) + \frac{1}{2m+1}\varsigma(-2m)$$

$$= \frac{1}{2m+1}\sum_{n=0}^{\infty}\frac{1}{n+1}\sum_{k=0}^{n}\binom{n}{k}(-1)^k (k+1)^{2m+1}\log(k+1)$$

Therefore, using (F.8a)

$$\varsigma'(-2m) = (-1)^m \frac{(2m)!}{2(2\pi)^{2m}}\varsigma(2m+1)$$

we obtain (which we also saw in Volume II(a))

(4.3.245) $\quad \varsigma(2m+1) = (-1)^m \frac{2(2\pi)^{2m}}{(2m+1)!}\sum_{n=0}^{\infty}\frac{1}{n+1}\sum_{k=0}^{n}\binom{n}{k}(-1)^k (k+1)^{2m+1}\log(k+1)$

and in particular with $m = 1$ we have

(4.3.246) $\quad \varsigma(3) = -\frac{4\pi^2}{3}\sum_{n=0}^{\infty}\frac{1}{n+1}\sum_{k=0}^{n}\binom{n}{k}(-1)^k (k+1)^3 \log(k+1)$

We may therefore deduce a formula for $\sum_{n=0}^{\infty} \frac{1}{n+1} \sum_{k=0}^{n} \binom{n}{k} (-1)^k k^3 \log(k+1)$ involving $\varsigma(3)$.

From (4.3.107a) we have

$$3\varsigma'(-2, x) - \varsigma(-2, x) = \sum_{k=0}^{\infty} \frac{1}{k+1} \sum_{j=0}^{k} \binom{k}{j} (-1)^j (x+j)^3 \log(x+j)$$

and letting $x = 1$ this becomes

$$3\varsigma'(-2) = \sum_{k=0}^{\infty} \frac{1}{k+1} \sum_{j=0}^{k} \binom{k}{j} (-1)^j (1+j)^3 \log(1+j)$$

From (F.8b) we have $\varsigma'(-2) = -\frac{\varsigma(3)}{4\pi^2}$ and hence we recover (4.3.246).

With $s = 1 - 2m$ in (4.3.239ii) we have

$$\varsigma'(1-2m) = \frac{1}{2m} \sum_{n=0}^{\infty} \frac{1}{n+1} \sum_{k=0}^{n} \binom{n}{k} (-1)^k (k+1)^{2m} \log(k+1) + \frac{1}{2m} \varsigma(1-2m)$$

and using (F.12b) $\varsigma(1-2m) = -\frac{B_{2m}}{2m}$ this becomes

$$\varsigma'(1-2m) = \frac{1}{2m} \sum_{n=0}^{\infty} \frac{1}{n+1} \sum_{k=0}^{n} \binom{n}{k} (-1)^k (k+1)^{2m} \log(k+1) - \frac{B_{2m}}{4m^2}$$

In this regard, we also recall (F.8)

$$2m\varsigma'(1-2m) - \left[ H_{2m-1}^{(1)} - \gamma - \log 2\pi \right] B_{2m} = \frac{\varsigma'(2m)}{\varsigma(2m)} B_{2m}$$

and we therefore obtain

$$\sum_{n=0}^{\infty} \frac{1}{n+1} \sum_{k=0}^{n} \binom{n}{k} (-1)^k (k+1)^{2m} \log(k+1) = \left[ H_{2m-1}^{(1)} - \gamma - \log 2\pi \right] B_{2m} + \frac{\varsigma'(2m)}{\varsigma(2m)} B_{2m} - \frac{B_{2m}}{2m}$$

A further differentiation results in

$$\varsigma''(s) = \frac{1}{s-1} \sum_{n=0}^{\infty} \frac{1}{n+1} \sum_{k=0}^{n} \binom{n}{k} \frac{(-1)^k \log^2(k+1)}{(k+1)^{s-1}} + \frac{1}{(s-1)^2} \sum_{n=0}^{\infty} \frac{1}{n+1} \sum_{k=0}^{n} \binom{n}{k} \frac{(-1)^k \log(k+1)}{(k+1)^{s-1}}$$

$$-\frac{1}{s-1}\varsigma'(s)+\frac{1}{(s-1)^2}\varsigma(s)$$

$$=\frac{1}{s-1}\sum_{n=0}^{\infty}\frac{1}{n+1}\sum_{k=0}^{n}\binom{n}{k}\frac{(-1)^k \log^2(k+1)}{(k+1)^{s-1}}-\frac{2}{s-1}\varsigma'(s)$$

and thus we get

(4.3.246a)

$$\varsigma''(0)=-\sum_{n=0}^{\infty}\frac{1}{n+1}\sum_{k=0}^{n}\binom{n}{k}(-1)^k(k+1)\log^2(k+1)-\log(2\pi)$$

$$=-\sum_{n=0}^{\infty}\frac{1}{n+1}\sum_{k=0}^{n}\binom{n}{k}(-1)^k k \log^2(k+1)+2\gamma_1-\log(2\pi)$$

We also note from Apostol's paper [14aa] that

(4.3.247) $\qquad \varsigma''(0)=\gamma_1+\frac{1}{2}\gamma^2-\frac{1}{24}\pi^2-\frac{1}{2}\log^2(2\pi)$

and we then obtain

(4.3.248) $\quad \sum_{n=0}^{\infty}\frac{1}{n+1}\sum_{k=0}^{n}\binom{n}{k}(-1)^k k \log^2(k+1)=\gamma_1-\log(2\pi)-\frac{1}{2}\gamma^2+\frac{1}{24}\pi^2+\frac{1}{2}\log^2(2\pi)$

With $s=2$ we have

$$\varsigma''(2)=\sum_{n=0}^{\infty}\frac{1}{n+1}\sum_{k=0}^{n}\binom{n}{k}\frac{(-1)^k \log^2(k+1)}{k+1}-2\varsigma'(2)$$

$$=-\gamma_2'(1)-2\varsigma'(2)$$

and therefore we have

(4.3.249) $\qquad \varsigma''(2)=-\gamma_2'(1)-2\gamma_1'(1)+2\varsigma(2)$

## A CONNECTION WITH LOGARITHMIC INTEGRALS

Using the integral definition of the gamma function for $s>0$

$$\Gamma(s) = \int_0^\infty e^{-x} x^{s-1} \, dt$$

and making the substitution $x = ut$ (where $u > 0$) we obtain for $u, s > 0$

$$\int_0^\infty e^{-ut} t^{s-1} \, dt = \frac{\Gamma(s)}{u^s}$$

and using the binomial theorem we get for $n \geq 0$

(4.3.250) $$\int_0^\infty e^{-ut} \left(1 - e^{-t}\right)^n t^{s-1} \, dt = \Gamma(s) \sum_{k=0}^n \binom{n}{k} \frac{(-1)^k}{(k+u)^s}$$

Making the summation gives us

$$\sum_{n=0}^\infty \frac{1}{n+1} \int_0^\infty e^{-ut} \left(1 - e^{-t}\right)^n t^{s-1} \, dt = \sum_{n=0}^\infty \int_0^\infty e^{-(n+1)x} dx \int_0^\infty e^{-ut} \left(1 - e^{-t}\right)^n t^{s-1} \, dt$$

and this may be expressed as

$$\sum_{n=0}^\infty \int_0^\infty e^{-(n+1)x} dx \int_0^\infty e^{-ut} \left(1 - e^{-t}\right)^n t^{s-1} \, dt = \sum_{n=0}^\infty \int_0^\infty \int_0^\infty e^{-ut} e^{-x} \left[e^{-x}\left(1 - e^{-t}\right)\right]^n t^{s-1} \, dt \, dx$$

$$= \int_0^\infty \int_0^\infty \frac{e^{-ut} e^{-x} t^{s-1}}{1 - e^{-x}[1 - e^{-t}]} \, dt \, dx$$

We have

$$\int_0^\infty \frac{e^{-x}}{1 - e^{-x}[1 - e^{-t}]} \, dx = \frac{1}{1 - e^{-t}} \log(1 - e^{-x}[1 - e^{-t}]) \Big|_0^\infty = \frac{t}{1 - e^{-t}}$$

and therefore we see that

$$\int_0^\infty \int_0^\infty \frac{e^{-ut} e^{-x} t^{s-1}}{1 - e^{-x}[1 - e^{-t}]} \, dt \, dx = \int_0^\infty \frac{e^{-ut} t^s}{1 - e^{-t}} \, dt = \int_0^\infty \frac{e^{-(u-1)t} t^s}{e^t - 1} \, dt$$

We note from [126, p.92] that for $\operatorname{Re}(s) > 0$ and $\operatorname{Re}(u) > 0$

(4.3.251) $$\int_0^\infty \frac{e^{-(u-1)t} t^s}{e^t - 1} \, dt = \Gamma(s+1) \varsigma(s+1, u)$$

and with $u = 1$ we obtain the familiar formula

$$(4.3.252) \quad \int_0^\infty \frac{t^{s-1}}{e^t - 1} dt = \Gamma(s)\varsigma(s,1) = \Gamma(s)\varsigma(s)$$

We therefore have

$$(4.3.253) \quad \int_0^\infty \frac{e^{-(u-1)t} t^s}{e^t - 1} dt = \Gamma(s) \sum_{n=0}^\infty \frac{1}{n+1} \sum_{k=0}^n \binom{n}{k} \frac{(-1)^k}{(k+u)^s}$$

This gives us another proof of the Hasse identity for $\operatorname{Re}(s) > 0$

$$\varsigma(s+1, u) = \frac{1}{s} \sum_{n=0}^\infty \frac{1}{n+1} \sum_{k=0}^n \binom{n}{k} \frac{(-1)^k}{(k+u)^s}$$

Integrating (4.3.251) with respect to $u$ we see that

$$\int_1^x du \int_0^\infty \frac{e^{-(u-1)t} t^s}{e^t - 1} dt = \int_0^\infty \frac{[1 - e^{-(x-1)t}] t^{s-1}}{e^t - 1} = \Gamma(s)[\varsigma(s) - \varsigma(s, x)]$$

We also note that

$$\int_0^\infty \frac{e^{-(u-1)t} t^s}{e^t - 1} dt = \Gamma(s+1)\varsigma(s+1, u) = \Gamma(s) s \varsigma(s+1, u) = -\Gamma(s) \frac{\partial}{\partial u} \varsigma(s, u)$$

$$= \Gamma(s) \sum_{n=0}^\infty \frac{1}{n+1} \sum_{k=0}^n \binom{n}{k} \frac{(-1)^k}{(k+u)^s}$$

We have

$$\int_0^\infty \frac{e^{-(u-1)t} t^s}{\Gamma(s)[e^t - 1]} dt = \sum_{n=0}^\infty \frac{1}{n+1} \sum_{k=0}^n \binom{n}{k} \frac{(-1)^k}{(k+u)^s}$$

and differentiating this with respect to $s$ results in

$$\int_0^\infty \frac{e^{-(u-1)t} t^s [\log t - \psi(s)]}{\Gamma(s)[e^t - 1]} dt = -\sum_{n=0}^\infty \frac{1}{n+1} \sum_{k=0}^n \binom{n}{k} \frac{(-1)^k \log(k+u)}{(k+u)^s}$$

Hence we have

(4.3.254)

$$\int_0^\infty \frac{e^{-(u-1)t} t^s \log t}{e^t - 1} dt = \psi(s)\Gamma(s+1)\varsigma(s+1,u) - \Gamma(s)\sum_{n=0}^\infty \frac{1}{n+1}\sum_{k=0}^n \binom{n}{k} \frac{(-1)^k \log(k+u)}{(k+u)^s}$$

With $u = 1$ in (4.3.254) we get

$$\int_0^\infty \frac{t^s \log t}{e^t - 1} dt = \psi(s)\Gamma(s+1)\varsigma(s+1) - \Gamma(s)\sum_{n=0}^\infty \frac{1}{n+1}\sum_{k=0}^n \binom{n}{k} \frac{(-1)^k \log(k+1)}{(k+1)^s}$$

which is equivalent to

$$\int_0^\infty \frac{t^s \log t}{e^t - 1} dt = \frac{d}{ds}[\Gamma(s+1)\varsigma(s+1)]$$

With $s = 1$ in (4.3.254) we get

$$\int_0^\infty \frac{e^{-(u-1)t} t \log t}{e^t - 1} dt = -\gamma\varsigma(2,u) - \sum_{n=0}^\infty \frac{1}{n+1}\sum_{k=0}^n \binom{n}{k} \frac{(-1)^k \log(k+u)}{k+u}$$

and using (4.3.70) this becomes

$$\int_0^\infty \frac{e^{-(u-1)t} t \log t}{e^t - 1} dt = -\gamma\psi'(u) - \sum_{n=0}^\infty \frac{1}{n+1}\sum_{k=0}^n \binom{n}{k} \frac{(-1)^k \log(k+u)}{k+u}$$

Reference to (4.3.210) then shows that

(4.3.255) $\int_0^\infty \frac{e^{-(u-1)t} t \log t}{e^t - 1} dt = \gamma_1'(u) - \gamma\psi'(u) = \gamma_1'(u) + \gamma\gamma_0'(u)$

Integrating this with respect to $u$ we get

$$\int_1^x du \int_0^\infty \frac{e^{-(u-1)t} t \log t}{e^t - 1} dt = \int_0^\infty \frac{[1 - e^{-(x-1)t}]\log t}{e^t - 1} dt = \gamma_1(x) - \gamma_1 - \gamma\psi(x) - \gamma^2$$

The substitution $y = e^{-t}$ gives us

(4.3.256) $\quad \int_0^1 \frac{1 - y^{x-1}}{1 - y} \log\log(1/y) \, dy = \gamma_1(x) - \gamma_1 - \gamma\psi(x) - \gamma^2$

and with $x = 1/2$ this becomes

(4.3.257)

$$\int_0^1 \frac{1 - 1/\sqrt{y}}{1-y} \log\log(1/y)\,dy = -\int_0^1 \frac{1/\sqrt{y}}{1+\sqrt{y}} \log\log(1/y)\,dy = \gamma_1\left(\frac{1}{2}\right) - \gamma_1 - \gamma\psi\left(\frac{1}{2}\right) - \gamma^2$$

We now recall Adamchik's result (C.57) in Volume VI (see reference [2a] and also Roach's paper [110b])

(4.3.258)

$$\int_0^1 \frac{x^{p-1}}{1+x^n} \log\log\left(\frac{1}{x}\right) dx = \frac{\gamma + \log(2n)}{2n}\left[\psi\left(\frac{p}{2n}\right) - \psi\left(\frac{n+p}{2n}\right)\right] + \frac{1}{2n}\left[\varsigma'\left(1,\frac{p}{2n}\right) - \varsigma'\left(1,\frac{n+p}{2n}\right)\right]$$

where, notwithstanding the suggestive notation, neither $p$ nor $n$ need be integers.

With $p = n$ we get

(4.3.259) $\int_0^1 \frac{x^{n-1}}{1+x^n} \log\log\left(\frac{1}{x}\right) dx = \frac{\gamma + \log(2n)}{2n}\left[\psi\left(\frac{1}{2}\right) - \psi(1)\right] + \frac{1}{2n}\left[\varsigma'\left(1,\frac{1}{2}\right) - \varsigma'(1,1)\right]$

From (4.3.228c), and also from Adamchik's paper [2a], we see that

(4.3.260) $\varsigma'\left(1,\frac{1}{2}\right) - \varsigma'(1,1) = \log^2 2 + 2\gamma\log 2$

and hence, using $\psi\left(\frac{1}{2}\right) = -\gamma - 2\log 2$, we have

(4.3.261) $\int_0^1 \frac{x^{n-1}}{1+x^n} \log\log\left(\frac{1}{x}\right) dx = -\frac{1}{2n}\left[\log^2 2 + 2\log 2 \log n\right]$

With $p = n = 1/2$ we get

(4.3.262) $\int_0^1 \frac{1/\sqrt{x}}{1+\sqrt{x}} \log\log\left(\frac{1}{x}\right) dx = \log^2 2$

We then see from (4.3.257) and (4.3.262) that

(4.3.263) $$\gamma_1\left(\frac{1}{2}\right) = \gamma_1 - \log^2 2 - 2\gamma \log 2$$

and this is the same as the result obtained by Coffey [45c] in (4.3.233).

From (4.3.256) we have

(4.3.264) $$\int_0^1 \frac{1-x^{-3/4}}{1-x} \log\log(1/x)\, dx = \gamma_1\left(\frac{1}{4}\right) - \gamma_1 - \gamma\psi\left(\frac{1}{4}\right) - \gamma^2$$

and using [126, p.22]

(4.3.265) $$\psi\left(\frac{1}{4}\right) = -\gamma - \frac{1}{2}\pi - 3\log 2$$

we get

(4.3.266) $$\int_0^1 \frac{1-x^{-3/4}}{1-x} \log\log(1/x)\, dx = \gamma_1\left(\frac{1}{4}\right) - \gamma_1 + \frac{1}{2}\gamma\pi + 3\gamma \log 2$$

We easily see that

$$\frac{1-x^{-3/4}}{1-x} = -\frac{x^{-3/4}\left(x^{3/4}-1\right)}{x-1} = -\frac{x^{-3/4}\left(x^{3/4}-1\right)}{\left(x^{1/2}-1\right)\left(x^{1/2}+1\right)} = -\frac{x^{-3/4}}{\left(x^{1/2}+1\right)}\left[x^{1/4} + \frac{x^{1/4}-1}{\left(x^{1/2}-1\right)}\right]$$

$$= -\frac{x^{-3/4}}{\left(x^{1/2}+1\right)}\left[x^{1/4} + \frac{1}{\left(x^{1/4}+1\right)}\right] = -\frac{x^{-1/2}}{\left(x^{1/2}+1\right)} - \frac{x^{-3/4}}{\left(x^{1/2}+1\right)\left(x^{1/4}+1\right)}$$

Using partial fractions we see that

$$\frac{x^{-3/4}}{\left(x^{1/2}+1\right)\left(x^{1/4}+1\right)} = \frac{x^{-3/4}}{2\left(x^{1/4}+1\right)} - \frac{x^{-1/2}}{2\left(x^{1/2}+1\right)} + \frac{x^{-3/4}}{2\left(x^{1/2}+1\right)}$$

and we therefore obtain

$$\frac{1-x^{-3/4}}{1-x} = -\frac{x^{-3/4}}{2(x^{1/4}+1)} - \frac{x^{-1/2}}{2(x^{1/2}+1)} - \frac{x^{-3/4}}{2(x^{1/2}+1)}$$

There is a more straightforward derivation using $y = x^{1/4}$ but, for posterity, I have left my slightly more complex struggle intact!

Letting $p = 1/4$ and $n = 1/2$ in (4.3.258) we get

$$\int_0^1 \frac{x^{-3/4}}{1+x^{1/2}} \log\log\left(\frac{1}{x}\right) dx = \gamma\left[\psi\left(\frac{1}{4}\right) - \psi\left(\frac{3}{4}\right)\right] + \left[\varsigma'\left(1,\frac{1}{4}\right) - \varsigma'\left(1,\frac{3}{4}\right)\right]$$

and with $p = 1/2$ and $n = 1/2$ in (4.3.258) we have already seen that

$$\int_0^1 \frac{x^{-1/2}}{1+x^{1/2}} \log\log\left(\frac{1}{x}\right) dx = \log^2 2$$

With $p = n = 1/4$ we get

$$\int_0^1 \frac{x^{-3/4}}{1+x^{1/4}} \log\log\left(\frac{1}{x}\right) dx = 6\log^2 2$$

We therefore obtain

$$\int_0^1 \frac{1-x^{-3/4}}{1-x} \log\log(1/x)\, dx = -\frac{7}{2}\log^2 2 - \frac{1}{2}\gamma\left[\psi\left(\frac{1}{4}\right) - \psi\left(\frac{3}{4}\right)\right] - \frac{1}{2}\left[\varsigma'\left(1,\frac{1}{4}\right) - \varsigma'\left(1,\frac{3}{4}\right)\right]$$

We then refer to (4.3.228b) and see that

$$\varsigma'\left(1,\frac{1}{4}\right) - \varsigma'\left(1,\frac{3}{4}\right) = -\left[\gamma_1\left(\frac{1}{4}\right) - \gamma_1\left(\frac{3}{4}\right)\right]$$

and using [126, p.22]

$$\psi\left(\frac{1}{4}\right) = -\gamma - \frac{1}{2}\pi - 3\log 2 \qquad \psi\left(\frac{3}{4}\right) = -\gamma + \frac{1}{2}\pi - 3\log 2$$

we get

(4.3.267) $$\int_0^1 \frac{1-x^{-3/4}}{1-x} \log\log(1/x)\, dx = -\frac{7}{2}\log^2 2 + \frac{1}{2}\gamma\pi + \frac{1}{2}\left[\gamma_1\left(\frac{1}{4}\right) - \gamma_1\left(\frac{3}{4}\right)\right]$$

Using (4.3.233a)

$$\gamma_1\left(\frac{1}{4}\right) - \gamma_1\left(\frac{3}{4}\right) = -\pi\left[\gamma + 4\log 2 + 3\log \pi - 4\log \Gamma\left(\frac{1}{4}\right)\right]$$

this then becomes

$$\int_0^1 \frac{1-x^{-3/4}}{1-x} \log\log(1/x)\, dx = -\frac{7}{2}\log^2 2 + \frac{1}{2}\gamma\pi - \frac{1}{2}\pi\left[\gamma + 4\log 2 + 3\log \pi - 4\log \Gamma\left(\frac{1}{4}\right)\right]$$

and equating this with (4.3.266) we get

(4.3.268)

$$\gamma_1\left(\frac{1}{4}\right) - \gamma_1 + \frac{1}{2}\gamma\pi + 3\gamma\log 2 = -\frac{7}{2}\log^2 2 + \frac{1}{2}\gamma\pi - \frac{1}{2}\pi\left[\gamma + 4\log 2 + 3\log \pi - 4\log \Gamma\left(\frac{1}{4}\right)\right]$$

This is in fact another proof of equation (4.3.233h).

$\square$

Differentiating (4.3.251) results in

(4.3.269) $$\int_0^\infty \frac{e^{-(u-1)t} t^s \log t}{e^t - 1}\, dt = \Gamma(s+1)\varsigma'(s+1,u) + \Gamma'(s+1)\varsigma(s+1,u)$$

and with $s = 1$ we obtain

(4.3.270) $$\int_0^\infty \frac{e^{-(u-1)t} t \log t}{e^t - 1}\, dt = \varsigma'(2,u) + \Gamma'(2)\varsigma(2,u)$$

Hence we have using the derivative of (4.3.253) combined with (4.3.254)

(4.3.271) $$\varsigma'(2,u) + \Gamma'(2)\varsigma(2,u) = \gamma_1'(u) - \gamma\psi'(u) = \gamma_1'(u) + \gamma\gamma_0'(u)$$

and integration over the interval $[1, x]$ produces

$$\int_1^x [\varsigma'(2,u) + \psi(2)\varsigma(2,u)]\, du = \sum_{n=0}^\infty \frac{1-\psi(2)}{n+u} + \sum_{n=0}^\infty \frac{\log(n+u)}{n+u}\Big|_1^x$$

$$= [1-\psi(2)]\sum_{n=0}^{\infty}\left[\frac{1}{n+x}-\frac{1}{n+1}\right]+\sum_{n=0}^{\infty}\left[\frac{\log(n+x)}{n+x}-\frac{\log(n+1)}{n+1}\right]$$

$$= -[1-\psi(2)][\gamma+\psi(x)]+\gamma_1(x)-\gamma_1$$

(4.3.272)
$$= \gamma_1(x)-\gamma_1-\gamma\psi(x)-\gamma^2$$

which we have seen previously.

Letting $u=1$ in (4.3.271) gives us (see also (4.4.42bi))

(4.3.273)
$$\int_0^\infty \frac{t\log t}{e^t-1}dt = \varsigma'(2)+(1-\gamma)\varsigma(2)$$

which may also be represented differently using (F.7)

$$\varsigma'(-1) = \frac{1}{12}(1-\gamma-\log 2\pi)+\frac{1}{2\pi^2}\varsigma'(2)$$

Differentiating (4.3.253) with respect to $u$ gives us

$$\int_0^\infty \frac{e^{-(u-1)t}t^{s+1}}{\Gamma(s)[e^t-1]}dt = s\sum_{n=0}^\infty \frac{1}{n+1}\sum_{k=0}^n \binom{n}{k}\frac{(-1)^k}{(k+u)^{s+1}}$$

which is of course equivalent to letting $s \to s+1$ in (4.3.253).

Alternatively, differentiating (4.3.253) with respect to $s$ results in

(4.3.274)
$$\int_0^\infty \frac{e^{-(u-1)t}t^s(\log t-\psi(s))}{e^t-1}dt = -\Gamma(s)\sum_{n=0}^\infty \frac{1}{n+1}\sum_{k=0}^n \binom{n}{k}\frac{(-1)^k\log(k+u)}{(k+u)^s}$$

and a further differentiation gives us

(4.3.275)
$$\int_0^\infty \frac{e^{-(u-1)t}t^s\{[\log t-\psi(s)]^2-\psi'(s)\}}{e^t-1}dt = \Gamma(s)\sum_{n=0}^\infty \frac{1}{n+1}\sum_{k=0}^n \binom{n}{k}\frac{(-1)^k\log^2(k+u)}{(k+u)^s}$$

Completing the summation of (4.3.250) starting at $n=1$ results in

$$\sum_{n=1}^\infty \frac{1}{n+1}\int_0^\infty e^{-ut}\left(1-e^{-t}\right)^n t^{s-1}dt = \sum_{n=1}^\infty \int_0^\infty e^{-(n+1)x}dx \int_0^\infty e^{-ut}\left(1-e^{-t}\right)^n t^{s-1}dt$$

$$= \sum_{n=1}^{\infty} \int_0^{\infty}\int_0^{\infty} e^{-ut} e^{-x} \left[ e^{-x}\left(1-e^{-t}\right)\right]^n t^{s-1}\, dt\, dx$$

$$= \int_0^{\infty}\int_0^{\infty} \frac{e^{-ut} e^{-2x}\left(1-e^{-t}\right) t^{s-1}}{1-e^{-x}[1-e^{-t}]}\, dt\, dx$$

We have with the obvious substitution $y = e^{-x}$

$$\int_0^{\infty} \frac{e^{-2x}}{1-e^{-x}[1-e^{-t}]}\, dx = \int_0^1 \frac{y}{1-[1-e^{-t}]y}\, dy = -\frac{y}{1-e^{-t}} - \frac{\log(1-[1-e^{-t}]y)}{[1-e^{-t}]^2}\bigg|_0^1$$

$$= -\frac{1}{1-e^{-t}} + \frac{t}{[1-e^{-t}]^2}$$

and we therefore see that

$$\int_0^{\infty}\int_0^{\infty} \frac{e^{-ut} e^{-2x}\left(1-e^{-t}\right) t^{s-1}}{1-e^{-x}[1-e^{-t}]}\, dt\, dx = \int_0^{\infty} \frac{e^{-ut} t^s}{[1-e^{-t}]^2}\, dt - \int_0^{\infty} \frac{e^{-ut} t^{s-1}}{1-e^{-t}}\, dt$$

$$= \int_0^{\infty} \frac{e^{-(u-2)t} t^s}{[e^t-1]^2}\, dt - \int_0^{\infty} \frac{e^{-(u-1)t} t^{s-1}}{e^t-1}\, dt$$

We now let $t \to \alpha t$ and obtain

$$I = \int_0^{\infty} \frac{e^{-(u-2)t} t^{s-1}}{e^t - 1}\, dt = \int_0^{\infty} \frac{e^{-(u-2)\alpha t} \alpha^s t^{s-1}}{e^{\alpha t} - 1}\, dt$$

Then we have

$$\frac{dI}{d\alpha} = \int_0^{\infty} \frac{t^{s-1}\left[(e^{\alpha t}-1)e^{-(u-2)\alpha t}\left\{s\alpha^{s-1} - (u-2)t\alpha^s\right\}\right] - e^{-(u-2)\alpha t}\alpha^s t^s e^{\alpha t}}{[e^{\alpha t}-1]^2}\, dt = 0$$

and with $\alpha = 1$ this becomes

$$\int_0^{\infty} \frac{t^{s-1}\left[(e^t-1)e^{-(u-2)t}\left\{s-(u-2)t\right\}\right] - e^{-(u-2)t} t^s e^t}{[e^t-1]^2}\, dt = 0$$

and we therefore see that

$$\int\limits_0^\infty \frac{e^{-(u-2)t}t^s e^t}{[e^t-1]^2}\,dt = s\int\limits_0^\infty \frac{t^{s-1}e^{-(u-2)t}}{[e^t-1]}\,dt - (u-2)\int\limits_0^\infty \frac{t^s e^{-(u-2)t}}{e^t-1}\,dt$$

With $u=2$ this becomes for $\mathrm{Re}(s)>1$

(4.3.276) $$\int\limits_0^\infty \frac{t^s e^t}{[e^t-1]^2}\,dt = s\int\limits_0^\infty \frac{t^{s-1}}{e^t-1}\,dt = s\Gamma(s)\varsigma(s)$$

which is reported in [126, p.103].

This may also be shown using integration by parts: we have

$$\int\limits_0^\infty \frac{e^t t^s}{[e^t-1]^2}\,dt = -\frac{t^s}{e^t-1}\bigg|_0^\infty + s\int\limits_0^\infty \frac{t^{s-1}}{e^t-1}\,dt$$

Application of L'Hôpital's rule shows us that

$$\lim_{t\to 0}\frac{t^s}{e^t-1}=0 \text{ if } s>0 \quad \text{and} \quad \lim_{t\to\infty}\frac{t^s}{e^t-1}=0$$

and therefore, as noted above, for $s>0$ we have

$$\int\limits_0^\infty \frac{e^t t^s}{[e^t-1]^2}\,dt = s\int\limits_0^\infty \frac{t^{s-1}}{e^t-1}\,dt$$

We have therefore shown that for $u,s>0$

$$\int\limits_0^\infty\int\limits_0^\infty \frac{e^{-ut}e^{-2x}(1-e^{-t})t^{s-1}}{1-e^{-x}[1-e^{-t}]}\,dt\,dx = s\int\limits_0^\infty \frac{t^{s-1}e^{-(u-2)t}}{e^t-1}\,dt - (u-2)\int\limits_0^\infty \frac{t^s e^{-(u-2)t}}{e^t-1}\,dt - \int\limits_0^\infty \frac{e^{-(u-1)t}t^{s-1}}{e^t-1}\,dt$$

$$=\Gamma(s)\sum_{n=1}^\infty \frac{1}{n+1}\sum_{k=0}^n \binom{n}{k}\frac{(-1)^k}{(k+u)^s}$$

and, using (4.3.251), for $u>1$, $s>1$ we have

$$\int\limits_0^\infty\int\limits_0^\infty \frac{e^{-ut}e^{-2x}(1-e^{-t})t^{s-1}}{1-e^{-x}[1-e^{-t}]}\,dt\,dx = \Gamma(s+1)\varsigma(s,u-1) - (u-2)\Gamma(s+1)\varsigma(s,u-1) - \Gamma(s)\varsigma(s,u)$$

$$= \Gamma(s) \sum_{n=1}^{\infty} \frac{1}{n+1} \sum_{k=0}^{n} \binom{n}{k} \frac{(-1)^k}{(k+u)^s}$$

Hence we have for $u > 1$, $s > 1$

$$\sum_{n=1}^{\infty} \frac{1}{n+1} \sum_{k=0}^{n} \binom{n}{k} \frac{(-1)^k}{(k+u)^s} = s\varsigma(s, u-1) - (u-2)s\varsigma(s, u-1) - \varsigma(s, u)$$

and with $u = 2$ we obtain

$$\sum_{n=1}^{\infty} \frac{1}{n+1} \sum_{k=0}^{n} \binom{n}{k} \frac{(-1)^k}{(k+2)^s} = s\varsigma(s, 1) - \varsigma(s, 2) = s\varsigma(s) - \varsigma(s, 2)$$

Differentiating (4.3.251) results in

$$\int_0^\infty \frac{e^{-(u-1)t} t^s \log^p t}{e^t - 1} = \frac{\partial^p}{\partial s^p} [\Gamma(s+1)\varsigma(s+1, u)] = (-1)^p \Gamma(s) \sum_{n=0}^{\infty} \frac{1}{n+1} \sum_{k=0}^{n} \binom{n}{k} \frac{(-1)^k \log^p (k+u)}{(k+u)^s}$$

and for convenience we may let $s \to s-1$ and write the above equation as

$$\frac{1}{\Gamma(s-1)} \int_0^\infty \frac{e^{-(u-1)t} t^{s-1} \log^p t}{e^t - 1} = \frac{1}{\Gamma(s-1)} \frac{\partial^p}{\partial s^p} [\Gamma(s)\varsigma(s, u)] = (-1)^p \sum_{n=0}^{\infty} \frac{1}{n+1} \sum_{k=0}^{n} \binom{n}{k} \frac{(-1)^k \log^p (k+u)}{(k+u)^{s-1}}$$

$\square$

In the book "Concrete Mathematics" [75, p.252] it is an exercise to prove that

(4.3.280) $$\Gamma(x+1) = \sum_{n=0}^{\infty} s_n \binom{x}{n}$$

where $s_0 = s_1 = 0$ and for $n \geq 2$ we have

$$s_n = (-1)^n \int_0^\infty e^{-t} (1 - e^{-t})^{n-1} t^{-1} dt$$

$$= (-1)^n \sum_{k=0}^{n} \binom{n-1}{k} (-1)^{k+1} \log(1+k)$$

and we note the connection with (4.3.250) which is valid for $s > 0$.

$$\int_0^\infty e^{-ut}\left(1-e^{-t}\right)^n t^{s-1}\, dt = \Gamma(s)\sum_{k=0}^{n}\binom{n}{k}\frac{(-1)^k}{(k+u)^s}$$

Therefore, since $\binom{n}{n+1}=0$, for $n\geq 1$ we have

$$\int_0^\infty e^{-t}\left(1-e^{-t}\right)^n t^{-1}\, dt = \sum_{k=0}^{n}\binom{n}{k}(-1)^{k+1}\log(1+k)$$

$\square$

In 1999 Coppo [46b] showed that there exists a polynomial $p_n(z)$ such that

(4.3.290) $$x^n = \int_0^\infty p_n(x-\log z)e^{-z}\, dz$$

and the following related analysis is extracted from Brede's dissertation [33b].

With $p_n(z) = \sum_{k=0}^{n} a_k z^k$ we have

$$\int_0^\infty p_n(x-\log z)e^{-z}\, dz = \int_0^\infty \sum_{k=0}^{n} a_k (x-\log z)^k e^{-z}\, dz$$

$$= \sum_{k=0}^{n} a_k \int_0^\infty \sum_{l=0}^{k}\binom{k}{l} x^l (-\log z)^{k-l} e^{-z}\, dz$$

$$= \sum_{k=0}^{n} a_k \sum_{l=0}^{k}\binom{k}{l} (-1)^{k-l}\Gamma^{(k-l)}(1)\, x^l$$

Hence we have

$$\int_0^\infty p_n(x-\log z)e^{-z}\, dz = \sum_{k=0}^{n}\left[\sum_{l=0}^{n-k}(-1)^l \binom{k+l}{l}\Gamma^{(l)}(1)\, a_{k+l}\right] x^k$$

and we then select the coefficients $a_k$ such that

$$\sum_{l=0}^{n-k}(-1)^l \binom{k+l}{l}\Gamma^{(l)}(1)\, a_{k+l} = \begin{cases} 1 & \text{for } k=n \\ 0 & \text{for } k \text{ less than } n \end{cases}$$

resulting in

$$x^n = \int_0^\infty p_n(x - \log z) e^{-z} dz$$

For example we have

(4.3.291) $p_0(z) = 1$

$$p_1(z) = z - \gamma$$

$$p_2(z) = z^2 - 2\gamma z + \gamma^2 - \varsigma(2)$$

Letting $x \to \log x$ we obtain

$$\log^n x = \int_0^\infty p_n(\log x - \log z) e^{-z} dz$$

and, with the substitution $z/x = -\log t$, we have

$$\frac{\log^n x}{x} = \int_0^1 p_n(-\log \log(1/t)) t^{x-1} dt$$

Referring to (4.3.225a)

$$\gamma_n = \sum_{k=1}^\infty \left[ \frac{\log^n k}{k} - \frac{\log^n(k+1) - \log^n k}{n+1} \right] = \sum_{k=1}^\infty \left[ \frac{\log^n k}{k} - \int_k^{k+1} \frac{\log^n t}{t} dt \right]$$

$$= \lim_{N \to \infty} \left[ \sum_{k=1}^N \frac{\log^n k}{k} - \int_1^{N+1} \frac{\log^n t}{t} dt \right]$$

we then see that

$$\gamma_n = \lim_{N \to \infty} \left[ \sum_{k=1}^N \int_0^1 p_n(-\log \log(1/t)) t^{k-1} dt - \int_1^{N+1} \int_0^1 p_n(-\log \log(1/t)) t^{x-1} dt\, dx \right]$$

$$= \lim_{N \to \infty} \int_0^1 p_n(-\log \log(1/t)) \left[ \sum_{k=0}^{N-1} t^k - \int_1^{N+1} t^x dx \right] dt$$

$$= \lim_{N \to \infty} \int_0^1 p_n(-\log \log(1/t)) \left[ \frac{t^N - 1}{t - 1} - \frac{t^N - 1}{\log t} \right] dt$$

Hence, as shown more completely in Brede's dissertation, we obtain as $N \to \infty$

(4.3.292) $\quad \gamma_n = \int_0^1 p_n(-\log \log(1/t)) \left[ \frac{1}{\log t} - \frac{1}{t - 1} \right] dt$

Since $\dfrac{1}{\log t} - \dfrac{1}{t-1} = \int_0^1 \dfrac{1 - u^t}{1 - u} du$ we may write this as a double integral

$$\gamma_n = \int_0^1 \int_0^1 p_n(-\log \log(1/t)) \frac{1 - u^t}{1 - u} du\, dt$$

From (4.3.68d) we saw that for $n \geq 1$

(4.3.293) $\quad \displaystyle\int_0^1 \frac{y^{u-1}(1-y)^n}{\log^r y} dy = \frac{1}{\Gamma(r)} \sum_{j=0}^{n} \binom{n}{j} (-1)^j (u+j)^{r-1} \log(u+j)$

and we have the summation

$$\sum_{n=1}^{\infty} \frac{1}{n+1} \int_0^1 \frac{y^{u-1}(1-y)^n}{\log^r y} dy = \frac{1}{\Gamma(r)} \sum_{n=1}^{\infty} \frac{1}{n+1} \sum_{j=0}^{n} \binom{n}{j} (-1)^j (u+j)^{r-1} \log(u+j)$$

We have

$$\sum_{n=1}^{\infty} \frac{(1-y)^n}{n+1} = -\frac{1}{1-y}[\log y + 1 - y]$$

and therefore we obtain

(4.3.294) $\quad \displaystyle\int_0^1 \left[ \frac{y^{u-1}}{(1-y)\log^{r-1} y} + \frac{y^{u-1}}{\log^r y} \right] dy = -\frac{1}{\Gamma(r)} \sum_{n=1}^{\infty} \frac{1}{n+1} \sum_{j=0}^{n} \binom{n}{j} (-1)^j (u+j)^{r-1} \log(u+j)$

Starting the summation at $n = 0$ gives us

$$\int_0^1 \left[ \frac{y^{u-1}}{(1-y)\log^{r-1} y} + \frac{y^{u-1}}{\log^r y} \right] dy = \frac{1}{\Gamma(r)} u^{r-1} \log u - \frac{1}{\Gamma(r)} \sum_{n=0}^{\infty} \frac{1}{n+1} \sum_{j=0}^{n} \binom{n}{j} (-1)^j (u+j)^{r-1} \log(u+j)$$

and reference to (4.3.107a) shows that this is related to $\varsigma'(2-r, u)$.

With $r=1$ we get

(4.3.295) $$\int_0^1 \left[\frac{y^{u-1}}{(1-y)} + \frac{y^{u-1}}{\log y}\right] dy = \log u - \sum_{n=0}^{\infty} \frac{1}{n+1} \sum_{j=0}^{n} \binom{n}{j} (-1)^j \log(u+j)$$

and using (4.3.74) this becomes

(4.3.296) $$\int_0^1 \left[\frac{y^{u-1}}{(1-y)} + \frac{y^{u-1}}{\log y}\right] dy = \log u - \psi(u)$$

With the substitution $y = 1/t$ in (4.3.293) we get

$$\int_0^1 \frac{y^{u-1}(1-y)^n}{\log^r y} dy = \int_1^{\infty} \frac{(1-1/t)^n}{t^{u+1} \log^r (1/t)} dt$$

and a summation gives us (starting at $n=1$)

$$\sum_{n=1}^{\infty} \frac{1}{n+1} \int_1^{\infty} \frac{(1-1/t)^n}{t^{u+1} \log^r (1/t)} dt = \int_1^{\infty} \left[\frac{t \log(1/t)}{1-t} - 1\right] \frac{1}{t^{u+1} \log^r (1/t)} dt$$

$$= \int_1^{\infty} \left[\frac{1}{t^u (1-t) \log^{r-1}(1/t)} - \frac{1}{t^{u+1} \log^r (1/t)}\right] dt$$

Differentiation of this with respect to $r$ results in

(4.3.297)
$$-\int_1^{\infty} \left[\frac{\log \log(1/t)}{t^u (1-t) \log^{r-1}(1/t)} - \frac{\log \log(1/t)}{t^{u+1} \log^r (1/t)}\right] dt = -\frac{\psi(r)}{\Gamma(r)} \sum_{n=1}^{\infty} \frac{1}{n+1} \sum_{j=0}^{n} \binom{n}{j} (-1)^j (u+j)^{r-1} \log(u+j)$$

$$+ \frac{1}{\Gamma(r)} \sum_{n=1}^{\infty} \frac{1}{n+1} \sum_{j=0}^{n} \binom{n}{j} (-1)^j (u+j)^{r-1} \log^2(u+j)$$

and for $r=1$ we get

(4.3.297a) $$-\int_1^{\infty} \left[\frac{\log \log(1/t)}{t^u (1-t)} - \frac{\log \log(1/t)}{t^{u+1} \log(1/t)}\right] dt = \gamma \sum_{n=1}^{\infty} \frac{1}{n+1} \sum_{j=0}^{n} \binom{n}{j} (-1)^j \log(u+j)$$

$$+ \sum_{n=1}^{\infty} \frac{1}{n+1} \sum_{j=0}^{n} \binom{n}{j} (-1)^j \log^2(u+j)$$

$$= \gamma \sum_{n=0}^{\infty} \frac{1}{n+1} \sum_{j=0}^{n} \binom{n}{j} (-1)^j \log(u+j)$$

$$+ \sum_{n=0}^{\infty} \frac{1}{n+1} \sum_{j=0}^{n} \binom{n}{j} (-1)^j \log^2(u+j)$$

$$-\gamma \log u - \log^2 u$$

Using (4.3.210) we see that this may be written in terms of the Stieltjes constants

(4.3.298) $\displaystyle\int_{1}^{\infty} \left[ \frac{1}{1-t} - \frac{1}{t \log(1/t)} \right] \frac{\log \log(1/t)}{t^u} dt = \gamma \gamma_0(u) + 2\gamma_1(u) + \gamma \log u + \log^2 u$

$$= 2\gamma_1(u) - \gamma \psi(u) + \gamma \log u + \log^2 u$$

and with $u=1$ we obtain

(4.3.299) $\displaystyle\int_{1}^{\infty} \left[ \frac{\log \log(1/t)}{t(1-t)} - \frac{\log \log(1/t)}{t^2 \log(1/t)} \right] dt = 2\gamma_1 + \gamma^2$

Further differentiations will result in integrals $\displaystyle\int_{1}^{\infty} \left[ \frac{\log^n \log(1/t)}{t^u (1-t)} - \frac{\log^n \log(1/t)}{t^{u+1} \log(1/t)} \right] dt$ involving the Stieltjes constants $\gamma_n(u)$.

Differentiation of (4.3.298) with respect to $u$ results in

(4.3.299a) $-\displaystyle\int_{1}^{\infty} \left[ \frac{\log t \log \log(1/t)}{t^u (1-t)} + \frac{\log \log(1/t)}{t^{u+1}} \right] dt = 2\gamma_1'(u) - \gamma \psi'(u)$

and with $u=1$ we get

(4.3.299b) $-\displaystyle\int_{1}^{\infty} \left[ \frac{\log t \log \log(1/t)}{t(1-t)} + \frac{\log \log(1/t)}{t^2} \right] dt = 2\gamma_1'(1) - \gamma \varsigma(2)$

Using (4.3.244)

$$\gamma_1'(1) = 2\pi^2 \varsigma'(-1) + \varsigma(2)(\gamma + \log 2\pi)$$

we may write this as

$$\text{(4.3.299c)} \quad -\int_1^\infty \left[\frac{\log t}{1-t} + \frac{1}{t}\right] \frac{\log\log(1/t)}{t}\, dt = 4\pi^2 \varsigma'(-1) + \varsigma(2)(\gamma + 2\log 2\pi)$$

## A HITCHHIKER'S GUIDE TO THE RIEMANN HYPOTHESIS

I was vaguely aware that the Stieltjes constants were relevant to the Riemann Hypothesis but I never really paid much attention to that aspect of analysis because I assumed that the mathematics was fiendishly difficult. However, having derived an interesting formula for $\gamma_p(u)$, I decided to dip into this exotic topic via Coffey's papers entitled "Toward verification of the Riemann Hypothesis: Application of the Li criterion" [45f] and "Polygamma theory, the Li/Keiper constants, and validity of the Riemann Hypothesis" [45g]. The following is a brief synopsis of some of the results contained in those papers, with some new material added.

Defining the Riemann xi function $\xi(s)$ as

$$\text{(4.3.300)} \quad \xi(s) = \frac{1}{2}s(s-1)\pi^{-s/2}\Gamma(s/2)\varsigma(s)$$

we see from the functional equation for $\varsigma(s)$ that (see Appendix F in Volume VI)

$$\text{(4.3.301)} \quad \xi(s) = \xi(1-s)$$

In 1996, Li [101i] defined the sequence of numbers $(\lambda_n)$ by

$$\text{(4.3.302)} \quad \lambda_n = \frac{1}{(n-1)!} \frac{d^n}{ds^n}[s^{n-1}\log\xi(s)]\bigg|_{s=1}$$

and proved that a necessary and sufficient condition for the non-trivial zeros $\rho$ of the Riemann zeta function to lie on the critical line $s = \frac{1}{2} + i\tau$ is that $\lambda_n$ is non-negative for every positive integer $n$. Earlier in 1992, Keiper [83a] showed that if the Riemann hypothesis is true, then $\lambda_n > 0$ for all $n \geq 1$.

Li also showed that

$$\text{(4.3.303)} \quad \lambda_n = \sum_\rho \left[1 - \left(1 - \frac{1}{\rho}\right)^n\right]$$

Taking logarithms of (4.3.300) gives us

(4.3.304) $$\log \xi(s) = -\log 2 + \log[s\Gamma(s/2)] - \frac{s}{2}\log \pi + \log[(s-1)\varsigma(s)]$$

and we also have the Maclaurin expansion about $s = 1$

(4.3.305) $$\log \xi(s) = -\log 2 - \sum_{k=1}^{\infty} (-1)^k \frac{\sigma_k}{k}(s-1)^k$$

which implies that

$$\frac{\xi'(s)}{\xi(s)} = -\sum_{k=1}^{\infty} (-1)^k \sigma_k (s-1)^{k-1}$$

The constant term in (4.3.305) arises because

$$\lim_{s \to 1} \xi(s) = \frac{1}{2}\pi^{-1/2}\Gamma(1/2)\lim_{s \to 1}[(s-1)\varsigma(s)] = \frac{1}{2}$$

and we note that the coefficients $\sigma_k$ are defined by

$$\sigma_k = \frac{(-1)^{k+1}}{(k-1)!}\frac{d^k}{ds^k}\log \xi(s)\bigg|_{s=1}$$

Comparing this with (4.3.302) we immediately see that $\lambda_1 = \sigma_1$.

In passing, we may also note that

$$\lim_{s \to 0}\xi(s) = -\frac{1}{2}\varsigma(0)\lim_{s \to 0}[s\Gamma(s/2)] = -\varsigma(0)\lim_{s \to 0}[s/2\Gamma(s/2)] = \frac{1}{2}$$

in accordance with the functional equation for $\xi(s)$.

Eliminating $\log \xi(s)$ from (4.3.304) and (4.3.305) results in

(4.3.306) $$\log[s\Gamma(s/2)] - \frac{s}{2}\log \pi + \log[(s-1)\varsigma(s)] = -\sum_{k=1}^{\infty}(-1)^k \frac{\sigma_k}{k}(s-1)^k$$

Using the following identity (see [1, p.256] which was also employed by Keiper [83a] and Coffey [45f])

$$\text{(4.3.307)} \quad \log[s\Gamma(s/2)] = \log 2 + \frac{1}{2}(\gamma-1) - \frac{1}{2}(\gamma-1)(s-1) + \sum_{k=2}^{\infty} \frac{[\varsigma(k)-1]}{k2^k}[1-(s-1)]^k$$

we get

$$\log[(s-1)\varsigma(s)] = \frac{s}{2}\log\pi - \sum_{k=1}^{\infty}(-1)^k \frac{\sigma_k}{k}(s-1)^k$$

(4.3.308)

$$-\log 2 - \frac{1}{2}(\gamma-1) + \frac{1}{2}(\gamma-1)(s-1) - \sum_{k=2}^{\infty}\frac{[\varsigma(k)-1]}{k2^k}[1-(s-1)]^k$$

From (4.3.112a) we have for the left-hand side of (4.3.308)

$$\lim_{s\to 1}\log[(s-1)\varsigma(s)] = \log\lim_{s\to 1}[(s-1)\varsigma(s)] = \log 1 = 0$$

and hence we get

$$\text{(4.3.309)} \quad \sum_{k=2}^{\infty}\frac{[\varsigma(k)-1]}{k2^k} = \frac{1}{2}(1-\gamma) + \log\left(\frac{\sqrt{\pi}}{2}\right)$$

in agreement with the formula reported in the book by Srivastava and Choi [126, p.174] and also Keiper's paper [83a]. In fact, Srivastava and Choi [126, p.173] also give the general result valid for $|s|<2$

$$\text{(4.3.310)} \quad \sum_{k=2}^{\infty}\frac{[\varsigma(k)-1]}{k}s^k = \log\Gamma(2-s) + (1-\gamma)s$$

Differentiating this gives us

$$\text{(4.3.311)} \quad \sum_{k=2}^{\infty}[\varsigma(k)-1]s^{k-1} = -\psi(2-s) - \gamma$$

and with $s = 1/2$ we see that

$$\text{(4.3.312)} \quad \sum_{k=2}^{\infty}\frac{[\varsigma(k)-1]}{2^k} = -\frac{1}{2}[\psi(3/2) - \gamma]$$

We have [126, p.20]

$$\text{(4.3.313)} \quad \psi\left(n+\frac{1}{2}\right) = -\gamma - 2\log 2 + 2\sum_{k=0}^{n-1}\frac{1}{2k+1}$$

and hence

(4.3.314) $$\sum_{k=2}^{\infty} \frac{[\varsigma(k)-1]}{2^k} = \log 2 - \frac{1}{2}$$

Similarly, letting $t=1$ in (4.3.310) we obtain the expansion

(4.3.315) $$\sum_{k=2}^{\infty} \frac{[\varsigma(k)-1]}{k} = 1 - \gamma$$

This may also be obtained by letting $s=0$ in (4.3.307).

In (4.3.106a) we let $g(s) = \log[(s-1)\varsigma(s)]$ and we note that $\lim_{s \to 1} g(s) = 0$. We have

$$g(s) = \log\left[\sum_{n=0}^{\infty} \frac{1}{n+1} \sum_{k=0}^{n} \binom{n}{k} \frac{(-1)^k}{(1+k)^{s-1}}\right] = \log h(s)$$

and we also see that $g(1) = 0$. Therefore we have the limit

$$\lim_{s \to 1}\left[\frac{d^k}{ds^k} \log[(s-1)\varsigma(s)]\right] = \lim_{s \to 1} \frac{d^k}{ds^k} g(s)$$

From (4.3.308) we see that

(4.3.316)

$$\frac{d}{ds}\log[(s-1)\varsigma(s)] = \frac{1}{2}(\log \pi + \gamma - 1) - \sum_{k=1}^{\infty}(-1)^k \sigma_k (s-1)^{k-1} + \sum_{k=2}^{\infty} \frac{[\varsigma(k)-1]}{2^k}[1-(s-1)]^{k-1}$$

and hence

$$g'(1) = \frac{d}{ds}\log[(s-1)\varsigma(s)]\bigg|_{s=1} = \frac{1}{2}(\log \pi + \gamma - 1) + \sigma_1 + \sum_{k=2}^{\infty} \frac{[\varsigma(k)-1]}{2^k}$$

$$= \frac{1}{2}(\log \pi + \gamma - 1) + \sigma_1 + \log 2 - \frac{1}{2}$$

With reference to $g(s) = \log h(s)$ and the fact that $h(1) = 1$, elementary calculus results in

$$h^{(p)}(s) = (-1)^p \sum_{n=0}^{\infty} \frac{1}{n+1} \sum_{k=0}^{n} \binom{n}{k} \frac{(-1)^k \log^p(1+k)}{(1+k)^{s-1}}$$

$$h^{(p)}(1) = (-1)^p \sum_{n=0}^{\infty} \frac{1}{n+1} \sum_{k=0}^{n} \binom{n}{k} (-1)^k \log^p(1+k) = (-1)^{p+1} p \gamma_{p-1} \quad \text{for } p \geq 1$$

$$g'(s) = \frac{h'(s)}{h(s)} \qquad g'(1) = h'(1) = \gamma$$

Therefore we have obtained

(4.3.317) $$\sigma_1 = -\frac{1}{2}\log \pi + \frac{1}{2}\gamma + 1 - \log 2$$

as reported, inter alia, by Coffey in [45f]. We have

(4.3.317a) $$\sigma_1 = \lambda_1 \simeq 0.023...$$

Differentiating (4.3.316) results in for $n \geq 1$

$$\frac{d^{n+1}}{ds^{n+1}} \log[(s-1)\varsigma(s)] = -\sum_{k=1}^{\infty} (-1)^k (k-1)(k-2)...(k-n)\sigma_k (s-1)^{k-n-1}$$

$$+ (-1)^n \sum_{k=2}^{\infty} \frac{[\varsigma(k)-1]}{2^k} (k-1)(k-2)...(k-n)[1-(s-1)]^{k-n-1}$$

Therefore we have

(4.3.318)
$$\lim_{s \to 1} \frac{d^{n+1}}{ds^{n+1}} \log[(s-1)\varsigma(s)] = (-1)^n n! \sigma_{n+1} + (-1)^n \sum_{k=2}^{\infty} \frac{[\varsigma(k)-1]}{2^k}(k-1)(k-2)...(k-n)$$

By differentiating (4.3.311) $n$ times we obtain

(4.3.319) $$\sum_{k=2}^{\infty} [\varsigma(k)-1](k-1)(k-2)...(k-n)s^{k-n-1} = (-1)^{n+1} \psi^{(n)}(2-s)$$

and, multiplying across by $s^{n+1}$ and letting $s = 1/2$, we get

(4.3.320) $$\sum_{k=2}^{\infty} \frac{[\varsigma(k)-1]}{2^k}(k-1)(k-2)...(k-n) = \frac{(-1)^{n+1}}{2^{n+1}} \psi^{(n)}(3/2)$$

Since $\psi^{(n)}(1+x) = \psi^{(n)}(x) + (-1)^n n! x^{-n-1}$ we have

(4.3.321) $$\psi^{(n)}(3/2) = \psi^{(n)}(1/2) + (-1)^n n! 2^{n+1}$$

We have

(4.3.322) $$\psi^{(n)}(1/2) = (-1)^{n+1} n! [2^{n+1} - 1]\varsigma(n+1)$$

and hence we get

(4.3.323) $$\psi^{(n)}(3/2) = (-1)^{n+1} n! [2^{n+1} - 1]\varsigma(n+1) + (-1)^n n! 2^{n+1}$$

Accordingly we see that

(4.3.324) $$\sum_{k=2}^{\infty} \frac{[\varsigma(k)-1]}{2^k}(k-1)(k-2)\ldots(k-n) = n!\left(1 - \frac{1}{2^{n+1}}\right)\varsigma(n+1) - n!$$

Hence we see that for $n \geq 1$

(4.3.324a)
$$\lim_{s \to 1} \frac{d^{n+1}}{ds^{n+1}} \log[(s-1)\varsigma(s)] = (-1)^n n! \sigma_{n+1} + (-1)^n n!\left(1 - \frac{1}{2^{n+1}}\right)\varsigma(n+1) - (-1)^n n!$$

and with $n = 1$ we get

$$g''(1) = \lim_{s \to 1} \frac{d^2}{ds^2} \log[(s-1)\varsigma(s)] = -\sigma_2 - \frac{3}{4}\varsigma(2) + 1$$

We have

$$g''(s) = \frac{h(s)h''(s) - [h'(s)]^2}{[h(s)]^2} \qquad g''(1) = h''(1) - [h'(1)]^2$$

Therefore we obtain

$$g''(1) = -2\gamma_1 - \gamma^2$$

$$= -\sigma_2 - \frac{3}{4}\varsigma(2) + 1$$

and hence we see that

(4.3.325) $$\sigma_2 = -\frac{3}{4}\varsigma(2) + 1 + 2\gamma_1 + \gamma^2$$

We have from Keiper's paper [83a]

(4.3.326) $$\sigma_n = \sum_{\rho} \frac{1}{\rho^n}$$

and since

$$\lambda_n = \sum_{\rho}\left[1-\left(1-\frac{1}{\rho}\right)^n\right] = -\sum_{\rho}\sum_{k=0}^{n}\binom{n}{k}\frac{(-1)^k}{\rho^k}$$

we see that

(4.3.327) $$\lambda_n = -\sum_{k=1}^{n}\binom{n}{k}(-1)^k \sigma_k$$

which implies that $\lambda_1 = \sigma_1$.

This gives us in accordance with [45f]

(4.3.328) $$\lambda_2 = 2\sigma_1 - \sigma_2 = -\log\pi + \gamma + 1 - 2\log 2 + \frac{3}{4}\varsigma(2) - 2\gamma_1 - \gamma^2$$

We see from (4.3.324a) that

(4.3.329) $$(-1)^n \sigma_{n+1} = \lim_{s\to 1}\frac{1}{n!}\frac{d^{n+1}}{ds^{n+1}}\log[(s-1)\varsigma(s)] - (-1)^n\left(1-\frac{1}{2^{n+1}}\right)\varsigma(n+1) + (-1)^n$$

and we have

$$\lambda_m = -\sum_{n=1}^{m}\binom{m}{n}(-1)^n \sigma_n$$

$$= \sum_{n=0}^{m}\binom{m}{n+1}(-1)^n \sigma_{n+1} = \sum_{n=1}^{m}\binom{m}{n+1}(-1)^n \sigma_{n+1} + m\sigma_1$$

since $\binom{m}{m+1} = 0$. Therefore we obtain

$$\lambda_m = \lim_{s\to 1}\sum_{n=1}^{m}\binom{m}{n+1}\frac{1}{n!}\frac{d^{n+1}}{ds^{n+1}}\log[(s-1)\varsigma(s)]$$

$$-\sum_{n=1}^{m}\binom{m}{n+1}(-1)^{n}\left(1-\frac{1}{2^{n+1}}\right)\varsigma(n+1)+\sum_{n=1}^{m}\binom{m}{n+1}(-1)^{n}+m\sigma_{1}$$

This may be written as

$$\lambda_{m}=\lim_{s\to 1}\sum_{n=1}^{m}\binom{m}{n+1}\frac{1}{n!}\frac{d^{n+1}}{ds^{n+1}}\log[(s-1)\varsigma(s)]+\sum_{r=2}^{m}\binom{m}{r}(-1)^{r}\left(1-\frac{1}{2^{r}}\right)\varsigma(r)+1-m+m\sigma_{1}$$

since $\sum_{n=1}^{m}\binom{m}{n+1}(-1)^{n+1}=1-m$. Employing Coffey's notation [45g] we have

(4.3.330)

$$\lambda_{m}=1-\frac{m}{2}[\log\pi-\gamma+2\log 2]+S_{1}(m)+\lim_{s\to 1}\sum_{n=1}^{m}\binom{m}{n+1}\frac{1}{n!}\frac{d^{n+1}}{ds^{n+1}}\log[(s-1)\varsigma(s)]$$

where

(4.3.331) $$S_{1}(m)=\sum_{r=2}^{m}\binom{m}{r}(-1)^{r}\left(1-\frac{1}{2^{r}}\right)\varsigma(r)$$

Using the Leibniz rule we see that

$$\lim_{s\to 1}\frac{1}{(m-1)!}\frac{d^{m}}{ds^{m}}\left[s^{m-1}\log[(s-1)\varsigma(s)]\right]=\lim_{s\to 1}\sum_{n=1}^{m}\binom{m}{n}\frac{1}{(n-1)!}\frac{d^{n}}{ds^{n}}\log[(s-1)\varsigma(s)]$$

$$=\lim_{s\to 1}\sum_{n=0}^{m}\binom{m}{n+1}\frac{1}{n!}\frac{d^{n+1}}{ds^{n+1}}\log[(s-1)\varsigma(s)]$$

$$=\lim_{s\to 1}\sum_{n=1}^{m}\binom{m}{n+1}\frac{1}{n!}\frac{d^{n+1}}{ds^{n+1}}\log[(s-1)\varsigma(s)]+\lim_{s\to 1}\binom{m}{1}\frac{d}{ds}\log[(s-1)\varsigma(s)]$$

Therefore we obtain

$$\lim_{s\to 1}\frac{1}{(m-1)!}\frac{d^{m}}{ds^{m}}\left[s^{m-1}\log[(s-1)\varsigma(s)]\right]=\lim_{s\to 1}\sum_{n=1}^{m}\binom{m}{n+1}\frac{1}{n!}\frac{d^{n+1}}{ds^{n+1}}\log[(s-1)\varsigma(s)]+m\gamma$$

and hence we derive Coffey's expression

$$\lambda_m = 1 - \frac{m}{2}[\log \pi - \gamma + 2\log 2] + S_1(m) + \lim_{s \to 1} \frac{1}{(m-1)!} \frac{d^m}{ds^m}\left[s^{m-1}\log[(s-1)\varsigma(s)]\right] - m\gamma$$

or equivalently

(4.3.332) $$\lambda_m = 1 - \frac{m}{2}[\log \pi + \gamma + 2\log 2] + S_1(m) + S_2(m)$$

where

(4.3.333) $$S_2(m) = \lim_{s \to 1} \sum_{n=1}^{m} \binom{m}{n} \frac{1}{(n-1)!} \frac{d^n}{ds^n} \log[(s-1)\varsigma(s)]$$

We have

(4.3.334) $$\varsigma(r, 3/2) = (2^r - 1)\varsigma(r) - 2^r$$

and therefore we see that

$$1 + \frac{\varsigma(r, 3/2)}{2^r} = \left(1 - \frac{1}{2^r}\right)\varsigma(r)$$

Hence we may write

$$S_1(m) = \sum_{r=2}^{m} \binom{m}{r}(-1)^r \left(1 + \frac{\varsigma(r, 3/2)}{2^r}\right)$$

$$= \sum_{r=0}^{m} \binom{m}{r}(-1)^r - \binom{m}{0} + \binom{m}{1} + \sum_{r=2}^{m} \binom{m}{r}(-1)^r \frac{\varsigma(r, 3/2)}{2^r}$$

$$= m - 1 + \sum_{r=2}^{m} \binom{m}{r}(-1)^r \frac{\varsigma(r, 3/2)}{2^r}$$

Coffey [45f] has shown that

$$\frac{m}{2}\log m + (\gamma - 1)\frac{m}{2} + \frac{1}{2} \leq S_1(m) \leq \frac{m}{2}\log m + (\gamma + 1)\frac{m}{2} - \frac{1}{2}$$

and hence we have

$$1 - \frac{m}{2}[\log \pi + \gamma + 2\log 2] + S_1(m) \geq \frac{m}{2}\log m + (\gamma - 1)\frac{m}{2} + \frac{1}{2} + 1 - \frac{m}{2}[\log \pi + \gamma + 2\log 2]$$

$$= \frac{m}{2}\left[\log\frac{m}{4\pi}-1\right]+\frac{3}{2}$$

Referring to Kölbig's paper [91aa] (and see also Coffey [45h]) we may write

(4.3.335) $$\frac{d^n}{ds^n}\left[s^{n-1}\log[(s-1)\varsigma(s)]\right]=\frac{d^n}{ds^n}\exp\log\left[s^{n-1}\log[(s-1)\varsigma(s)]\right]$$

and then express this in terms of the complete Bell polynomials (albeit this operation does not really simplify the output).

The constants $\eta_n$ are defined by reference to the logarithmic derivative of the Riemann zeta function

(4.3.336) $$\frac{d}{ds}[\log\varsigma(s)]=\frac{\varsigma'(s)}{\varsigma(s)}=-\frac{1}{s-1}-\sum_{k=0}^{\infty}\eta_k(s-1)^k \qquad |s-1|<3$$

and we may also note that

$$\frac{d}{ds}\log[(s-1)\varsigma(s)]=\frac{\varsigma'(s)}{\varsigma(s)}+\frac{1}{s-1}=-\sum_{k=0}^{\infty}\eta_k(s-1)^k$$

We then see that

$$\frac{d^{n+1}}{ds^{n+1}}\log[(s-1)\varsigma(s)]=-\sum_{k=0}^{\infty}\eta_k k(k-1)\ldots(k-n+1)(s-1)^{k-n}$$

and hence we get

(4.3.337) $$\lim_{s\to 1}\frac{d^{n+1}}{ds^{n+1}}\log[(s-1)\varsigma(s)]=-n!\eta_n$$

Referring to (4.3.329)

$$(-1)^n\sigma_{n+1}=\lim_{s\to 1}\frac{1}{n!}\frac{d^{n+1}}{ds^{n+1}}\log[(s-1)\varsigma(s)]-(-1)^n\left(1-\frac{1}{2^{n+1}}\right)\varsigma(n+1)+(-1)^n$$

we then obtain

$$(-1)^n\sigma_{n+1}=-\eta_n-(-1)^n\left(1-\frac{1}{2^{n+1}}\right)\varsigma(n+1)+(-1)^n$$

(4.3.338) $$\sigma_{n+1} = (-1)^{n+1}\eta_n - \left(1 - \frac{1}{2^{n+1}}\right)\varsigma(n+1) + 1$$

in accordance with Theorem 3 in Coffey's paper [45f].

In 2004 Maślanka [101c] reports that the following expression was also implicitly given in the 1999 paper by Bombieri and Lagarias [24b]

$$\lambda_m = 1 - \frac{m}{2}[\log \pi + \gamma + 2\log 2] + S_1(m) + S_2(m)$$

and Maślanka makes the decomposition

$$\lambda_m = \bar{\lambda}_m + \tilde{\lambda}_m$$

where the "trend" $\bar{\lambda}_m$ is strictly increasing and the "oscillations" $\tilde{\lambda}_m$ are given by

$$\bar{\lambda}_m = 1 - \frac{m}{2}[\log \pi + \gamma + 2\log 2] + S_1(m)$$

$$\tilde{\lambda}_m = -\sum_{n=1}^{m} \binom{m}{n} \eta_{n-1}$$

From (4.3.337) we see that

$$\lim_{s \to 1} \frac{1}{(n-1)!} \frac{d^n}{ds^n} \log[(s-1)\varsigma(s)] = -\eta_{n-1}$$

and therefore we have

$$-\sum_{n=1}^{m} \binom{m}{n} \eta_{n-1} = \lim_{s \to 1} \sum_{n=1}^{m} \binom{m}{n} \frac{1}{(n-1)!} \frac{d^n}{ds^n} \log[(s-1)\varsigma(s)]$$

This then shows that

(4.3.339) $$\tilde{\lambda}_m = S_2(m) = -\sum_{n=1}^{m} \binom{m}{n} \eta_{n-1}$$

and we have

$$S_2(m) = \lim_{s \to 1} \sum_{n=1}^{m} \binom{m}{n} \frac{1}{(n-1)!} \frac{d^n}{ds^n} \log[(s-1)\varsigma(s)]$$

$$\frac{d^n}{ds^n}\log[(s-1)\varsigma(s)] = \frac{d^n}{ds^n}\log \sum_{m=0}^{\infty}\frac{1}{m+1}\sum_{p=0}^{m}\binom{m}{p}\frac{(-1)^p}{(1+p)^{s-1}}$$

In 1968 Gould [123c] reminded the mathematical community of the "not well-known" formula for the $n$ th derivative of a composite function $f(z)$ where $z$ is a function of $x$, namely

(4.3.340) $$D_x^{(n)}f(z) = \sum_{k=1}^{n} D_z^{(k)}f(z)\frac{(-1)^k}{k!}\sum_{j=1}^{k}(-1)^j\binom{k}{j}z^{k-j} \quad , \text{ for } n \geq 1$$

This differentiation algorithm was also reported by Gould [73a] in 1972 in a somewhat different form in the case where $z$ is a function of $s$

(4.3.341) $$D_s^{(n)}f(z) = \sum_{k=0}^{n} D_z^{(k)}f(z)\frac{(-1)^k}{k!}\sum_{j=0}^{k}(-1)^j\binom{k}{j}z^{k-j}D_s^{(n)}z^j \quad , \text{ for } n \geq 1$$

This expression is frequently easier to handle than the di Bruno algorithm and it may be worthwhile using it in an attempt to simplify, and shed further light on, the structure of $\tilde{\lambda}_m = S_2(m)$.

Letting $f(z) = \log[z(s)]$ and $z(s) = \sum_{q=0}^{\infty}\frac{1}{q+1}\sum_{r=0}^{q}\binom{q}{r}\frac{(-1)^r}{(1+r)^{s-1}}$.

we see that

$$z(1) = 1 \text{ and } f(1) = \log[z(1)] = 0 \quad D_z^{(k)}f(z) = D_z^{(k)}\log[z(s)] = (-1)^{k-1}(k-1)!z^{-k}$$

and we have

$$D_s^{(n)}f(z)\Big|_{s=1} = \sum_{k=1}^{n}(-1)^{k-1}(k-1)!z^{-k}\frac{(-1)^k}{k!}\sum_{j=0}^{k}(-1)^j\binom{k}{j}z^{k-j}D_s^{(n)}z^j\Big|_{s=1}$$

$$= -\sum_{k=1}^{n}\frac{1}{k}\sum_{j=0}^{k}(-1)^j\binom{k}{j}z^{-j}D_s^{(n)}z^j\Big|_{s=1}$$

Letting $s = \frac{1}{1-z}$ we have $(s-1)\varsigma(s) \to \frac{z}{1-z}\varsigma\left(\frac{1}{1-z}\right)$ and note that the function

$F(z) = \log\left[\dfrac{z}{1-z}\varsigma\left(\dfrac{1}{1-z}\right)\right]$ has been considered by Smith [120ai] and Coffey [45g] in connection with the Riemann Hypothesis. In [45g] Coffey reports that

(4.3.341) $\qquad F(z) = -\sum_{n=1}^{\infty}\dfrac{(-1)^n}{n}\left[\sum_{j=1}^{\infty}\dfrac{(-1)^j}{j!}\gamma_j z^{j+1}\left(\sum_{p=0}^{\infty}z^p\right)^{j+1}\right]^n$

and an alternative proof is shown below using the formula for the Stieltjes constants $\gamma_j$ (it should however be noted that in the expression used by Coffey the summation starts at $j=0$)

$$= -\sum_{n=1}^{\infty}\dfrac{(-1)^n}{n}\left[\sum_{j=1}^{\infty}\dfrac{(-1)^j}{j!}\left\{-\dfrac{1}{j+1}\sum_{m=0}^{\infty}\dfrac{1}{m+1}\sum_{k=0}^{m}\binom{m}{k}(-1)^k\log^{j+1}(1+k)\right\}z^{j+1}\left(\sum_{p=0}^{\infty}z^p\right)^{j+1}\right]^n$$

$$= -\sum_{n=1}^{\infty}\dfrac{(-1)^n}{n}\left[\sum_{m=0}^{\infty}\dfrac{1}{m+1}\sum_{k=0}^{m}\binom{m}{k}(-1)^k\sum_{j=1}^{\infty}\dfrac{(-1)^{j+1}\log^{j+1}(1+k)\left(\sum_{p=0}^{\infty}z^p\right)^{j+1}z^{j+1}}{(j+1)!}\right]^n$$

Since $\sum_{p=0}^{\infty}z^p = \dfrac{1}{1-z}$ for $|z|<1$ this becomes

$$= -\sum_{n=1}^{\infty}\dfrac{(-1)^n}{n}\left[\sum_{m=0}^{\infty}\dfrac{1}{m+1}\sum_{k=0}^{m}\binom{m}{k}(-1)^k\sum_{j=1}^{\infty}\dfrac{(-1)^{j+1}\log^{j+1}(1+k)\left(\dfrac{z}{1-z}\right)^{j+1}}{(j+1)!}\right]^n$$

$$= -\sum_{n=1}^{\infty}\dfrac{(-1)^n}{n}\left[\sum_{m=0}^{\infty}\dfrac{1}{m+1}\sum_{k=0}^{m}\binom{m}{k}(-1)^k\exp\left\{-\left(\dfrac{z}{1-z}\right)\log(1+k)\right\}\right]^n$$

$$= -\sum_{n=1}^{\infty}\dfrac{(-1)^n}{n}\left[\sum_{m=0}^{\infty}\dfrac{1}{m+1}\sum_{k=0}^{m}\binom{m}{k}\dfrac{(-1)^k}{(1+k)^{\frac{z}{1-z}}}\right]^n$$

and using the Hasse formula for the Riemann zeta function this becomes

$$= -\sum_{n=1}^{\infty} \frac{(-1)^n}{n} \left[ \frac{z}{1-z} \varsigma\left(\frac{1}{1-z}\right) \right]^n$$

$$= \log\left[ \frac{z}{1-z} \varsigma\left(\frac{1}{1-z}\right) \right]$$

□

With reference to Maślanka's paper [101c] we consider the sum

$$\sum_{n=0}^{\infty} \gamma_n s^n = -\sum_{n=0}^{\infty} s^n \frac{1}{n+1} \sum_{m=0}^{\infty} \frac{1}{m+1} \sum_{k=0}^{m} \binom{m}{k} (-1)^k \log^{n+1}(1+k)$$

$$= -\frac{1}{s} \sum_{m=0}^{\infty} \frac{1}{m+1} \sum_{k=0}^{m} \binom{m}{k} (-1)^k \sum_{n=0}^{\infty} \frac{s^{n+1} \log^{n+1}(1+k)}{n+1}$$

$$= \frac{1}{s} \sum_{m=0}^{\infty} \frac{1}{m+1} \sum_{k=0}^{m} \binom{m}{k} (-1)^k \log[1 - s\log(1+k)]$$

□

Another approach to the Riemann Hypothesis is outlined below. We slightly rearrange the definition of the Riemann xi function as below

$$\xi(s) = \frac{1}{2} s(s-1) \pi^{-s/2} \Gamma(s/2) \varsigma(s) = \pi^{-s/2} \frac{s}{2} \Gamma\left(\frac{s}{2}\right)(s-1)\varsigma(s) = \pi^{-s/2} \Gamma\left(1 + \frac{s}{2}\right)(s-1)\varsigma(s)$$

and upon taking logarithms we have

$$\log \xi(s) = -\frac{s}{2} \log \pi + \log \Gamma\left(1 + \frac{s}{2}\right) + \log[(s-1)\varsigma(s)]$$

(the presence of the $(s-1)\varsigma(s)$ term indicates a source for the Stieltjes constants)

Using the Hadamard representation of the Riemann zeta function

(4.3.343) $$\varsigma(s) = \frac{e^{cs}}{2(s-1)\Gamma\left(1+\frac{s}{2}\right)} \prod_\rho \left(1 - \frac{s}{\rho}\right) e^{s/\rho}$$

where $c = \log(2\pi) - 1 - \gamma/2$, Snowden [120aa] showed in 2001 that

$$(4.3.344) \quad \log[(s-1)\varsigma(s)] = -\log 2 + s[\log(2\pi)-1] + \sum_{k=2}^{\infty}\left[\frac{(-1)^{k+1}}{2^k}\varsigma(k)-\sigma_k\right]\frac{s^k}{k}$$

We also have the Maclaurin expansion for the gamma function (E.22n) in Volume VI

$$\log\Gamma(1+x) = -\gamma x + \sum_{k=2}^{\infty}\frac{(-1)^k}{k}\varsigma(k)x^k$$

Therefore we have for $|s| \leq 2$

$$(4.3.345) \quad \log\Gamma\left(1+\frac{s}{2}\right) = -\frac{1}{2}\gamma s + \sum_{k=2}^{\infty}\frac{(-1)^k}{k2^k}\varsigma(k)s^k$$

and this gives us the expansion

$$\log\xi(s) = -\log 2 + s[\log(2\pi)-1] - \frac{s}{2}\log\pi - \frac{1}{2}\gamma s + \sum_{k=2}^{\infty}\frac{(-1)^k}{k2^k}\varsigma(k)s^k - \sum_{k=2}^{\infty}\left[\frac{(-1)^k}{2^k}\varsigma(k)+\sigma_k\right]\frac{s^k}{k}$$

which becomes

$$(4.3.346) \quad \log\xi(s) = \frac{1}{2}[\log\pi-\gamma]s - \log 2 - s[1-\log 2] - \sum_{k=2}^{\infty}\frac{\sigma_k}{k}s^k$$

In view of (4.3.317) this may be written as

$$(4.3.347) \quad \log\xi(s) = -\sigma_1 s - \log 2 - \sum_{k=2}^{\infty}\frac{\sigma_k}{k}s^k = -\log 2 - \sum_{k=1}^{\infty}\frac{\sigma_k}{k}s^k$$

and, having regard to the functional equation (4.3.301) $\xi(s)=\xi(1-s)$, we immediately see that (4.3.347) is simply another representation of (4.3.305). Indeed we have

$$(4.3.348) \quad \sum_{k=1}^{\infty}\frac{\sigma_k}{k}s^k = \sum_{k=1}^{\infty}(-1)^k\frac{\sigma_k}{k}(s-1)^k$$

from which we see that

$$(4.3.349) \quad \sum_{k=1}^{\infty}\frac{\sigma_k}{k} = 0$$

Letting $s \to s-1$ we see that

$$\sum_{k=1}^{\infty}\frac{\sigma_k}{k}(s-1)^k = \sum_{k=1}^{\infty}(-1)^k\frac{\sigma_k}{k}(s-2)^k$$

and with $s=2$ (valid?) we obtain

$$\sum_{k=1}^{\infty}(-1)^k\frac{\sigma_k}{k} = 0$$

and hence we have

$$\sum_{k=1}^{\infty}\frac{\sigma_{2k}}{k} = \sum_{k=1}^{\infty}\frac{\sigma_{2k+1}}{2k+1} = 0$$

Differentiating (4.3.348) results in (we need to prove that term by term differentiation is valid)

(4.3.350) $$\sum_{k=1}^{\infty}\sigma_k s^{k-1} = \sum_{k=1}^{\infty}(-1)^k \sigma_k (s-1)^{k-1}$$

and with $s=1$ we obtain

(4.3.351) $$\sum_{k=1}^{\infty}\sigma_k = -\sigma_1$$

as reported by Keiper [83a].

Let us now multiply (4.3.350) by $s$

$$\sum_{k=1}^{\infty}\sigma_k s^k = \sum_{k=1}^{\infty}(-1)^k \sigma_k s(s-1)^{k-1}$$

and differentiate to obtain

(4.3.352) $$\sum_{k=1}^{\infty}k\sigma_k s^{k-1} = \sum_{k=1}^{\infty}(-1)^k \sigma_k \left[s(k-1)(s-1)^{k-2} + (s-1)^{k-1}\right]$$

Letting $s=1$ we have

$$\sum_{k=1}^{\infty}k\sigma_k = \sigma_2 - \sigma_1$$

We now repeat the same operation with (4.3.352)

$$\sum_{k=1}^{\infty} k\sigma_k s^k = \sum_{k=1}^{\infty} (-1)^k \sigma_k \left[ s^2(k-1)(s-1)^{k-2} + s(s-1)^{k-1} \right]$$

$$\sum_{k=1}^{\infty} k^2 \sigma_k s^{k-1} = \sum_{k=1}^{\infty} (-1)^k \sigma_k \left[ s^2(k-1)(k-2)(s-1)^{k-3} + 3s(k-1)(s-1)^{k-2} + (s-1)^{k-1} \right]$$

Letting $s = 1$ we have

$$\sum_{k=1}^{\infty} k^2 \sigma_k = -\sigma_3 + 3\sigma_2 - \sigma_1$$

Having regard to the definition of the Li/Keiper constants

$$\lambda_n = \frac{1}{(n-1)!} \frac{d^n}{ds^n} [s^{n-1} \log \xi(s)] \Big|_{s=1}$$

we consider

$$s^{n-1} \log \xi(s) = -s^{n-1} \log 2 - \sum_{k=1}^{\infty} \frac{\sigma_k}{k} s^{k+n-1} = -s^{n-1} \log 2 - \sum_{k=1}^{\infty} (-1)^k \frac{\sigma_k}{k} s^{n-1} (s-1)^k$$

whereupon we obtain

$$\frac{d^n}{ds^n} [s^{n-1} \log \xi(s)] = -\sum_{k=1}^{\infty} \frac{(k+n-1)(k+n-2)\ldots k \sigma_k}{k} s^{k-1}$$

$$\frac{d^n}{ds^n} [s^{n-1} \log \xi(s)] = -\frac{d^n}{ds^n} \sum_{k=1}^{\infty} (-1)^k \frac{\sigma_k}{k} s^{n-1} (s-1)^k$$

For $n = 1$ we get

$$\frac{d}{ds} [\log \xi(s)] = -\frac{d}{ds} \sum_{k=1}^{\infty} (-1)^k \frac{\sigma_k}{k} (s-1)^k$$

$$= -\sum_{k=1}^{\infty} (-1)^k \frac{k\sigma_k}{k} (s-1)^{k-1}$$

and therefore we see that $\lambda_1 = \sigma_1$.

For $n = 2$ we have

$$\frac{d^2}{ds^2} [s \log \xi(s)] = -\frac{d^2}{ds^2} \sum_{k=1}^{\infty} (-1)^k \frac{\sigma_k}{k} s(s-1)^k$$

We have

$$\frac{d}{ds}\sum_{k=1}^{\infty}(-1)^k \frac{\sigma_k}{k} s(s-1)^k = \sum_{k=1}^{\infty}(-1)^k \frac{\sigma_k}{k} sk(s-1)^{k-1} + \sum_{k=1}^{\infty}(-1)^k \frac{\sigma_k}{k}(s-1)^k$$

$$\frac{d^2}{ds^2}\sum_{k=1}^{\infty}(-1)^k \frac{\sigma_k}{k} s(s-1)^k = \sum_{k=1}^{\infty}(-1)^k \frac{\sigma_k}{k} sk(k-1)(s-1)^{k-2} + 2\sum_{k=1}^{\infty}(-1)^k \frac{\sigma_k}{k} k(s-1)^{k-1}$$

which evaluated at $s=1$ becomes $\sigma_2 - 2\sigma_1$. Hence we obtain

$$\lambda_2 = -\sigma_2 + 2\sigma_1$$

More generally, using the Leibniz differentiation formula we obtain (4.3.327)

$$\lambda_n = -\sum_{k=1}^{n}\binom{n}{k}(-1)^k \sigma_k$$

## A MULTITUDE OF IDENTITIES AND THEOREMS

This part contains some ancillary theorems and miscellaneous comments.

**(i)** An alternative proof of (4.2.1) is set out below.

**Theorem 4.1:**

(4.4.1) $$\frac{n!}{x(1+x)...(n+x)} = \int_0^1 t^{x-1}(1-t)^n dt = \sum_{k=0}^{n}\binom{n}{k}\frac{(-1)^k}{k+x}$$

where $x > 0$.

**Proof:**

Consider the function $g(x)$ defined by

(4.4.2) $$g(x) = \frac{n!}{x(1+x)...(n+x)}$$

Using (4.3.3) we have

$$（4.4.3） \quad = \frac{n!\Gamma(x)}{\Gamma(x+n+1)}$$

$$（4.4.4） \quad = \frac{\Gamma(n+1)\Gamma(x)}{\Gamma(x+n+1)}$$

$$（4.4.5） \quad = B(n+1, x)$$

$$（4.4.6） \quad = B(x, n+1)$$

The function $g(x)$ is therefore related to the beta function defined by the integral

$$（4.4.7） \quad B(x, y) = \int_0^1 t^{x-1}(1-t)^{y-1} dt \quad , (\operatorname{Re}(x) > 0, \operatorname{Re}(y) > 0)$$

which, in turn, is related to the gamma function via the identity [115, p.193]

$$（4.4.8） \quad B(x, y) = \frac{\Gamma(x)\Gamma(y)}{\Gamma(x+y)}$$

However, because this paper is deliberately pitched at an elementary level (well, it started out that way!), we shall not assume any prior knowledge of the fundamental identity (4.4.8). Instead, let us consider a function $h(x)$ defined by

$$（4.4.9） \quad h(x) = \int_0^1 t^{x-1}(1-t)^n dt$$

We can use integration by parts to evaluate this integral

$$= \frac{(1-t)^n t^x}{x}\bigg|_0^1 + \frac{n}{x}\int_0^1 (1-t)^n t^x dt$$

$$= \frac{n}{x}\int_0^1 (1-t)^{n-1} t^x dt$$

$$= \frac{n(n-1)}{x(x+1)}\int_0^1 (1-t)^{n-2} t^{x+1} dt$$

$$= \frac{n(n-1)\ldots 1}{x(x+1)\ldots(x+n-1)}\int_0^1 t^{x+n-1} dt$$

and therefore we obtain

(4.4.10) $$h(x) = \frac{n!}{x(x+1)...(x+n)}$$

From (4.4.2) and (4.4.10) we therefore see that $g(x) = h(x)$.

Let us now evaluate the integral (4.4.9) in a different way using the binomial theorem (4.1.2). We have

(4.4.11) $$\int_0^1 t^{x-1}(1-t)^n dt = \int_0^1 t^{x-1} \sum_{k=0}^n \binom{n}{k}(-1)^k t^k dt$$

(4.4.11a) $$= \sum_{k=0}^n \binom{n}{k} \frac{(-1)^k}{k+x}$$

Hence we deduce

(4.4.11b) $$g(x) = \frac{n!}{x(1+x)...(n+x)} = \int_0^1 t^{x-1}(1-t)^n dt = \sum_{k=0}^n \binom{n}{k}\frac{(-1)^k}{k+x}$$

(4.4.11c) $$= \frac{\Gamma(n+1)\Gamma(x)}{\Gamma(x+n+1)} = B(n+1,x) = B(x,n+1)$$

From (4.4.11) we also note that

$$\int_0^1 t^{x-1}(1-yt)^n dt = \int_0^1 \sum_{k=0}^n \binom{n}{k}(-1)^k y^k t^{k+x-1} dt$$

and hence

(4.4.11d) $$\int_0^1 t^{x-1}(1-yt)^n dt = \sum_{k=0}^n \binom{n}{k}(-1)^k \frac{y^k}{k+x}$$

Differentiation with respect to $x$ results in

(4.4.12) $$\int_0^1 t^{x-1}(1-yt)^n \log^s t\, dt = (-1)^s s! \sum_{k=0}^n \binom{n}{k}(-1)^k \frac{y^k}{(k+x)^{s+1}}$$

Completing the summation we find

(4.4.13a) $$\int_0^1 t^{x-1} Li_p(1-yt) \log^s t \, dt = (-1)^s s! \sum_{n=1}^{\infty} \frac{1}{n^p} \sum_{k=0}^{n} \binom{n}{k} (-1)^k \frac{y^k}{(k+x)^{s+1}}$$

(4.4.13b) $$\int_0^1 \frac{t^{x-1} \log^s t}{1-u(1-yt)} \, dt = (-1)^s s! \sum_{n=0}^{\infty} u^n \sum_{k=0}^{n} \binom{n}{k} (-1)^k \frac{y^k}{(k+x)^{s+1}}$$

Letting $u = 1/2$ in (4.4.13b) results in

(4.4.13c) $$\int_0^1 \frac{t^{x-1} \log^s t}{1+yt} \, dt = (-1)^s s! \sum_{n=0}^{\infty} \frac{1}{2^{n+1}} \sum_{k=0}^{n} \binom{n}{k} (-1)^k \frac{y^k}{(k+x)^{s+1}}$$

Letting $y = 1$ gives us

(4.4.13d) $$\int_0^1 \frac{t^{x-1} \log^s t}{1+t} \, dt = (-1)^s s! \sum_{n=0}^{\infty} \frac{1}{2^{n+1}} \sum_{k=0}^{n} \binom{n}{k} \frac{(-1)^k}{(k+x)^{s+1}}$$

With $x = 1$ this becomes

(4.4.13e) $$\int_0^1 \frac{\log^s t}{1+t} \, dt = (-1)^s s! \sum_{n=0}^{\infty} \frac{1}{2^{n+1}} \sum_{k=0}^{n} \binom{n}{k} \frac{(-1)^k}{(k+1)^{s+1}}$$

We have from (3.86ii) in Volume I

$$\int_0^1 \frac{\log^s t}{1+t} \, dt = (-1)^s s! \varsigma_a(s+1)$$

and we have therefore another derivation of the Hasse/Sondow identity (3.11)

With $y = -1$ we have

(4.4.13f) $$\int_0^1 \frac{t^{x-1} \log^s t}{1-t} \, dt = (-1)^s s! \sum_{n=0}^{\infty} \frac{1}{2^{n+1}} \sum_{k=0}^{n} \binom{n}{k} \frac{1}{(k+x)^{s+1}}$$

We have from (4.4.231)

$$\int_0^1 \frac{\log^s t}{1-t} \, dt = \int_0^1 \frac{\log^s(1-t)}{t} \, dt = (-1)^s \varsigma(s+1) \Gamma(s+1) \qquad , s \geq 1$$

and we therefore obtain

(4.4.14) $$\varsigma(s+1) = \sum_{n=0}^{\infty} \frac{1}{2^{n+1}} \sum_{k=0}^{n} \binom{n}{k} \frac{1}{(k+x)^{s+1}}$$

Further integrals may be obtained by multiplying (4.4.13b) by $\log^p y$ and integrating using (3.237i) ff.

The Wolfram Integrator gives us $\int \frac{t^{x-1} \log^s t}{1+t} dt$ in terms of generalised hypergeometric functions.

Let us make a bold step and assume that (4.4.13b) is valid where $s$ is a real number (and is not simply restricted to positive integers): we then **conjecture** that

$$\int_0^1 \frac{t^{x-1} \log^s t}{1-u(1-yt)} dt = \cos(\pi s)\Gamma(s+1) \sum_{n=0}^{\infty} u^n \sum_{k=0}^{n} \binom{n}{k} (-1)^k \frac{y^k}{(k+x)^{s+1}}$$

and differentiating with respect to $s$ gives us

$$\int_0^1 \frac{t^{x-1} \log^s t \log(\log t)}{1-u(1-yt)} dt = \cos(\pi s)\Gamma(s+1) \sum_{n=0}^{\infty} u^n \sum_{k=0}^{n} \binom{n}{k} (-1)^k \frac{y^k \log(k+x)}{(k+x)^{s+1}}$$

$$+ \cos(\pi s)\Gamma'(s+1) \sum_{n=0}^{\infty} u^n \sum_{k=0}^{n} \binom{n}{k} (-1)^k \frac{y^k}{(k+x)^{s+1}} - \pi \sin(\pi s)\Gamma(s+1) \sum_{n=0}^{\infty} u^n \sum_{k=0}^{n} \binom{n}{k} (-1)^k \frac{y^k}{(k+x)^{s+1}}$$

With $s=0$ we get

$$\int_0^1 \frac{t^{x-1} \log(\log t)}{1-u(1-yt)} dt = -\sum_{n=0}^{\infty} u^n \sum_{k=0}^{n} \binom{n}{k} (-1)^k \frac{y^k \log(k+x)}{k+x} - \Gamma'(1) \sum_{n=0}^{\infty} u^n \sum_{k=0}^{n} \binom{n}{k} (-1)^k \frac{y^k}{k+x}$$

With $s=1$ we get

$$\int_0^1 \frac{t^{x-1} \log t \log(\log t)}{1-u(1-yt)} dt = -\sum_{n=0}^{\infty} u^n \sum_{k=0}^{n} \binom{n}{k} (-1)^k \frac{y^k \log(k+x)}{(k+x)^2} - \Gamma'(2) \sum_{n=0}^{\infty} u^n \sum_{k=0}^{n} \binom{n}{k} (-1)^k \frac{y^k}{(k+x)^2}$$

Integration with respect to $u$ over the interval $[0,1]$ gives us

$$-\int_0^1 \frac{t^{x-1}\log^s t \log(1-yt)\log(\log t)}{1-yt}dt = \cos(\pi s)\Gamma(s+1)\sum_{n=0}^{\infty}\frac{1}{n+1}\sum_{k=0}^{n}\binom{n}{k}(-1)^k \frac{y^k \log(k+x)}{(k+x)^{s+1}}$$

$$+[\cos(\pi s)\Gamma'(s+1) - \pi\sin(\pi s)\Gamma(s+1)]\sum_{n=0}^{\infty}\frac{1}{n+1}\sum_{k=0}^{n}\binom{n}{k}(-1)^k \frac{y^k}{(k+x)^{s+1}}$$

and with $y=1$ this becomes

$$-\int_0^1 \frac{t^{x-1}\log^s t \log(1-t)\log(\log t)}{1-t}dt = \cos(\pi s)\Gamma(s+1)\sum_{n=0}^{\infty}\frac{1}{n+1}\sum_{k=0}^{n}\binom{n}{k}(-1)^k \frac{\log(k+x)}{(k+x)^{s+1}}$$

$$+[\cos(\pi s)\Gamma'(s+1) - \pi\sin(\pi s)\Gamma(s+1)]\sum_{n=0}^{\infty}\frac{1}{n+1}\sum_{k=0}^{n}\binom{n}{k}\frac{(-1)^k}{(k+x)^{s+1}}$$

Referring to (4.3.107a) we see that

$$\sum_{n=0}^{\infty}\frac{1}{n+1}\sum_{k=0}^{n}\binom{n}{k}(-1)^k \frac{\log(k+x)}{(k+x)^{s+1}} = -(s+1)\varsigma'(s+2,x) - \varsigma(s+2,x)$$

and we therefore conjecture that

$$\int_0^1 \frac{t^{x-1}\log^s t \log(1-t)\log(\log t)}{1-t}dt = \cos(\pi s)\Gamma(s+1)[(s+1)\varsigma'(s+2,x) + \varsigma(s+2,x)]$$

$$-[\cos(\pi s)\Gamma'(s+1) - \pi\sin(\pi s)\Gamma(s+1)](s+1)\varsigma(s+2,x)$$

With $x=1$ we get

$$\int_0^1 \frac{\log^s t \log(1-t)\log(\log t)}{1-t}dt = \cos(\pi s)\Gamma(s+1)[(s+1)\varsigma'(s+2) + \varsigma(s+2)]$$

$$-[\cos(\pi s)\Gamma'(s+1) - \pi\sin(\pi s)\Gamma(s+1)](s+1)\varsigma(s+2)$$

and with $s=0$ we might have

$$\int_0^1 \frac{\log(1-t)\log(\log t)}{1-t}dt = \varsigma'(2) + \varsigma(2) + \gamma$$

Alternatively, we may write (4.4.13b) as

$$\int_0^1 \frac{t^{x-1}(-\log t)^s}{1-u(1-yt)}\,dt = \Gamma(s+1)\sum_{n=0}^{\infty} u^n \sum_{k=0}^{n} \binom{n}{k}(-1)^k \frac{y^k}{(k+x)^{s+1}}$$

and differentiating with respect to $s$ gives us

$$\int_0^1 \frac{t^{x-1}\log^s t \log(-\log t)}{1-u(1-yt)}\,dt = \Gamma(s+1)\sum_{n=0}^{\infty} u^n \sum_{k=0}^{n} \binom{n}{k}(-1)^k \frac{y^k \log(k+x)}{(k+x)^{s+1}}$$

$$+\Gamma'(s+1)\sum_{n=0}^{\infty} u^n \sum_{k=0}^{n} \binom{n}{k}(-1)^k \frac{y^k}{(k+x)^{s+1}}$$

$\square$

Let us now perform an alternative integration by parts of (4.4.9). We have for $x > 1$

$$\int_0^1 t^{x-1}(1-t)^n\,dt = -t^{x-1}\frac{(1-t)^{n+1}}{n+1}\bigg|_0^1 + \frac{(x-1)}{n+1}\int_0^1 t^{x-2}(1-t)^{n+1}\,dt$$

$$= \frac{(x-1)}{n+1}\int_0^1 t^{x-2}(1-t)^{n+1}\,dt$$

Continuing the integration by parts we get for $x > k$

$$\int_0^1 t^{x-1}(1-t)^n\,dt = \frac{(x-1)(x-2)\ldots(x-k)}{(n+1)(n+2)\ldots(n+k)}\int_0^1 t^{x-k-1}(1-t)^{n+k}\,dt$$

$$= \frac{(x-1)(x-2)\ldots(x-k)}{(n+1)(n+2)\ldots(n+k)}\int_0^1 t^{x-k-1}\sum_{j=0}^{n+k}\binom{n+k}{j}(-1)^j t^j\,dt$$

$$= \frac{(x-1)(x-2)\ldots(x-k)}{(n+1)(n+2)\ldots(n+k)}\sum_{j=0}^{n+k}\binom{n+k}{j}\frac{(-1)^j}{x-k+j}$$

Therefore, for $x > n$, we obtain

$$\frac{n!}{x(1+x)\ldots(n+x)} = \frac{(x-1)(x-2)\ldots(x-n)}{(n+1)(n+2)\ldots(2n)}\sum_{j=0}^{2n}\binom{2n}{j}\frac{(-1)^j}{x-n+j}$$

and this gives us

(4.4.14a) $$1 = \frac{x(x^2-1^2)(x^2-2^2)...(x^2-n^2)}{(2n)!} \sum_{j=0}^{2n} \binom{2n}{j} \frac{(-1)^j}{x-n+j}$$

An alternative derivation is shown below.

$$(1-y)^{2n} = \sum_{j=0}^{2n} \binom{2n}{j}(-1)^j y^j$$

$$y^{x-n-1}(1-y)^{2n} = \sum_{j=0}^{2n} \binom{2n}{j}(-1)^j y^{j+x-n-1}$$

$$\int_0^1 y^{x-n-1}(1-y)^{2n} dy = \sum_{j=0}^{2n} \binom{2n}{j} \frac{(-1)^j}{x-n+j}$$

and we see that

$$\int_0^1 y^{x-n-1}(1-y)^{2n} dy = \frac{(2n)!}{(x-n)(x-n+1)...(x-n+2n)}$$

$$= \frac{(2n)!}{x(x^2-1^2)(x^2-2^2)...(x^2-n^2)}$$

Let us now generalise (4.4.9) by letting $n \to y$. We have for $0 < x < 1$ and $y > -1$

$$\int_0^x u^{a-1}(1-u)^y du = \int_0^x u^{a-1} \sum_{k=0}^{\infty} \binom{y}{k}(-1)^k u^k du$$

$$= \sum_{k=0}^{\infty} \binom{y}{k} \frac{(-1)^k}{k+a} x^{a+k}$$

Letting $x \to 1$ and employing Abel's theorem [13, p.245] we obtain

(4.4.14ai) $$\int_0^1 u^{a-1}(1-u)^y du = \sum_{k=0}^{\infty} \binom{y}{k} \frac{(-1)^k}{k+a} = \frac{\Gamma(a)\Gamma(y+1)}{\Gamma(a+y+1)}$$

It may also be noted that

(4.4.14aii) $$\sum_{k=0}^{\infty} \binom{y}{k} \frac{(-1)^k}{k+a} = \sum_{k=0}^{\infty} \frac{(-y)_k (a)_k}{(a+1)_k k!} = \frac{1}{a} {}_2F_1(-y,a;a+1;1)$$

and then employ the Gauss identity for the hypergeometric function.

Differentiation of (4.4.14ai) with respect to $a$ results in

(4.4.14aiv) $\int_0^1 u^{a-1}(1-u)^y \log u \, du = -\sum_{k=0}^{\infty} \binom{y}{k} \frac{(-1)^k}{(k+a)^2} = \frac{\Gamma(a)\Gamma(y+1)}{\Gamma(a+y+1)}[\psi(a)-\psi(a+y+1)]$

(4.4.14av) $\int_0^1 u^{a-1}(1-u)^y \log^s u \, du = (-1)^s s! \sum_{k=0}^{\infty} \binom{y}{k} \frac{(-1)^k}{(k+a)^{s+1}}$

Integration gives us

(4.4.14avi) $\int_0^1 \frac{(u^{t-1}-1)(1-u)^y}{\log u} du = \sum_{k=0}^{\infty} \binom{y}{k}(-1)^k \log\frac{k+t}{k+1}$

We have previously mentioned the Flajolet and Sedgewick paper "Mellin Transforms and Asymptotics: Finite Differences and Rice's Integrals" [68] in (3.13): their paper employed the following lemma:

Let $\varphi(z)$ be analytic in a domain that contains the half line $[n_0, \infty)$. Then the differences of the sequence $\{\varphi(k)\}$ admit the integral representation

(4.4.14b) $\sum_{k=0}^{n} \binom{n}{k}(-1)^k \varphi(k) = \frac{(-1)^n}{2\pi i} \int_C \varphi(z) \frac{n!}{z(z-1)\ldots(z-n)} dz$

where $C$ is a positively oriented closed curve that lies in the domain of analyticity of $\varphi(z)$, encircles $[n_0, n)$ and does not include any of the integers $0, 1, \ldots, n_0 - 1$.

The similarity between the structure of $g(x)$ and the integrand of (4.4.14b) clearly indicates why the Flajolet and Sedgewick paper was relevant in my earlier analysis.

Upon differentiation of the first part of (4.4.13) we have

$$g'(x) = n! \frac{d}{dx}\left[\frac{1}{x(1+x)\ldots(n+x)}\right] = -\sum_{k=0}^{n} \binom{n}{k} \frac{(-1)^k}{(k+x)^2}$$

The Pochhamer symbol is defined by

$$(x)_n = x(1+x)...(n-1+x) = \frac{\Gamma(x+n)}{\Gamma(x)}$$

and differentiation results in

$$\frac{d}{dx}(x)_n = \frac{\Gamma(x)\Gamma'(x+n) - \Gamma(x+n)\Gamma'(x)}{\Gamma^2(x)} = (x)_n [\psi(x+n) - \psi(x)]$$

We see that $g(x) = \frac{n!}{(x)_{n+1}}$ and hence

$$g'(x) = \frac{d}{dx}\frac{n!}{(x)_{n+1}} = -n!\frac{(x)'_{n+1}}{[(x)_{n+1}]^2} = -\frac{n!}{(x)_{n+1}}[\psi(x+n+1) - \psi(x)]$$

Accordingly we get

(4.4.14c) $\quad g'(x) = -\frac{n!\Gamma(x)}{\Gamma(x+n+1)}[\psi(x+n+1) - \psi(x)] = B(x, n+1)[\psi(x+n+1) - \psi(x)]$

$$= -\sum_{k=0}^{n}\binom{n}{k}\frac{(-1)^k}{(k+x)^2}$$

and in particular we have

$$g'(1) = -\frac{n!\Gamma(1)}{\Gamma(n+2)}[\psi(n+2) - \psi(1)] = -\frac{1}{n+1}H_{n+1}$$

Therefore we obtain

$$\sum_{k=0}^{n}\binom{n}{k}\frac{(-1)^k}{(k+1)^2} = \frac{H_{n+1}}{n+1}$$

which we also saw in (4.2.16)

Upon differentiation of (4.4.13) we see that

(4.4.15) $\quad g^{(s-1)}(x) = (-1)^{s-1}(s-1)!\sum_{k=0}^{n}\binom{n}{k}\frac{(-1)^k}{(k+x)^s}$

$$(4.4.16) \qquad = \int_0^1 t^{x-1}(1-t)^n \log^{s-1} t\, dt$$

Now let us consider the series

$$(4.4.17) \qquad H = \sum_{n=1}^{\infty} \frac{g^{(s-1)}(x)}{2^n}$$

$$(4.4.18) \qquad = \sum_{n=1}^{\infty} \frac{1}{2^n} \int_0^1 t^{x-1}(1-t)^n \log^{s-1} t\, dt$$

$$(4.4.19) \qquad = \int_0^1 \frac{t^{x-1}(1-t)\log^{s-1} t}{1+t}\, dt$$

where we completed the summation by employing the geometric series. Using the binomial theorem to expand the denominator of (4.4.19) we obtain

$$(4.4.20) \qquad H = \sum_{k=0}^{\infty} \int_0^1 (-1)^k t^{k+x-1}(1-t)\log^{s-1} t\, dt$$

$$= \sum_{k=0}^{\infty} \int_0^1 (-1)^k t^{k+x-1} \log^{s-1} t\, dt - \sum_{k=0}^{\infty} \int_0^1 (-1)^k t^{k+x} \log^{s-1} t\, dt$$

We have the elementary integral

$$\int_0^1 x^\alpha dx = \frac{1}{\alpha+1}$$

and by parametric differentiation we deduce

$$(4.4.21) \qquad \frac{\partial^{s-1}}{\partial \alpha^{s-1}} \int_0^1 x^\alpha dx = \int_0^1 x^\alpha \log^{s-1} x\, dx = \frac{(-1)^{s-1}(s-1)!}{(\alpha+1)^s}$$

Using (4.4.21) we can write (4.4.20) as

$$(4.4.22) \qquad H = (-1)^{s-1}(s-1)! \left\{ \sum_{k=0}^{\infty} \frac{(-1)^k}{(k+x)^s} - \sum_{k=0}^{\infty} \frac{(-1)^k}{(k+x+1)^s} \right\}$$

Using (4.4.15) we may also express $H$ as

(4.4.23) $$H = (-1)^{s-1}(s-1)!\sum_{n=1}^{\infty}\frac{1}{2^n}\sum_{k=0}^{n}\binom{n}{k}\frac{(-1)^k}{(k+x)^s}$$

and we therefore have an expression involving the alternating Hurwitz-Lerch zeta function (see (4.4.82))

(4.4.24) $$\sum_{n=1}^{\infty}\frac{1}{2^n}\sum_{k=0}^{n}\binom{n}{k}\frac{(-1)^k}{(k+x)^s} = \left\{\sum_{k=0}^{\infty}\frac{(-1)^k}{(k+x)^s} - \sum_{k=0}^{\infty}\frac{(-1)^k}{(k+x+1)^s}\right\}$$

This simplifies to (where the series now starts at $n=0$)

(4.4.24a) $$\sum_{n=0}^{\infty}\frac{1}{2^{n+1}}\sum_{k=0}^{n}\binom{n}{k}\frac{(-1)^k}{(k+x)^s} = \sum_{k=0}^{\infty}\frac{(-1)^k}{(k+x)^s}$$

This identity is also employed in (4.4.112w).

Differentiating (4.4.24a) with respect to $s$ we obtain (it should be noted that we have only proved (4.4.24a) for the case where $s$ is an integer).

(4.4.24aa) $$\sum_{n=0}^{\infty}\frac{1}{2^{n+1}}\sum_{k=0}^{n}\binom{n}{k}\frac{(-1)^k \log(k+x)}{(k+x)^s} = \sum_{k=0}^{\infty}\frac{(-1)^k \log(k+x)}{(k+x)^s}$$

Completing the summation of (4.4.15) we obtain

$$\sum_{n=1}^{\infty}\frac{1}{n^p}\int_0^1 t^{x-1}(1-t)^n \log^{s-1} t\, dt = (-1)^{s-1}(s-1)!\sum_{n=1}^{\infty}\frac{1}{n^s}\sum_{k=0}^{n}\binom{n}{k}\frac{(-1)^k}{(k+x)^s}$$

and hence we get

(4.4.24aaa) $$\int_0^1 t^{x-1} Li_p(1-t)\log^{s-1} t\, dt = (-1)^{s-1}(s-1)!\sum_{n=1}^{\infty}\frac{1}{n^p}\sum_{k=0}^{n}\binom{n}{k}\frac{(-1)^k}{(k+x)^s}$$

We may note that

$$\frac{d^{s-1}}{dx^{s-1}}\int_0^1 t^{x-1} Li_p(1-t)\, dt = \int_0^1 t^{x-1} Li_p(1-t)\log^{s-1} t\, dt$$

The Wolfram Integrator evaluates $\int t^{x-1} Li_2(1-t)\, dt$ in terms of generalised hypergeometric functions. See also (4.4.100qiv).

Let us now integrate equation (4.4.1), noting that it is not defined for $x = 0$

$$\int_\alpha^\beta \frac{n!}{x(1+x)\ldots(n+x)} dx = \int_\alpha^\beta dx \int_0^1 t^{x-1}(1-t)^n dt = \int_\alpha^\beta \sum_{k=0}^n \binom{n}{k} \frac{(-1)^k}{k+x} dx$$

Therefore we obtain for $\beta, \alpha > 0$

(4.4.24b) $$\int_\alpha^\beta dx \int_0^1 t^{x-1}(1-t)^n dt = \int_0^1 (1-t)^n dt \int_\alpha^\beta t^{x-1} dx = \int_0^1 \frac{(t^{\beta-1} - t^{\alpha-1})(1-t)^n}{\log t} dt$$

and

$$\int_\alpha^\beta \sum_{k=0}^n \binom{n}{k} \frac{(-1)^k}{k+x} dx = \sum_{k=0}^n \binom{n}{k} (-1)^k \log \frac{k+\beta}{k+\alpha}$$

Therefore we have

(4.4.24ba) $$\int_0^1 \frac{(t^{\beta-1} - t^{\alpha-1})(1-t)^n}{\log t} dt = \sum_{k=0}^n \binom{n}{k} (-1)^k \log \frac{k+\beta}{k+\alpha}$$

With $n = 0$ we have

(4.4.24c) $$\int_0^1 \frac{(t^{\beta-1} - t^{\alpha-1})}{\log t} dt = \log \frac{\beta}{\alpha}$$

and a version of the above formula is employed in (4.4.112c).

We also have the series for $\left|\frac{1-t}{p}\right| < 1$

$$\sum_{n=0}^\infty \frac{1}{p^n} \int_0^1 \frac{(t^{\beta-1} - t^{\alpha-1})(1-t)^n}{\log t} dt = \sum_{n=0}^\infty \frac{1}{p^n} \sum_{k=0}^n \binom{n}{k} (-1)^k \log \frac{k+\beta}{k+\alpha}$$

and hence after completing the summation we obtain

(4.4.24d) $$\int_0^1 \frac{(t^{\beta-1} - t^{\alpha-1})}{[1-(1-t)/p]\log t} dt = \sum_{n=0}^\infty \frac{1}{p^n} \sum_{k=0}^n \binom{n}{k} (-1)^k \log \frac{k+\beta}{k+\alpha}$$

With $p = 2$ we have

(4.4.24e) $$\int_0^1 \frac{(t^{\beta-1}-t^{\alpha-1})}{(1+t)\log t} dt = \sum_{n=0}^{\infty} \frac{1}{2^{n+1}} \sum_{k=0}^{n} \binom{n}{k}(-1)^k \log\frac{k+\beta}{k+\alpha}$$

As regards the above equation, reference should also be made to (3.36h), (4.4.99c), (4.4.112g) and (4.4.112ga). See also (4.4.112aiii) in Volume III.

Letting $u = 1/p$ in (4.4.24d) we get

(4.4.24f) $$\int_0^1 \frac{(t^{\beta-1}-t^{\alpha-1})}{[1-(1-t)u]\log t} dt = \sum_{n=0}^{\infty} u^n \sum_{k=0}^{n} \binom{n}{k}(-1)^k \log\frac{k+\beta}{k+\alpha}$$

and integrating with respect to $u$ results in

$$\int_0^x du \int_0^1 \frac{(t^{\beta-1}-t^{\alpha-1})}{[1-(1-t)u]\log t} dt = \sum_{n=0}^{\infty} \frac{x^{n+1}}{n+1} \sum_{k=0}^{n} \binom{n}{k}(-1)^k \log\frac{k+\beta}{k+\alpha}$$

Reversing the order of integration gives us

$$\int_0^x du \int_0^1 \frac{(t^{\beta-1}-t^{\alpha-1})}{[1-(1-t)u]\log t} dt = \int_0^1 \frac{(t^{\beta-1}-t^{\alpha-1})}{\log t} dt \int_0^x \frac{1}{[1-(1-t)u]} du$$

$$= -\int_0^1 \frac{(t^{\beta-1}-t^{\alpha-1})\log[1-(1-t)x]}{(1-t)\log t} dt$$

and hence we have

(4.4.24g) $$\int_0^1 \frac{(t^{\beta-1}-t^{\alpha-1})\log[1-(1-t)x]}{(1-t)\log t} dt = \sum_{n=0}^{\infty} \frac{x^{n+1}}{n+1} \sum_{k=0}^{n} \binom{n}{k}(-1)^{k+1} \log\frac{k+\beta}{k+\alpha}$$

Note the close connection between (4.4.24g) and the series given in (4.3.72a) for the digamma function (and refer also to (4.4.24j)).

With $x = 1/2$ this becomes

$$\int_0^1 \frac{(t^{\beta-1}-t^{\alpha-1})\log[(1+t)/2]}{(1-t)\log t} dt = \sum_{n=0}^{\infty} \frac{1}{(n+1)2^{n+1}} \sum_{k=0}^{n} \binom{n}{k}(-1)^{k+1} \log\frac{k+\beta}{k+\alpha}$$

and with $\beta = 2$, $\alpha = 1$ this becomes

(4.4.24h) $$\int_0^1 \frac{\log[(1+t)/2]}{\log t} dt = \sum_{n=0}^{\infty} \frac{1}{(n+1)2^{n+1}} \sum_{k=0}^{n} \binom{n}{k} (-1)^k \log \frac{k+2}{k+1}$$

Letting $x = 1$ in (4.4.24g) we obtain

(4.4.24hi) $$\int_0^1 \frac{(t^{\beta-1} - t^{\alpha-1})}{(1-t)} dt = \sum_{n=0}^{\infty} \frac{1}{n+1} \sum_{k=0}^{n} \binom{n}{k} (-1)^{k+1} \log \frac{k+\beta}{k+\alpha}$$

and with $\beta = 2, \alpha = 1$ this becomes

$$-1 = \sum_{n=0}^{\infty} \frac{1}{n+1} \sum_{k=0}^{n} \binom{n}{k} (-1)^{k+1} \log \frac{k+2}{k+1}$$

We will see in (4.4.92a) that

$$\gamma = \sum_{n=0}^{\infty} \frac{1}{n+1} \sum_{k=0}^{n} \binom{n}{k} (-1)^{k+1} \log(k+1)$$

and hence we have

$$\sum_{n=0}^{\infty} \frac{1}{n+1} \sum_{k=0}^{n} \binom{n}{k} (-1)^{k+1} \log(k+2) = \gamma - 1$$

Since $\frac{t^N - 1}{1-t} = -\sum_{j=1}^{N-1} t^j$, we have $\int_0^1 \frac{t^N - 1}{1-t} dt = -H_N$ and therefore

$$\int_0^1 \frac{t^N - 1}{1-t} dt = -H_N = \sum_{n=0}^{\infty} \frac{1}{n+1} \sum_{k=0}^{n} \binom{n}{k} (-1)^{k+1} \log \frac{k+N+1}{k+1}$$

This then gives us

(4.4.24i) $$\sum_{n=0}^{\infty} \frac{1}{n+1} \sum_{k=0}^{n} \binom{n}{k} (-1)^{k+1} \log(k+N+1) = \gamma - H_N$$

We may generalise this by noting that

$$\int_0^1 \frac{(t^{\beta-1} - t^{\alpha-1})}{1-t} dt = \int_0^1 \frac{1-t^{\alpha-1}}{1-t} dt - \int_0^1 \frac{1-t^{\beta-1}}{1-t} dt = \psi(\alpha) - \psi(\beta)$$

and hence from (4.4.24hi) we have again derived (4.3.72a)

(4.4.24j) $$\sum_{n=0}^{\infty}\frac{1}{n+1}\sum_{k=0}^{n}\binom{n}{k}(-1)^{k+1}\log\frac{k+\beta}{k+\alpha}=\psi(\alpha)-\psi(\beta)$$

Using (4.4.92a) this gives us another derivation of (4.3.74)

$$\sum_{n=0}^{\infty}\frac{1}{n+1}\sum_{k=0}^{n}\binom{n}{k}(-1)^{k}\log(k+\alpha)=\psi(\alpha)$$

Integrating (4.4.24g) gives us

$$\int_{0}^{1}\frac{(t^{\beta-1}-t^{\alpha-1})\log[1-(1-t)u]}{(1-t)\log t}\left(u-\frac{1}{1-t}\right)dt-u\int_{0}^{1}\frac{(t^{\beta-1}-t^{\alpha-1})}{(1-t)\log t}dt=$$

$$\sum_{n=0}^{\infty}\frac{u^{n+2}}{(n+1)(n+2)}\sum_{k=0}^{n}\binom{n}{k}(-1)^{k+1}\log\frac{k+\beta}{k+\alpha}$$

and, using (4.4.24g), this may be written as

$$-\int_{0}^{1}\frac{(t^{\beta-1}-t^{\alpha-1})\log[1-(1-t)u]}{(1-t)^{2}\log t}dt=$$

$$\sum_{n=0}^{\infty}\frac{u^{n+2}}{(n+1)(n+2)}\sum_{k=0}^{n}\binom{n}{k}(-1)^{k+1}\log\frac{k+\beta}{k+\alpha}-\sum_{n=0}^{\infty}\frac{u^{n+2}}{n+1}\sum_{k=0}^{n}\binom{n}{k}(-1)^{k+1}\log\frac{k+\beta}{k+\alpha}+\int_{0}^{1}\frac{(t^{\beta-1}-t^{\alpha-1})}{(1-t)\log t}dt$$

Differentiating (4.4.24g) with respect to $\beta$ we get

(4.4.24k) $$\int_{0}^{1}\frac{t^{\beta-1}\log[1-(1-t)x]}{1-t}dt=\sum_{n=0}^{\infty}\frac{x^{n+1}}{n+1}\sum_{k=0}^{n}\binom{n}{k}\frac{(-1)^{k+1}}{k+\beta}$$

With $\beta = N+1$ and $x=1$ we have

$$\int_{0}^{1}\frac{t^{N}\log t}{(1-t)}dt=\sum_{n=0}^{\infty}\frac{1}{n+1}\sum_{k=0}^{n}\binom{n}{k}\frac{(-1)^{k+1}}{k+N+1}$$

We also have from (4.4.238aa)

$$\int \frac{t^N \log t}{(1-t)} dt = \sum_{j=1}^{N} \frac{t^j}{j^2} - \log t \sum_{j=1}^{N} \frac{t}{j} + Li_2(1-t)$$

and hence

$$\int_0^1 \frac{t^N \log t}{(1-t)} dt = \sum_{j=1}^{N} \frac{1}{j^2} - Li_2(1) = H_N^{(2)} - \varsigma(2)$$

Therefore we have shown that

(4.4.24l) $$H_N^{(2)} - \varsigma(2) = \sum_{n=0}^{\infty} \frac{1}{n+1} \sum_{k=0}^{n} \binom{n}{k} \frac{(-1)^{k+1}}{k+N+1}$$

and, for example, with $N=1$ we have

(4.4.24m) $$1 - \varsigma(2) = \sum_{n=0}^{\infty} \frac{1}{n+1} \sum_{k=0}^{n} \binom{n}{k} \frac{(-1)^{k+1}}{k+2}$$

Differentiating (4.4.24k) with respect to $\beta$ we get

(4.4.24n) $$\int_0^1 \frac{t^{\beta-1} \log t \log[1-(1-t)x]}{(1-t)} dt = \sum_{n=0}^{\infty} \frac{x^{n+1}}{n+1} \sum_{k=0}^{n} \binom{n}{k} \frac{(-1)^k}{(k+\beta)^2}$$

With $\beta = N+1$ and $x=1$ we get

$$\int_0^1 \frac{t^N \log^2 t}{(1-t)} dt = \sum_{n=0}^{\infty} \frac{1}{n+1} \sum_{k=0}^{n} \binom{n}{k} \frac{(-1)^k}{(k+N+1)^2}$$

As will be seen in equation (4.4.235) of Volume IV we have

$$\int \frac{x^k \log^2 x}{1-x} dx = -2 \sum_{j=1}^{k} \frac{x^j}{j^3} + 2 \log x \sum_{j=1}^{k} \frac{x^j}{j^2} - \log^2 x \sum_{j=1}^{k} \frac{x^j}{j} - \log(1-x) \log^2 x - 2 \log x \, Li_2(x) + 2Li_3(x)$$

and more specifically

$$\int_0^1 \frac{t^N \log^2 t}{(1-t)} dt = -2H_N^{(2)} + 2\varsigma(3)$$

(this is also valid for $N=0$, where we take $H_0^{(2)} = 0$).

Hence we get

(4.4.24o) $$\sum_{n=0}^{\infty}\frac{1}{n+1}\sum_{k=0}^{n}\binom{n}{k}\frac{(-1)^k}{(k+N+1)^2}=2\left[\varsigma(3)-H_N^{(2)}\right]$$

Differentiating (4.4.24n) $p$ times we have

$$\int_0^1\frac{t^{\beta-1}\log^{p+1}t\log[1-(1-t)x]}{(1-t)}dt=(-1)^p(p+1)!\sum_{n=0}^{\infty}\frac{x^{n+1}}{n+1}\sum_{k=0}^{n}\binom{n}{k}\frac{(-1)^k}{(k+\beta)^{p+2}}$$

and with $x=1$ we get

$$\int_0^1\frac{t^{\beta-1}\log^{p+2}t}{(1-t)}dt=(-1)^p(p+1)!\sum_{n=0}^{\infty}\frac{1}{n+1}\sum_{k=0}^{n}\binom{n}{k}\frac{(-1)^k}{(k+\beta)^{p+2}}$$

Differentiating (4.3.68) we see that

$$\int_0^1\frac{t^{\beta-1}\log^{p+2}t}{(1-t)}dt=-\psi^{(p+2)}(\beta)$$

and hence we have

(4.4.24p) $$\psi^{(p+2)}(\beta)=(-1)^{p+1}(p+1)!\sum_{n=0}^{\infty}\frac{1}{n+1}\sum_{k=0}^{n}\binom{n}{k}\frac{(-1)^k}{(k+\beta)^{p+2}}$$

We know from [126, p.22] that

$$\psi^{(n)}(z)=(-1)^{n+1}n!\sum_{k=0}^{\infty}\frac{1}{(k+z)^{n+1}}=(-1)^{n+1}n!\varsigma(n+1,z)$$

and therefore we see that

$$\psi^{(p+2)}(\beta)=(-1)^{p+1}(p+1)!\sum_{n=0}^{\infty}\frac{1}{n+1}\sum_{k=0}^{n}\binom{n}{k}\frac{(-1)^k}{(k+\beta)^{p+2}}=(-1)^{p+1}(p+2)!\varsigma(p+3,\beta)$$

Hence we get

(4.4.24pi) $$\varsigma(p+2,\beta)=\frac{1}{p+1}\sum_{n=0}^{\infty}\frac{1}{n+1}\sum_{k=0}^{n}\binom{n}{k}\frac{(-1)^k}{(k+\beta)^{p+2}}$$

which is the formula (3.12a) originally derived by Hasse.

As will be noted in (4.4.135a), Larcombe et al. [95] showed that

$$m\binom{m+n}{n}\sum_{k=0}^{n}\binom{n}{k}\frac{(-1)^k}{(m+k)^2} = \sum_{k=m}^{m+n}\frac{1}{k}$$

and hence we have

$$\sum_{n=0}^{\infty}\frac{1}{n+1}\sum_{k=0}^{n}\binom{n}{k}\frac{(-1)^k}{(k+N+1)^2} = \frac{1}{N+1}\sum_{n=0}^{\infty}\frac{1}{n+1}\binom{N+1+n}{n}^{-1}\sum_{k=N+1}^{N+1+n}\frac{1}{k} = 2\left[\varsigma(3) - H_N^{(2)}\right]$$

With $N = 1$ we obtain

$$\sum_{n=0}^{\infty}\frac{1}{n+1}\sum_{k=0}^{n}\binom{n}{k}\frac{(-1)^k}{(k+2)^2} = \frac{1}{2}\sum_{n=0}^{\infty}\frac{1}{n+1}\binom{n+2}{n}^{-1}\sum_{k=2}^{n+2}\frac{1}{k}$$

and the right-hand side is equivalent to

$$\frac{1}{2}\sum_{n=0}^{\infty}\frac{1}{n+1}\binom{n+2}{n}^{-1}\sum_{k=2}^{n+2}\frac{1}{k} = \sum_{n=0}^{\infty}\frac{H_{n+2}-1}{(n+1)^2(n+2)}$$

This then gives us

(4.4.24q) $$\sum_{n=0}^{\infty}\frac{H_{n+2}-1}{(n+1)^2(n+2)} = 2\left[\varsigma(3) - 1\right]$$

which, from a structural point of view, initially appears to be unexpected. However, structurally similar series have been found by Ogreid and Osland (see [105(ii)] and [105(iii)]): for example

$$\sum_{n=1}^{\infty}\frac{H_n}{n^2(n+1)} = 2\varsigma(3) - \varsigma(2)$$

$$\sum_{n=1}^{\infty}\frac{H_n}{n(n+1)} = \varsigma(2)$$

Similar series are also given in the 2005 paper "Six Gluon Open Superstring Disk Amplitude, Multiple Hypergeometric Series and Euler-Zagier Sums" by Oprisa and Stieberger [105(iv), p.54].

Using the Larcombe identity (4.4.135) we get another rather unexpected result after making the substitution in (4.4.24k)

(4.4.24r) $$\varsigma(2)-1 = \sum_{n=0}^{\infty}\frac{1}{n+1}\sum_{k=0}^{n}\binom{n}{k}\frac{(-1)^k}{k+2} = \frac{1}{2}\sum_{n=0}^{\infty}\frac{1}{(n+1)^2(n+2)}$$

Integrating (4.4.24n) with respect to $x$ we get

(4.4.24s) $$\int_0^u dx\int_0^1 \frac{t^{\beta-1}\log t \log[1-(1-t)x]}{1-t}dt = \sum_{n=0}^{\infty}\frac{u^{n+2}}{(n+1)(n+2)}\sum_{k=0}^{n}\binom{n}{k}\frac{(-1)^k}{(k+\beta)^2}$$

and therefore we obtain

(4.4.24t) $$u\int_0^1 \frac{t^{\beta-1}\log t \log[1-(1-t)u]}{1-t}dt - \int_0^1 \frac{t^{\beta-1}\log t \log[1-(1-t)u]}{(1-t)^2}dt - u\int_0^1 \frac{t^{\beta-1}\log t}{1-t}dt$$

$$= \sum_{n=0}^{\infty}\frac{u^{n+2}}{(n+1)(n+2)}\sum_{k=0}^{n}\binom{n}{k}\frac{(-1)^k}{(k+\beta)^2}$$

Dividing (4.4.24n) by $x$ and then integrating with respect to $x$ we get

(4.4.24u) $$\int_0^1 \frac{t^{\beta-1}\log t\, Li_2[(1-t)u]}{(1-t)}dt = -\sum_{n=0}^{\infty}\frac{u^n}{(n+1)^2}\sum_{k=0}^{n}\binom{n}{k}\frac{(-1)^k}{(k+\beta)^2}$$

With $\beta=1$ and $u=1$ we get

$$\int_0^1 \frac{\log t\, Li_2(1-t)}{1-t}dt = -\sum_{n=0}^{\infty}\frac{1}{(n+1)^2}\sum_{k=0}^{n}\binom{n}{k}\frac{(-1)^k}{(k+1)^2}$$

Since

$$\int_0^1 \frac{\log t\, Li_2(1-t)}{(1-t)}dt = \frac{1}{2}[Li_2(1-t)]^2\Big|_0^1 = -\frac{1}{2}\varsigma^2(2)$$

we see that

(4.4.24v) $$\sum_{n=0}^{\infty}\frac{1}{(n+1)^2}\sum_{k=0}^{n}\binom{n}{k}\frac{(-1)^k}{(k+1)^2} = \frac{1}{2}\varsigma^2(2)$$

From Larcombe's formula we see that

$$\sum_{k=0}^{n}\binom{n}{k}\frac{(-1)^k}{(k+1)^2}=\frac{H_{n+1}}{n+1}$$

and hence we have

$$\sum_{n=0}^{\infty}\frac{1}{(n+1)^2}\sum_{k=0}^{n}\binom{n}{k}\frac{(-1)^k}{(k+1)^2}=\sum_{n=0}^{\infty}\frac{H_{n+1}}{(n+1)^3}$$

which gives us

(4.4.24vi) $\quad\displaystyle\sum_{n=1}^{\infty}\frac{H_n}{n^3}=\frac{1}{2}\varsigma^2(2)$

Similarly, integrating (4.4.24u) with respect to $u$ we get

(4.4.24w) $\quad\displaystyle\int_{0}^{u}\frac{du}{u}\int_{0}^{1}\frac{t^{\beta-1}\log t\, Li_2[(1-t)u]}{(1-t)}dt=-\sum_{n=0}^{\infty}\frac{u^n}{(n+1)^3}\sum_{k=0}^{n}\binom{n}{k}\frac{(-1)^k}{(k+\beta)^2}$

and therefore we have

$$\int_{0}^{1}\frac{t^{\beta-1}\log t\, Li_3[(1-t)u]}{(1-t)}dt=-\sum_{n=0}^{\infty}\frac{u^n}{(n+1)^3}\sum_{k=0}^{n}\binom{n}{k}\frac{(-1)^k}{(k+\beta)^2}$$

With $\beta=1$ and $u=1$ we get

(4.4.24x) $\quad\displaystyle\int_{0}^{1}\frac{\log t\, Li_3(1-t)}{1-t}dt=-\sum_{n=0}^{\infty}\frac{1}{(n+1)^3}\sum_{k=0}^{n}\binom{n}{k}\frac{(-1)^k}{(k+1)^2}=-\sum_{n=1}^{\infty}\frac{H_n}{n^4}$

It is clear that this may be generalised to

(4.4.24y) $\quad\displaystyle\int_{0}^{1}\frac{t^{\beta-1}\log t\, Li_p[(1-t)u]}{1-t}dt=-\sum_{n=0}^{\infty}\frac{u^n}{(n+1)^p}\sum_{k=0}^{n}\binom{n}{k}\frac{(-1)^k}{(k+\beta)^2}$

and differentiation of (4.4.24y) with respect to $\beta$ results in

(4.4.24z) $\quad\displaystyle\int_{0}^{1}\frac{t^{\beta-1}\log^{q+1} t\, Li_p[(1-t)u]}{(1-t)}dt=(q+1)!(-1)^{q+1}\sum_{n=0}^{\infty}\frac{u^n}{(n+1)^p}\sum_{k=0}^{n}\binom{n}{k}\frac{(-1)^k}{(k+\beta)^{q+2}}$

With $\beta = 1$, $q = 1$, $u = 1$ and $p = 2$ we get (see (4.4.64ci) et seq. in Volume III)

$$\int_0^1 \frac{\log^2 t \, Li_2(1-t)}{(1-t)} dt = 2 \sum_{n=0}^{\infty} \frac{1}{(n+1)^2} \sum_{k=0}^{n} \binom{n}{k} \frac{(-1)^k}{(k+1)^3}$$

and integration by parts gives us

$$\int \frac{\log^2 t \, Li_2(1-t)}{(1-t)} dt = \int \log t \, \frac{\log t \, Li_2(1-t)}{(1-t)} dt$$

$$= \frac{1}{2} \log t \left[Li_2(1-t)\right]^2 - \frac{1}{2} \int \frac{\left[Li_2(1-t)\right]^2}{t} dt$$

Therefore we have

$$\int_a^1 \frac{\log^2 t \, Li_2(1-t)}{(1-t)} dt = \frac{1}{2} \log a \left[Li_2(1-a)\right]^2 - \frac{1}{2} \int_a^1 \frac{\varsigma^2(2)}{t} dt - \frac{1}{2} \int_a^1 \frac{\left[Li_2(1-t)\right]^2}{t} dt + \frac{1}{2} \int_a^1 \frac{\varsigma^2(2)}{t} dt$$

$$= \frac{1}{2} \left[Li_2(1-a)\right]^2 \log a - \frac{1}{2} \varsigma^2(2) \log a - \frac{1}{2} \int_a^1 \frac{\left[Li_2(1-t)\right]^2 - \varsigma^2(2)}{t} dt$$

and, in the limit as $a \to 0$, we get

$$\int_0^1 \frac{\log^2 t \, Li_2(1-t)}{(1-t)} dt = \frac{1}{2} \int_0^1 \frac{\varsigma^2(2) - \left[Li_2(1-t)\right]^2}{t} dt$$

In passing we note that if $\lim_{a \to o} f'(a)$ is finite then

$$\lim_{a \to o} f(a) \log a = \lim_{a \to o} \frac{f(a)}{a} a \log a$$

$$\lim_{a \to o} \frac{f(a)}{a} = \lim_{a \to o} \frac{f'(a)}{1}$$

and therefore $\lim_{a \to o} f(a) \log a = 0$ if $f'(0)$ is finite.

Hence we see that (see (4.4.64ci) et seq. in Volume III)

(4.4.24zi) $$\int_0^1 \frac{\varsigma^2(2) - \left[Li_2(1-t)\right]^2}{t} dt = 4 \sum_{n=0}^{\infty} \frac{1}{(n+1)^2} \sum_{k=0}^{n} \binom{n}{k} \frac{(-1)^k}{(k+1)^3}$$

We note from the Larcombe identity (4.4.135b) that

$$2m\binom{m+n}{n}\sum_{k=0}^{n}\binom{n}{k}\frac{(-1)^k}{(m+k)^3} = \left(\sum_{k=m}^{m+n}\frac{1}{k}\right)^2 + \sum_{k=m}^{m+n}\frac{1}{k^2}$$

and thus we may write (4.4.24zi) in terms involving the generalised harmonic numbers.

We have

$$\int_0^1 \frac{\varsigma^2(2)-[Li_2(1-t)]^2}{t}\,dt = \int_0^1 \frac{(\varsigma(2)+[Li_2(1-t)])(\varsigma(2)-[Li_2(1-t)])}{t}\,dt$$

$$\int \frac{(\varsigma(2)+[Li_2(1-t)])}{t}\,dt = \varsigma(2)\log t + Li_2(1-t)\log t + \log(1-t)\log^2 t + 2Li_2(t)\log t - 2Li_3(t)$$

and with integration by parts we obtain

$$\int_0^1 \frac{\varsigma^2(2)-[Li_2(1-t)]^2}{t}\,dt =$$

$$(\varsigma(2)-[Li_2(1-t)])\Big[\varsigma(2)\log t + Li_2(1-t)\log t + \log(1-t)\log^2 t + 2Li_2(t)\log t - 2Li_3(t)\Big]\Big|_0^1$$

$$-\int_0^1 \Big[\varsigma(2)\log t + Li_2(1-t)\log t + \log(1-t)\log^2 t + 2Li_2(t)\log t - 2Li_3(t)\Big]\frac{\log t}{1-t}\,dt$$

$$= -2\varsigma(2)\varsigma(3) - \int_0^1 \Big[\varsigma(2)\log t + Li_2(1-t)\log t + \log(1-t)\log^2 t + 2Li_2(t)\log t - 2Li_3(t)\Big]\frac{\log t}{1-t}\,dt$$

We have

$$\int \frac{\log^2 t}{1-t}\,dt = -\log(1-t)\log^2 t - 2Li_2(t)\log t + 2Li_3(t)$$

but the other integrals are more intractable. I considered trying Euler's identity (1.6c)

$$\varsigma(2) = \log x \log(1-x) + Li_2(x) + Li_2(1-x)$$

which gives us

$$[Li_2(1-x)]^2 = \varsigma^2(2) - 2\varsigma(2)\log x \log(1-x) - 2\varsigma(2)Li_2(x)$$

$$+ \log^2 x \log^2(1-x) - 2\log x \log(1-x)Li_2(x) + [Li_2(x)]^2$$

However, the Wolfram Integrator was unable to evaluate the resulting integrals.

From (4.4.38b) we see that

$$Li_p[(1-t)u] = \frac{(1-t)u(-1)^{p-1}}{\Gamma(p)} \int_0^1 \frac{\log^{p-1} y}{(1-(1-t)uy)} dy$$

and substitution in (4.4.24z) gives us

$$\int_0^1 \frac{t^{\beta-1} \log^{q+1} t \, Li_p[(1-t)u]}{(1-t)} dt = \frac{u(-1)^{p-1}}{\Gamma(p)} \int_0^1 \frac{t^{\beta-1} \log^{q+1} t}{(1-t)} dt \int_0^1 \frac{(1-t)\log^{p-1} y}{1-(1-t)uy} dy$$

$$= \frac{u(-1)^{p-1}}{\Gamma(p)} \int_0^1 \int_0^1 t^{\beta-1} \log^{q+1} t \frac{\log^{p-1} y}{1-(1-t)uy} dt\, dy$$

For example, we now let $\beta = 1$, $u = 1$ and $p = 3$ to get

$$\int_0^1 \frac{\log^2 y}{1-Ay} dy = \frac{1}{A}\left[\log(1-Ax)\log^2 x + 2Li_2(Ax)\log x - 2Li_3(Ax)\right]\Big|_0^1$$

$$= -\frac{2}{A} Li_3(A)$$

However, I believe that the analysis is becoming rather circular at this stage! Having reached the end of the alphabet, I will leave further consideration of this to another day!

The following result is well-known.

**(ii) Theorem 4.2:**

(4.4.25) $$Li_s(x) = \frac{x}{\Gamma(s)} \int_0^\infty \frac{u^{s-1}}{e^u - x} du$$

**Proof:**

The gamma function is defined in (4.3.1) as

(4.4.26) $$\Gamma(s) = \int_0^\infty t^{s-1} e^{-t} dt \quad , \operatorname{Re}(s) > 0$$

and, using the substitution $t = ku$, we obtain

(4.4.27) $$\Gamma(s) = k^s \int_0^\infty u^{s-1} e^{-ku} du$$

Hence we have

(4.4.28) $$\frac{1}{k^s} = \frac{1}{\Gamma(s)} \int_0^\infty u^{s-1} e^{-ku} du$$

We now consider the finite sum set out below

(4.4.29) $$S_n(x) = \sum_{k=1}^{n} \binom{n}{k} \frac{x^k}{k^s}$$

and combine (4.4.28) and (4.4.29) to obtain

(4.4.30) $$S_n(x) = \sum_{k=1}^{n} \binom{n}{k} \frac{x^k}{k^s} = \sum_{k=1}^{n} \binom{n}{k} x^k \cdot \frac{1}{\Gamma(s)} \int_0^\infty u^{s-1} e^{-ku} du$$

(4.4.31) $$= \frac{1}{\Gamma(s)} \int_0^\infty u^{s-1} \sum_{k=1}^{n} \left\{ \binom{n}{k} e^{-ku} x^k \right\} du$$

(4.4.32) $$= \frac{1}{\Gamma(s)} \int_0^\infty u^{s-1} \sum_{k=1}^{n} \left\{ \binom{n}{k} \left[ e^{-u} x \right]^k \right\} du$$

From the binomial theorem we have

$$\sum_{k=1}^{n} \binom{n}{k} A^k = (1+A)^n - 1$$

and hence (4.4.32) becomes

(4.4.33) $$S_n(x) = \sum_{k=1}^{n} \binom{n}{k} \frac{x^k}{k^s} = \frac{1}{\Gamma(s)} \int_0^\infty u^{s-1} \left\{ (1 + xe^{-u})^n - 1 \right\} du$$

We now revisit the expression for the polylogarithm derived in (3.55)

$$(4.4.34) \qquad 2Li_s(x) = \sum_{n=1}^{\infty} \frac{1}{2^n} \sum_{k=1}^{n} \binom{n}{k} \frac{x^k}{k^s}$$

and substitute the relation obtained for $S_n(x)$ to deduce

$$(4.4.35) \qquad 2Li_s(x) = \sum_{n=1}^{\infty} \frac{1}{2^n} \frac{1}{\Gamma(s)} \int_0^{\infty} \sum_{k=1}^{\infty} u^{s-1} \left\{ \left(1 + xe^{-u}\right)^n - 1 \right\} du$$

$$(4.4.36) \qquad = \frac{1}{\Gamma(s)} \int_0^{\infty} u^{s-1} \sum_{n=1}^{\infty} \frac{\left\{\left(1 + xe^{-u}\right)^n - 1\right\}}{2^n} du$$

(where we have assumed that interchanging the order of summation and integration is permissible). Use of the geometric series simplifies (4.4.36) to

$$(4.4.37) \qquad Li_s(x) = \frac{x}{\Gamma(s)} \int_0^{\infty} \frac{u^{s-1} e^{-u}}{1 - xe^{-u}} du = \frac{x}{\Gamma(s)} \int_0^{\infty} \frac{u^{s-1}}{e^u - x} du$$

After deriving the above formula, I subsequently found out that it had been discovered by Appell many years ago: equation (4.4.37) was in fact posed as an exercise in Whittaker and Watson's book [135, p.280] and, inter alia, is also reported in [126, p.121].

By letting $x = 1$ in (4.4.37) we obtain the familiar formula [25, p.222]

$$(4.4.38) \qquad Li_s(1) = \varsigma(s) = \frac{1}{\Gamma(s)} \int_0^{\infty} \frac{u^{s-1}}{e^u - 1} du$$

This simpler expression can be derived directly in the following manner. Using (4.4.28) and making the summation we have

$$\sum_{k=1}^{\infty} \frac{1}{k^s} = \frac{1}{\Gamma(s)} \sum_{k=1}^{\infty} \int_0^{\infty} u^{s-1} e^{-ku} du$$

$$= \frac{1}{\Gamma(s)} \int_0^{\infty} u^{s-1} \sum_{k=1}^{\infty} e^{-ku} du$$

After noting that $e^{-u} < 1 \ \forall \ u \in (0, \infty)$, we sum the geometric series and obtain (4.4.38). This identity is used in the proof of the functional equation for the Riemann zeta function in Appendix F of Volume VI.

When I started work on this in 2003 I was concentrating on the definition of the polylogarithm in (4.4.34) and I completely missed the rather more obvious approach outlined below.

$$\sum_{k=1}^{\infty} \frac{x^k}{k^s} = \frac{1}{\Gamma(s)} \int_0^{\infty} \sum_{k=1}^{\infty} u^{s-1} \left(xe^{-u}\right)^k du = \frac{x}{\Gamma(s)} \int_0^{\infty} \frac{u^{s-1}}{e^u - x} du$$

Making the substitution $u = \log t$ in (4.4.37) we have

(4.4.38a) $$Li_s(x) = \frac{x}{\Gamma(s)} \int_1^{\infty} \frac{\log^{s-1} t}{(t-x)t} dt$$

and the substitution $t = 1/y$ gives us

(4.4.38b) $$Li_s(x) = \frac{x(-1)^{s-1}}{\Gamma(s)} \int_0^1 \frac{\log^{s-1} y}{(1-xy)} dy$$

We also note the double integrals

(4.4.38c) $$Li_{s+1}(u) = \int_0^u \frac{Li_s(x)}{x} dx = \frac{(-1)^{s-1}}{\Gamma(s)} \int_0^u \int_0^1 \frac{\log^{s-1} y}{1-xy} dy dx$$

(4.4.38d) $$\varsigma(s+1) = \int_0^1 \frac{Li_s(x)}{x} dx = \frac{(-1)^{s-1}}{\Gamma(s)} \int_0^1 \int_0^1 \frac{\log^{s-1} y}{1-xy} dy dx$$

and after carrying out the integration on $x$ we get

(4.4.38e) $$Li_{s+1}(u) = \frac{(-1)^s}{\Gamma(s)} \int_0^1 \frac{\log^{s-1} y \log(1-uy)}{y} dy$$

(4.4.38f) $$\varsigma(s+1) = \frac{(-1)^s}{\Gamma(s)} \int_0^1 \frac{\log^{s-1} y \log(1-y)}{y} dy$$

(and in this regard, see (4.4.91ga)).

Letting $x = 1$ in (4.4.38b) we get

(4.4.39) $$Li_s(1) = \frac{(-1)^{s-1}}{\Gamma(s)} \int_0^1 \frac{\log^{s-1} y}{1-y} dy$$

and hence

(4.4.40) $$\varsigma(s) = \frac{(-1)^{s-1}}{\Gamma(s)} \int_0^1 \frac{\log^{s-1} y}{1-y} dy$$

Letting $x = -1$ in (4.4.38b) we have

(4.4.41) $$Li_s(-1) = \frac{(-1)^{s-2}}{\Gamma(s)} \int_0^1 \frac{\log^{s-1} y}{1+y} dy$$

and hence

(4.4.42) $$\varsigma_a(s) = \frac{(-1)^{s-1}}{\Gamma(s)} \int_0^1 \frac{\log^{s-1} y}{1+y} dy = \frac{1}{\Gamma(s)} \int_0^\infty \frac{u^{s-1}}{e^u + 1} du$$

(equation (4.4.42) was employed in (3.83) as the starting point of Amore's variational method). See also (3.86j).

From (4.4.38e) we obtain

$$\int_0^t \frac{Li_{s+1}(u)}{u} du = \frac{(-1)^s}{\Gamma(s)} \int_0^t \frac{du}{u} \int_0^1 \frac{\log^{s-1} y \log(1-uy)}{y} dy$$

and hence we get

(4.4.42i) $$Li_{s+2}(t) = \frac{(-1)^{s+1}}{\Gamma(s)} \int_0^1 \frac{\log^{s-1} y \, Li_2(ty)}{y} dy$$

Integration by parts gives the following for $s = 2$

$$\log y \, Li_3(ty) - Li_4(ty) = \int \frac{\log y \, Li_2(ty)}{y} dy$$

which concurs with the above. Equation (4.4.42i) is employed in Appendix E of Volume VI.

We also note that

$$\int_0^t \frac{Li_s(x)}{x}\,dx = \frac{1}{\Gamma(s)} \int_0^t dx \int_0^\infty \frac{u^{s-1}}{e^u - x}\,du$$

$$= \frac{1}{\Gamma(s)} \int_0^\infty u^{s-1}[u - \log(e^u - t)]\,du$$

and we then obtain

$$Li_{s+1}(t) = -\frac{1}{\Gamma(s)} \int_0^\infty u^{s-1} \log(1 - te^{-u})\,du$$

Differentiating (4.4.37) with respect to $s$ we obtain

$$\frac{\partial}{\partial s} Li_s(x) = -x \frac{\psi(s)}{\Gamma(s)} \int_0^\infty \frac{u^{s-1}}{e^u - x}\,du + \frac{x}{\Gamma(s)} \int_0^\infty \frac{u^{s-1} \log u}{e^u - x}\,du$$

and therefore

$$-\sum_{n=1}^\infty \frac{x^n \log n}{n^s} = -Li_s(x)\psi(s) + \frac{x}{\Gamma(s)} \int_0^\infty \frac{u^{s-1} \log u}{e^u - x}\,du$$

Letting $x = 1$ we obtain

$$\int_0^\infty \frac{u^{s-1} \log u}{e^u - 1}\,du = \Gamma(s)\varsigma'(s) + \varsigma(s)\Gamma(s)\psi(s)$$

(4.4.42a)
$$= \Gamma(s)\varsigma'(s) + \varsigma(s)\Gamma'(s) = \frac{d}{ds}[\Gamma(s)\varsigma(s)]$$

With the substitution $u = \log t$ we have

$$\int_0^\infty \frac{u^{s-1} \log u}{e^u - 1}\,du = \int_1^\infty \frac{\log^{s-1} t \log \log t}{t(t-1)}\,dt = \Gamma(s)\varsigma'(s) + \varsigma(s)\Gamma'(s)$$

and with $t = 1/y$ this becomes

$$(-1)^s \int_0^1 \frac{\log^{s-1} y \log \log(1/y)}{1 - y}\,dy = \Gamma(s)\varsigma'(s) + \varsigma(s)\Gamma'(s)$$

In [6a] Adamchik reports the following integral

(4.4.42b) $$\int_0^\infty \frac{u \log u}{e^{2\pi u} - 1} du = \frac{1}{2}\varsigma'(-1)$$

and making the substitution $x = 2\pi u$ this becomes

$$\int_0^\infty \frac{u \log u}{e^{2\pi u} - 1} du = \frac{1}{4\pi^2} \int_0^\infty \frac{x \log x}{e^x - 1} dx - \frac{\log 2\pi}{4\pi^2} \int_0^\infty \frac{x}{e^x - 1} dx$$

Reference to (4.4.37) shows that

$$\int_0^\infty \frac{x}{e^x - 1} dx = \varsigma(2)$$

From (4.4.42a) above we have

(4.4.42bi) $$\int_0^\infty \frac{x \log x}{e^x - 1} dx = \Gamma(2)\varsigma'(2) + \varsigma(2)\Gamma'(2)$$

Therefore, since $\Gamma(2) = 1$, we get

$$\int_0^\infty \frac{u \log u}{e^{2\pi u} - 1} du = \frac{1}{4\pi^2}[\varsigma'(2) + \varsigma(2)\Gamma'(2)] - \frac{\log 2\pi}{4\pi^2}\varsigma(2)$$

In (E.20a) of Volume VI it is shown that $\Gamma'(2) = 1 - \gamma$ and hence we get

$$\int_0^\infty \frac{u \log u}{e^{2\pi u} - 1} du = \frac{1}{4\pi^2}[\varsigma'(2) + (1-\gamma)\varsigma(2)] - \frac{\log 2\pi}{24}$$

Using the functional equation for $\varsigma(s)$ in (F.7) of Volume VI it is shown that

$$\varsigma'(-1) = \frac{1}{12}(1 - \gamma - \log 2\pi) + \frac{1}{2\pi^2}\varsigma'(2)$$

and combining these equations we obtain the same result as Adamchik

$$\int_0^\infty \frac{u \log u}{e^{2\pi u} - 1} du = \frac{1}{2}\varsigma'(-1)$$

We also have from (4.4.25)

$$Li_s(-x) = -\frac{x}{\Gamma(s)} \int_0^\infty \frac{u^{s-1}}{e^u + x} du$$

and therefore we get

$$Li_s(-1) = -\frac{1}{\Gamma(s)} \int_0^\infty \frac{u^{s-1}}{e^u + 1} du$$

and, since $\varsigma_a(s) = -Li_s(-1)$, this is equivalent to

$$\varsigma_a(s) = \frac{1}{\Gamma(s)} \int_0^\infty \frac{u^{s-1}}{e^u + 1} du$$

Differentiating with respect to $s$ we obtain

(4.4.42c) $$\varsigma_a'(s) = -\frac{\Gamma'(s)}{[\Gamma(s)]^2} \int_0^\infty \frac{u^{s-1}}{e^u + 1} du + \frac{1}{\Gamma(s)} \int_0^\infty \frac{u^{s-1} \log u}{e^u + 1} du$$

Therefore we have

$$-\int_0^\infty \frac{u^{s-1} \log u}{e^u + 1} du = \Gamma(s)\varsigma_a'(s) - \psi(s)\varsigma_a(s)$$

See also the analysis following (C.61) in Volume VI.

Let us now consider a generalised version of Theorem 4.2.

**(ii) Theorem 4.2(a):**

(4.4.43) $$\sum_{n=1}^\infty t^n \sum_{k=1}^n \binom{n}{k} \frac{x^k}{(k+y)^s} = \frac{xt}{(1-t)\Gamma(s)} \int_0^\infty \frac{u^{s-1} e^{-(y+1)u}}{1 - (1 + xe^{-u})t} du$$

$$= \frac{1}{1-t} \Phi\left(\frac{xt}{(1-t)}, s, y\right) - \frac{1}{(1-t)y^s}$$

$$= \frac{xt}{(1-t)^2} \Phi\left(\frac{xt}{(1-t)}, s, y+1\right)$$

where $\Phi(z,s,a) = \sum_{n=0}^{\infty} \frac{z^n}{(n+a)^s}$ is the Hurwitz-Lerch zeta function defined in (4.4.44fi).

We shall see in (4.4.44) that

$$\sum_{n=0}^{\infty} t^n \sum_{k=0}^{n} \binom{n}{k} \frac{x^k}{(k+y)^s} = \frac{1}{\Gamma(s)} \int_0^{\infty} \frac{u^{s-1} e^{-yu}}{\left[1-\left(1+xe^{-u}\right)t\right]} du$$

**Proof:**

We again start with the gamma function

$$\Gamma(s) = \int_0^{\infty} t^{s-1} e^{-t} dt \quad , \operatorname{Re}(s) > 0$$

and, using the substitution $t = u(k+y)$, we obtain

(4.4.43a) $$\frac{1}{(k+y)^s} = \frac{1}{\Gamma(s)} \int_0^{\infty} u^{s-1} e^{-u(k+y)} du$$

We now consider the finite sum set out below

(4.4.43b) $$S_n(x,y) = \sum_{k=1}^{n} \binom{n}{k} \frac{x^k}{(k+y)^s}$$

and combine (4.4.43a) and (4.4.43b) to obtain

$$S_n(x,y) = \sum_{k=1}^{n} \binom{n}{k} \frac{x^k}{(k+y)^s} = \sum_{k=1}^{n} \binom{n}{k} x^k \cdot \frac{1}{\Gamma(s)} \int_0^{\infty} u^{s-1} e^{-u(k+y)} du$$

$$= \frac{1}{\Gamma(s)} \int_0^{\infty} u^{s-1} \sum_{k=1}^{n} \left\{ \binom{n}{k} e^{-u(k+y)} x^k \right\} du$$

(4.4.43c) $$= \frac{1}{\Gamma(s)} \int_0^{\infty} u^{s-1} \sum_{k=1}^{n} \left\{ \binom{n}{k} \left[e^{-u} x\right]^k \right\} e^{-yu} du$$

Applying the binomial theorem we have

(4.4.43d) $$S_n(x, y) = \frac{1}{\Gamma(s)} \int_0^\infty u^{s-1} \left\{ (1 + xe^{-u})^n - 1 \right\} e^{-yu} du$$

We now consider the following function $Z_s(x, y, t)$ which, for $|t| < 1$, is defined by

(4.4.43e) $$Z_s(x, y, t) = \frac{1}{2} \sum_{n=1}^\infty t^n \sum_{k=1}^n \binom{n}{k} \frac{x^k}{(k+y)^s}$$

and substitute the relation obtained for $S_n(x, y)$ to deduce

$$Z_s(x, y, t) = \frac{1}{2} \sum_{n=1}^\infty t^n \frac{1}{\Gamma(s)} \int_0^\infty u^{s-1} \left\{ (1 + xe^{-(u+y)})^n - 1 \right\} e^{-yu} du$$

(4.4.43f) $$= \frac{1}{2} \frac{1}{\Gamma(s)} \int_0^\infty u^{s-1} \sum_{n=1}^\infty \left\{ (1 + xe^{-u})^n - 1 \right\} t^n e^{-yu} du$$

(where we have again assumed that the conditions required for interchanging the order of summation and integration are satisfied). Use of the geometric series simplifies (4.4.43f) to

(4.4.43g) $$Z_s(x, y, t) = \frac{1}{2} \sum_{n=1}^\infty t^n \sum_{k=1}^n \binom{n}{k} \frac{x^k}{(k+y)^s} = \frac{1}{2} \frac{xt}{(1-t)\Gamma(s)} \int_0^\infty \frac{u^{s-1} e^{-(y+1)u}}{1 - (1 + xe^{-u})t} du$$

We now recall the identity in (3.67d)

$$\sum_{n=1}^\infty t^n \sum_{k=1}^n \binom{n}{k} \frac{x^k}{k^s} = \frac{1}{1-t} Li_s\left(\frac{xt}{1-t}\right)$$

and with $y = 0$ this shows that

$$\frac{xt}{\Gamma(s)} \int_0^\infty \frac{u^{s-1} e^{-u}}{1 - (1 + xe^{-u})t} du = Li_s\left(\frac{xt}{1-t}\right)$$

With $t = 1/2$ we recover (4.4.25). An alternative proof is shown below. Letting $y = 0$ in (4.4.43g) we have

(4.4.43h) $$Z_s(x, 0, t) = \frac{1}{2} \sum_{n=1}^\infty t^n \sum_{k=1}^n \binom{n}{k} \frac{x^k}{k^s} = \frac{1}{2} \frac{xt}{(1-t)\Gamma(s)} \int_0^\infty \frac{u^{s-1} e^{-u}}{1 - (1 + xe^{-u})t} du$$

and therefore

$$\sum_{n=1}^{\infty} t^n \sum_{k=1}^{n} \binom{n}{k} \frac{x^k}{k^s} = \frac{xt}{(1-t)^2 \Gamma(s)} \int_0^{\infty} \frac{u^{s-1}}{e^u - [xt/(1-t)]} du$$

With $t = 1/2$ we recover (4.4.37).

From (4.4.25) we see that

$$Li_s[xt/(1-t)] = \frac{[xt/(1-t)]}{\Gamma(s)} \int_0^{\infty} \frac{u^{s-1}}{e^u - [xt/(1-t)]} du$$

and we accordingly obtain (3.67d) again

(4.4.43i) $$\sum_{n=1}^{\infty} t^n \sum_{k=1}^{n} \binom{n}{k} \frac{x^k}{k^s} = \frac{1}{1-t} Li_s\left(\frac{xt}{1-t}\right)$$

With $x = -1$ we have

(4.4.43ia) $$\sum_{n=1}^{\infty} t^n \sum_{k=1}^{n} \binom{n}{k} \frac{(-1)^k}{k^s} = \frac{1}{1-t} Li_s\left(\frac{-t}{1-t}\right)$$

and employing Landen's identity (3.111) in the case where $s = 2$

$$Li_2\left(\frac{-t}{1-t}\right) = -\frac{1}{2}\log^2(1-t) - Li_2(t)$$

we see that

(4.4.43ib) $$\sum_{n=1}^{\infty} t^n \sum_{k=1}^{n} \binom{n}{k} \frac{(-1)^k}{k^2} = -\frac{1}{2}\frac{\log^2(1-t)}{1-t} - \frac{Li_2(t)}{1-t} = \frac{1}{1-t} Li_2\left(\frac{-t}{1-t}\right)$$

Dividing (4.4.43ib) by $t$ and then integrating gives us

$$\int_0^x \frac{Li_2(t)}{t(1-t)} dt = -\log(1-x)Li_2(x) - \log x \log^2(1-x)$$

$$- 2Li_2(1-x)\log(1-x) + 2Li_3(1-x) + Li_3(x) - 2\varsigma(3)$$

$$\int_0^x \frac{\log^2(1-t)}{t(1-t)} dt = 2\log(1-x) Li_2(1-x) + \log x \log^2(1-x)$$

$$-\frac{1}{3}\log^3(1-x) - 2Li_3(1-x) + Li_3(x) + 2\varsigma(3)$$

and we obtain

(4.4.43ic)  $\displaystyle\sum_{n=1}^{\infty} \frac{x^n}{n} \sum_{k=1}^{n} \binom{n}{k} \frac{(-1)^k}{k^2} =$

$$\log(1-x) Li_2(x) + Li_2(1-x)\log(1-x) + \frac{1}{2}\log x \log^2(1-x) + \frac{1}{6}\log^3(1-x)$$

$$-Li_3(1-x) - \frac{3}{2} Li_3(x) + \varsigma(3)$$

Using Landen's identity (3.115) in the case where $s = 3$

$$Li_3\left(\frac{-t}{1-t}\right) = \varsigma(2)\log(1-t) - \frac{1}{2}\log t \log^2(1-t) - Li_3(1-t) + \varsigma(3) + \frac{1}{6}\log^3(1-t) - Li_3(t)$$

gives us

(4.4.43id)

$$\sum_{n=1}^{\infty} t^n \sum_{k=1}^{n} \binom{n}{k} \frac{(-1)^k}{k^3} = \varsigma(2)\log(1-t) - \frac{1}{2}\log t \log^2(1-t) - Li_3(1-t) + \varsigma(3) + \frac{1}{6}\log^3(1-t) - Li_3(t)$$

Letting $t = x$ in (4.4.43i) we obtain

(4.4.43ie)  $\displaystyle\sum_{n=1}^{\infty} x^n \sum_{k=1}^{n} \binom{n}{k} \frac{x^k}{k^s} = \frac{1}{1-x} Li_s\left(\frac{x^2}{1-x}\right)$

and $x = 1/2$ gives us

(4.4.43if)  $\displaystyle\sum_{n=1}^{\infty} \frac{1}{2^{n+1}} \sum_{k=1}^{n} \binom{n}{k} \frac{1}{2^k k^s} = Li_s\left(\frac{1}{2}\right)$

We have with $t = 1/2$ in (4.4.43i)

$$\sum_{n=1}^{\infty} \frac{1}{2^{n+1}} \sum_{k=1}^{n} \binom{n}{k} \frac{x^k}{k^s} = Li_s(x)$$

which we derived by alternative means in (3.64) in Volume I.

Substituting the Flajolet and Sedgewick identity (3.16c)

$$-\sum_{k=1}^{n} \binom{n}{k} \frac{(-1)^k}{k^3} = \frac{1}{6}\left(H_n^{(1)}\right)^3 + \frac{1}{2} H_n^{(1)} H_n^{(2)} + \frac{1}{3} H_n^{(3)}$$

in (4.4.43id) gives us for $t < 1$

(4.4.43ig)

$$-\sum_{n=1}^{\infty} \left[\frac{1}{6}\left(H_n^{(1)}\right)^3 + \frac{1}{2} H_n^{(1)} H_n^{(2)} + \frac{1}{3} H_n^{(3)}\right] t^n =$$

$$\varsigma(2)\log(1-t) - \frac{1}{2}\log t \log^2(1-t) - Li_3(1-t) + \varsigma(3) + \frac{1}{6}\log^3(1-t) - Li_3(t)$$

Dividing this by $t$ and integrating will result in more expressions for Euler sums.

We have from (4.4.43g) with $t = 1/2$

(4.4.43j) $$Z_s(x, y, 1/2) = \sum_{n=1}^{\infty} \frac{1}{2^{n+1}} \sum_{k=1}^{n} \binom{n}{k} \frac{x^k}{(k+y)^s} = \frac{x}{\Gamma(s)} \int_0^{\infty} \frac{u^{s-1} e^{-yu}}{e^u - x} du$$

and, with $x = 1$ and $y \to y-1$, this becomes

(4.4.43k) $$Z_s(1, y-1, 1/2) = \sum_{n=1}^{\infty} \frac{1}{2^{n+1}} \sum_{k=1}^{n} \binom{n}{k} \frac{1}{(k+y-1)^s} = \frac{1}{\Gamma(s)} \int_0^{\infty} \frac{u^{s-1} e^{-(y-1)u}}{e^u - 1} du$$

From [126, p.92] we have the following formula for the Hurwitz zeta function

(4.4.43l) $$\varsigma(s, y) = \frac{1}{\Gamma(s)} \int_0^{\infty} \frac{u^{s-1} e^{-(y-1)u}}{e^u - 1} du$$

and hence we obtain

(4.4.43m) $$\varsigma(s, y) = \sum_{n=1}^{\infty} \frac{1}{2^{n+1}} \sum_{k=1}^{n} \binom{n}{k} \frac{1}{(k+y-1)^s}$$

$$\varsigma(s,1) = \varsigma(s) = \sum_{n=1}^{\infty} \frac{1}{2^{n+1}} \sum_{k=1}^{n} \binom{n}{k} \frac{1}{k^s}$$

Letting $y = 1$ in (4.4.43j) we get

$$Z_s(x,1,1/2) = \sum_{n=1}^{\infty} \frac{1}{2^{n+1}} \sum_{k=1}^{n} \binom{n}{k} \frac{x^k}{(k+1)^s} = \frac{x}{\Gamma(s)} \int_0^{\infty} \frac{u^{s-1} e^{-u}}{e^u - x} du$$

and using partial fractions this becomes

$$= \frac{1}{\Gamma(s)} \int_0^{\infty} \frac{u^{s-1}}{e^u - x} du - \frac{1}{\Gamma(s)} \int_0^{\infty} u^{s-1} e^{-u} du$$

Reference to (4.4.25) and (4.4.26) then shows

$$= \frac{Li_s(x)}{x} - 1$$

and therefore we get

(4.4.43ma) $\quad Z_s(x,1,1/2) = \sum_{n=1}^{\infty} \frac{1}{2^{n+1}} \sum_{k=1}^{n} \binom{n}{k} \frac{x^k}{(k+1)^s} = \frac{Li_s(x)}{x} - 1$

(4.4.43mb) $\quad Z_s(1,1,1/2) = \sum_{n=1}^{\infty} \frac{1}{2^{n+1}} \sum_{k=1}^{n} \binom{n}{k} \frac{1}{(k+1)^s} = \varsigma(s) - 1$

Hence we see that

(4.4.43mc) $\quad \sum_{n=1}^{\infty} \frac{1}{2^{n+1}} \sum_{k=1}^{n} \binom{n}{k} \left[ \frac{1}{k^s} - \frac{1}{(k+1)^s} \right] = 1$

and with $s = 2$ we have

(4.4.43md) $\quad \sum_{n=1}^{\infty} \frac{1}{2^{n+1}} \sum_{k=1}^{n} \binom{n}{k} \frac{2k+1}{k^2(k+1)^2} = 1$

Differentiation of (4.4.43mc) results in

(4.4.43me) $\quad \sum_{n=1}^{\infty} \frac{1}{2^{n+1}} \sum_{k=1}^{n} \binom{n}{k} \frac{\log k}{k^s} = \sum_{n=1}^{\infty} \frac{1}{2^{n+1}} \sum_{k=1}^{n} \binom{n}{k} \frac{\log(k+1)}{(k+1)^s}$

As before, we get with $y = 0$ in (4.4.43j)

(4.4.43n) $$Z_s(x,0,1/2) = Li_s(x) = \sum_{n=1}^{\infty} \frac{1}{2^{n+1}} \sum_{k=1}^{n} \binom{n}{k} \frac{x^k}{k^s} = \frac{x}{\Gamma(s)} \int_0^{\infty} \frac{u^{s-1}}{e^u - x} du$$

where we have recalled the series expansion for the polylogarithm in (3.55) and made reference to (4.4.25).

Letting $y = 1$ in (4.4.43g) we get

(4.4.43o) $$Z_s(x,1,t) = \frac{1}{2} \sum_{n=1}^{\infty} t^n \sum_{k=1}^{n} \binom{n}{k} \frac{x^k}{(k+1)^s} = \frac{1}{2} \frac{xt}{(1-t)\Gamma(s)} \int_0^{\infty} \frac{u^{s-1} e^{-2u}}{1-(1+xe^{-u})t} du$$

We now divide (4.4.43g) by $t$ and integrate to obtain

(4.4.43p) $$\sum_{n=1}^{\infty} \frac{q^n}{n} \sum_{k=1}^{n} \binom{n}{k} \frac{x^k}{(k+y)^s} = \int_0^q \frac{x}{(1-t)\Gamma(s)} \int_0^{\infty} \frac{u^{s-1} e^{-(y+1)u}}{1-(1+xe^{-u})t} du\, dt$$

Reversing the order of integration we get

$$\int_0^q \frac{x}{(1-t)\Gamma(s)} \int_0^{\infty} \frac{u^{s-1} e^{-(y+1)u}}{1-(1+xe^{-u})t} du\, dt = \frac{x}{\Gamma(s)} \int_0^{\infty} u^{s-1} e^{-(y+1)u}\, du \int_0^q \frac{1}{[1-(1+xe^{-u})t](1-t)} dt$$

Straightforward integration gives us

$$\int_0^q \frac{1}{[1-(1+xe^{-u})t](1-t)} dt = \frac{1}{1-A} \int_0^q \left\{ \frac{-A}{1-At} + \frac{1}{1-t} \right\} dt \text{ where } A = (1+xe^{-u})$$

$$= \frac{1}{1-A} [\log(1-Aq) - \log(1-q)]$$

and hence

$$-\frac{1}{\Gamma(s)} \int_0^{\infty} u^{s-1} e^{-yu} \left\{ \log[1-(1+xe^{-u})q] - \log(1-q) \right\} du = \sum_{n=1}^{\infty} \frac{q^n}{n} \sum_{k=1}^{n} \binom{n}{k} \frac{x^k}{(k+y)^s}$$

From (4.4.28) we have $\int_0^{\infty} u^{s-1} e^{-yu} du = \frac{\Gamma(y)}{y^s}$ and therefore we get for $y > 0$

(4.4.43q)

$$\frac{\Gamma(y)}{y^s}\log(1-q) - \int_0^\infty u^{s-1} e^{-yu} \log\left[1-\left(1+xe^{-u}\right)q\right] du = \Gamma(s) \sum_{n=1}^\infty \frac{q^n}{n} \sum_{k=1}^n \binom{n}{k} \frac{x^k}{(k+y)^s}$$

Letting $q = 1/2$ we get

(4.4.43r) $$-\int_0^\infty u^{s-1} e^{-yu} \log\left[1 - xe^{-u}\right] du = \Gamma(s) \sum_{n=1}^\infty \frac{1}{n 2^n} \sum_{k=1}^n \binom{n}{k} \frac{x^k}{(k+y)^s}$$

With the substitution $p = e^{-u}$ this becomes

(4.4.43s) $$\int_0^1 p^{y-1} \log^{s-1} p \log(1-xp) dp = (-1)^s \Gamma(s) \sum_{n=1}^\infty \frac{1}{n 2^n} \sum_{k=1}^n \binom{n}{k} \frac{x^k}{(k+y)^s}$$

and with $x = 1$ we have

(4.4.43sa) $$\int_0^1 p^{y-1} \log^{s-1} p \log(1-p) dp = (-1)^s \Gamma(s) \sum_{n=1}^\infty \frac{1}{n 2^n} \sum_{k=1}^n \binom{n}{k} \frac{1}{(k+y)^s}$$

We immediately see that

$$\frac{d^{s-1}}{dy^{s-1}} \int_0^1 p^{y-1} \log(1-p) dp = \int_0^1 p^{y-1} \log^{s-1} p \log(1-p) dp$$

and this then directs our attention to the Beta function.

From (4.4.7) we have

$$B(y,z) = \int_0^1 p^{y-1}(1-p)^{z-1} dp \quad , (\text{Re}(y) > 0, \text{Re}(z) > 0)$$

$$= \frac{\Gamma(y)\Gamma(z)}{\Gamma(y+z)}$$

Therefore we have

$$\frac{\partial}{\partial z} \frac{\partial^{s-1}}{\partial y^{s-1}} B(y,z) \bigg|_{z=1} = \int_0^1 p^{y-1} \log^{s-1} p \log(1-p) dp$$

and accordingly we get

$$\frac{\partial}{\partial z}\frac{\partial^{s-1}}{\partial y^{s-1}}\frac{\Gamma(y)\Gamma(z)}{\Gamma(y+z)}\bigg|_{z=1} = (-1)^s \Gamma(s) \sum_{n=1}^{\infty} \frac{1}{n2^n} \sum_{k=1}^{n} \binom{n}{k} \frac{1}{(k+y)^s}$$

Differentiation results in

$$\frac{\partial}{\partial y}\frac{\Gamma(y)\Gamma(z)}{\Gamma(y+z)} = \frac{\Gamma(y+z)\Gamma'(y)\Gamma(z) - \Gamma(y)\Gamma(z)\Gamma'(y+z)}{\Gamma^2(y+z)}$$

(4.4.43t) $$= B(y,z)\big[\psi(y) - \psi(y+z)\big]$$

The second derivative results in

(4.4.43u) $$\frac{\partial^2}{\partial y^2}\frac{\Gamma(y)\Gamma(z)}{\Gamma(y+z)} = B(y,z)\big[\psi'(y) - \psi'(y+z)\big] + B'(y,z)\big[\psi(y) - \psi(y+z)\big]$$

$$= B(y,z)\Big\{\big[\psi'(y) - \psi'(y+z)\big] + \big[\psi(y) - \psi(y+z)\big]^2\Big\}$$

Similarly we see that

(4.4.43v) $$\frac{\partial^3}{\partial y^3}\frac{\Gamma(y)\Gamma(z)}{\Gamma(y+z)} =$$

$$B(y,z)\Big\{\big[\psi(y) - \psi(y+z)\big]^3 + 3\big[\psi(y) - \psi(y+z)\big]\big[\psi'(y) - \psi'(y+z)\big] + \big[\psi''(y) - \psi''(y+z)\big]\Big\}$$

Differentiation of (4.4.43t) with respect to $z$ then gives us for $s = 2$

$$\frac{\partial}{\partial z}\frac{\partial}{\partial y}\frac{\Gamma(y)\Gamma(z)}{\Gamma(y+z)} = -B(y,z)\frac{\partial}{\partial z}\psi(y+z) + B(y,z)\big[\psi(z) - \psi(y+z)\big]\big[\psi(y) - \psi(y+z)\big]$$

and hence since $B(y,1) = \dfrac{\Gamma(y)\Gamma(1)}{\Gamma(y+1)} = \dfrac{1}{y}$ we get

$$\frac{\partial}{\partial z}\frac{\partial}{\partial y}\frac{\Gamma(y)\Gamma(z)}{\Gamma(y+z)}\bigg|_{z=1} = \frac{1}{y}\Big\{-\psi'(y+1) + \big[\psi(1) - \psi(y+1)\big]\big[\psi(y) - \psi(y+1)\big]\Big\}$$

We have

$$\psi(y)-\psi(y+1)=-\frac{1}{y}$$

$$\psi(1)-\psi(y+1)=-\gamma-\psi(y+1)$$

From (E.16) we see that $\psi'(y+1)=\sum_{k=0}^{\infty}\frac{1}{(y+1+k)^2}$ and we therefore obtain

$$\frac{\partial}{\partial z}\frac{\partial}{\partial y}\frac{\Gamma(y)\Gamma(z)}{\Gamma(y+z)}\bigg|_{z=1}=\frac{1}{y}\left\{-\sum_{k=0}^{\infty}\frac{1}{(y+1+k)^2}+\frac{\gamma+\psi(y+1)}{y}\right\}$$

and hence

(4.4.43w) $\quad \dfrac{1}{y}\left\{-\sum_{k=0}^{\infty}\dfrac{1}{(y+1+k)^2}+\dfrac{\gamma+\psi(y+1)}{y}\right\}=\sum_{n=1}^{\infty}\dfrac{1}{n2^n}\sum_{k=1}^{n}\binom{n}{k}\dfrac{1}{(k+y)^2}$

$$=\int_0^1 p^{y-1}\log p\log(1-p)dp$$

With $y=1$, and using $\psi(2)=-\gamma+1$, we obtain

$$1-\sum_{k=0}^{\infty}\frac{1}{(2+k)^2}=\sum_{n=1}^{\infty}\frac{1}{n2^n}\sum_{k=1}^{n}\binom{n}{k}\frac{1}{(k+1)^2}$$

and this may be written as

(4.4.43x) $\quad 2-\varsigma(2)=\sum_{n=1}^{\infty}\dfrac{1}{n2^n}\sum_{k=1}^{n}\binom{n}{k}\dfrac{1}{(k+1)^2}=\int_0^1 \log x\log(1-x)dx$

This is readily verified as follows: employing integration by parts we see that

$$\int \log x\log(1-x)dx=2x+(1-x)\log(1-x)+\left[-x-(1-x)\log(1-x)\right]\log x-Li_2(x)$$

and hence we get

$$\int_0^1 \log x\log(1-x)dx=2-\varsigma(2)$$

It should be noted that, having evaluated $\int \log x \log(1-x)dx$, it is easy to determine

$$\int \log^2 x \log(1-x)dx = \int \log x [\log x \log(1-x)]dx$$ with integration by parts.

In addition we may consider the integral

$$\int_0^1 x^{y-1} \log(1-x)dx = -\sum_{n=1}^{\infty} \frac{1}{n} \int_0^1 x^{n+y-1} dx = -\sum_{n=1}^{\infty} \frac{1}{n(n+y)}$$

and differentiating $s-1$ times with respect to $y$ we obtain

$$\int_0^1 x^{y-1} \log^{s-1} x \log(1-x)dx = (-1)^s (s-1)! \sum_{n=1}^{\infty} \frac{1}{n(n+y)^s}$$

By the geometric series we have

$$\sum_{r=0}^{s-1} (x+1)^r = \frac{1-(x+1)^s}{x}$$

and dividing by $(x+1)^s$ we easily see that

$$\frac{1}{(x+1)^s} \sum_{r=0}^{s-1} (x+1)^r = \sum_{r=0}^{s-1} \frac{1}{(x+1)^{s-r}} = \frac{1}{x(x+1)^s} - \frac{1}{x}$$

and accordingly we obtain the partial fraction expansion

$$\frac{1}{x(x+1)^s} = \frac{1}{x} - \sum_{r=0}^{s-1} \frac{1}{(x+1)^{s-r}} = \frac{1}{x} - \frac{1}{x+1} - \sum_{r=0}^{s-2} \frac{1}{(x+1)^{s-r}}$$

We therefore have

$$\sum_{n=1}^{\infty} \frac{1}{n(n+1)^s} = \sum_{n=1}^{\infty} \left[\frac{1}{n} - \frac{1}{n+1}\right] - \left\{\sum_{n=1}^{\infty} \frac{1}{(n+1)^s} + \sum_{n=1}^{\infty} \frac{1}{(n+1)^{s-1}} + \ldots + \sum_{n=1}^{\infty} \frac{1}{(n+1)^2}\right\}$$

$$= 1 - \sum_{r=2}^{s} [\varsigma(r)-1] = s - \sum_{r=2}^{s} \varsigma(r)$$

and hence we obtain

$$\int_0^1 \log^{s-1} x \log(1-x)dx = (-1)^s (s-1)!\left[ s - \sum_{r=2}^{s} \varsigma(r) \right]$$

Incidentally, as $s \to \infty$ we get Goldbach's theorem

$$\sum_{r=2}^{\infty} [\varsigma(r) - 1] = 1$$

as noted in [126, p.142].

The above partial fraction expansion may be generalised as follows (see [30a] for more details)

$$\frac{1}{x^t(x+\alpha)^s} = \sum_{r=0}^{t-1} \binom{s+r+1}{s-1} \frac{(-1)^r}{x^{t-r}\alpha^{s+r}} + \sum_{r=0}^{s-1} \binom{t+r+1}{t-1} \frac{(-1)^s}{\alpha^{t+r}(x+\alpha)^{s-r}}$$

and in particular we get

$$\frac{1}{n(n+y)^s} = \frac{1}{y^s}\left[\frac{1}{n} - \frac{1}{n+y}\right] - \left[\frac{1}{y(n+y)^s} + \frac{1}{y^2(n+y)^{s-1}} + \ldots + \frac{1}{y^{s-1}(n+y)^2}\right]$$

Therefore we have

$$\sum_{n=1}^{\infty} \frac{1}{n(n+y)^s} = \frac{1}{y^s}\sum_{n=1}^{\infty}\left[\frac{1}{n} - \frac{1}{n+y}\right] - \left\{\frac{1}{y}\sum_{n=1}^{\infty}\frac{1}{(n+y)^s} + \frac{1}{y^2}\sum_{n=1}^{\infty}\frac{1}{(n+y)^{s-1}} + \ldots + \frac{1}{y^{s-1}}\sum_{n=1}^{\infty}\frac{1}{(n+y)^2}\right\}$$

We have from [135, p.14]

$$\sum_{n=1}^{\infty}\left[\frac{1}{n} - \frac{1}{n+y}\right] = y\sum_{n=1}^{\infty}\frac{1}{n(n+y)} = \psi(y) + \gamma + \frac{1}{y}$$

and $\quad \displaystyle\sum_{n=1}^{\infty}\frac{1}{(n+y)^r} = \sum_{n=0}^{\infty}\frac{1}{(n+y)^r} - \frac{1}{y^r} = \varsigma(r, y) - \frac{1}{y^r}$

Therefore we obtain

$$\sum_{n=1}^{\infty}\frac{1}{n(n+y)^s} = \frac{1}{y^s}\left[\psi(y) + \gamma + \frac{1}{y}\right] - \sum_{r=2}^{s}\left[\varsigma(r, y) - \frac{1}{y^r}\right]\frac{1}{y^{s-r+2}}$$

We also have

$$\frac{\partial}{\partial z}\frac{\partial^2}{\partial y^2}\frac{\Gamma(y)\Gamma(z)}{\Gamma(y+z)} = \frac{\partial}{\partial z}B(y,z)\left\{[\psi'(y)-\psi'(y+z)]+[\psi(y)-\psi(y+z)]^2\right\}$$

$$= B(y,z)\left\{-\psi''(y+z)+2[\psi(y)-\psi(y+z)][\psi'(y)-\psi'(y+z)]\right\}$$

$$+B(y,z)[\psi(z)-\psi(y+z)]\left\{[\psi'(y)-\psi'(y+z)]+[\psi(y)-\psi(y+z)]^2\right\}$$

Therefore we see that

$$\left.\frac{\partial}{\partial z}\frac{\partial^2}{\partial y^2}\frac{\Gamma(y)\Gamma(z)}{\Gamma(y+z)}\right|_{z=1} =$$

$$\frac{1}{y}\left\{-\psi''(y+1)-2[\psi(y)-\psi(y+1)]\psi'(y+1)\right\}$$

$$+\frac{1}{y}[\psi(1)-\psi(y+1)]\left\{[\psi'(y)-\psi'(y+1)]+[\psi(y)-\psi(y+1)]^2\right\}$$

$$= \frac{2}{y}\left\{\sum_{k=0}^{\infty}\frac{1}{(y+1+k)^3}+\frac{1}{y}\sum_{k=0}^{\infty}\frac{1}{(y+1+k)^2}-\frac{1}{y^2}[\gamma+\psi(y+1)]\right\}$$

Therefore we obtain

(4.4.43y)

$$\frac{1}{y}\left\{\sum_{k=0}^{\infty}\frac{1}{(y+1+k)^3}+\frac{1}{y}\sum_{k=0}^{\infty}\frac{1}{(y+1+k)^2}-\frac{1}{y^2}[\gamma+\psi(y+1)]\right\} = -\sum_{n=1}^{\infty}\frac{1}{n2^n}\sum_{k=1}^{n}\binom{n}{k}\frac{1}{(k+y)^3}$$

With $y=1$, and using $\psi(2)=-\gamma+1$, we obtain

$$1-\left\{\sum_{k=0}^{\infty}\frac{1}{(k+2)^3}+\sum_{k=0}^{\infty}\frac{1}{(k+2)^3}\right\} = \sum_{n=1}^{\infty}\frac{1}{n2^n}\sum_{k=1}^{n}\binom{n}{k}\frac{1}{(k+1)^3}$$

or alternatively

$$3-\varsigma(2)-\varsigma(3) = \sum_{n=1}^{\infty}\frac{1}{n2^n}\sum_{k=1}^{n}\binom{n}{k}\frac{1}{(k+1)^3}$$

We may also obtain (4.4.43y) by differentiating (4.4.43w) whence we get:

$$-\frac{1}{y^2}\left\{-\sum_{k=0}^{\infty}\frac{1}{(y+1+k)^2}+\frac{\gamma+\psi(y+1)}{y}\right\}+\frac{1}{y}\left\{2\sum_{k=0}^{\infty}\frac{1}{(y+1+k)^3}+\frac{y\psi'(y+1)-\gamma-\psi(y+1)}{y^2}\right\}=$$

$$-2\sum_{n=1}^{\infty}\frac{1}{n2^n}\sum_{k=1}^{n}\binom{n}{k}\frac{1}{(k+y)^3}$$

In (4.4.43i) let $t \to -t$ and $x = -1$ divide by $t$ and integrate to obtain

(4.4.43z) $$\sum_{n=1}^{\infty}(-1)^n\frac{u^n}{n}\sum_{k=1}^{n}\binom{n}{k}\frac{(-1)^k}{k^s}=\int_0^u\frac{1}{t(1+t)}Li_s\left(\frac{t}{1+t}\right)dt$$

Reference to Volume I shows that

$$\int_0^u\frac{1}{t(1+t)}Li_s\left(\frac{t}{1+t}\right)dt = Li_{s+1}\left(\frac{u}{1+u}\right)$$

and we accordingly obtain

(4.4.43za) $$Li_{s+1}\left(\frac{u}{1+u}\right)=\sum_{n=1}^{\infty}(-1)^n\frac{u^n}{n}\sum_{k=1}^{n}\binom{n}{k}\frac{(-1)^k}{k^s}$$

With, for example, $u = 1/2$ (4.4.43za) becomes

(4.4.43zb) $$Li_{s+1}(1/3)=\sum_{n=1}^{\infty}\frac{(-1)^n}{n2^n}\sum_{k=1}^{n}\binom{n}{k}\frac{(-1)^k}{k^s}$$

and compare this with (4.4.43if).

Letting $t = \frac{u}{1+u}$, we have

(4.4.43zc) $$Li_{s+1}(t)=\sum_{n=1}^{\infty}\frac{1}{n}(-1)^n\left(\frac{t}{1-t}\right)^n\sum_{k=1}^{n}\binom{n}{k}\frac{(-1)^k}{k^s}$$

Note that in (4.4.45) we have shown that

$$Li_{s+1}(t)=\sum_{n=1}^{\infty}\frac{1}{n2^n}\sum_{k=1}^{n}\binom{n}{k}\frac{t^k}{k^s}$$

and we therefore have

(4.4.43zd) $$Li_{s+1}(t) = \sum_{n=1}^{\infty} \frac{1}{n}(-1)^n \left(\frac{t}{1-t}\right)^n \sum_{k=1}^{n} \binom{n}{k}\frac{(-1)^k}{k^s} = \sum_{n=1}^{\infty} \frac{1}{n2^n} \sum_{k=1}^{n} \binom{n}{k}\frac{t^k}{k^s}$$

Letting $s = 1$ in (4.4.43za) and noting that $\sum_{k=1}^{n} \binom{n}{k}\frac{(-1)^{k+1}}{k} = H_n^{(1)}$ we get

(4.4.43ze) $$Li_2\left(\frac{u}{1+u}\right) = \sum_{n=1}^{\infty} (-1)^n \frac{H_n^{(1)}}{n} u^n$$

We now recall (3.31) and (3.111a)

$$\frac{1}{2}\log^2(1-u) + Li_2(u) = \sum_{n=1}^{\infty} \frac{H_n^{(1)}}{n} u^n \qquad , |u| < 1$$

$$Li_2\left(\frac{u}{1+u}\right) = -\frac{1}{2}\log^2(1+u) - Li_2(-u)$$

and note the equivalence via (3.111b).

Alternatively, letting $s = 2$ in (4.4.43za) and noting that

$$\sum_{k=1}^{n} \binom{n}{k}\frac{(-1)^k}{k^2} = \frac{1}{2}\left[\left(H_n^{(1)}\right)^2 + H_n^{(2)}\right]$$

we get

(4.4.43zf) $$Li_3\left(\frac{u}{1+u}\right) = \frac{1}{2}\sum_{n=1}^{\infty} \frac{(-1)^n}{n}\left[\left(H_n^{(1)}\right)^2 + H_n^{(2)}\right] u^n$$

With $u = 1$ we obtain

(4.4.43zg) $$Li_3(1/2) = \frac{1}{2}\sum_{n=1}^{\infty} \frac{(-1)^n}{n}\left[\left(H_n^{(1)}\right)^2 + H_n^{(2)}\right]$$

We also note from (3.46a) that

$$\sum_{n=1}^{\infty} \frac{(-1)^n}{n}\left(H_n^{(1)}\right)^2 u^n = -\frac{1}{3}\log^3(1+u) + Li_3(-u) - Li_2(-u)\log(1+u)$$

and (3.34) gives us

$$\int_0^u \frac{Li_2(x)}{x(1-x)} dx = \sum_{n=1}^{\infty} \frac{H_n^{(2)}}{n} u^n$$

We have

$$\int_0^u \frac{Li_2(x)}{x(1-x)} dx = \int_0^u \frac{Li_2(x)}{x} dx + \int_0^u \frac{Li_2(x)}{1-x} dx$$

$$\int_0^u \frac{Li_2(x)}{1-x} dx = -\log(1-x) Li_2(x) \Big|_0^u - \int_1^u \frac{\log^2(1-x)}{x} dx$$

The final integral is evaluated in (4.4.100gii) in Volume III and we see that

$$\int_0^u \frac{Li_2(x)}{x(1-x)} dx = Li_3(u) - \log(1-u) Li_2(u) - \log u \log^2(1-u) - 2\log(1-u) Li_2(1-u) + 2Li_3(1-u) - 2\varsigma(3)$$

and therefore we have

(4.4.43zh)

$$\sum_{n=1}^{\infty} \frac{H_n^{(2)}}{n} u^n = Li_3(u) - \log(1-u) Li_2(u) - \log u \log^2(1-u) - 2\log(1-u) Li_2(1-u) + 2Li_3(1-u) - 2\varsigma(3)$$

Letting $u \to -u$ we have an expression involving complex numbers

$$\sum_{n=1}^{\infty} (-1)^n \frac{H_n^{(2)}}{n} u^n = Li_3(-u) - \log(1+u) Li_2(-u) - \log(-u) \log^2(1+u)$$

$$- 2\log(1+u) Li_2(1+u) + 2Li_3(1+u) - 2\varsigma(3)$$

From (4.4.44ic) we have

$$\frac{1}{2} \sum_{n=1}^{\infty} \frac{u^n}{n} \left[ \left(H_n^{(1)}\right)^2 + H_n^{(2)} \right] = -\frac{1}{2} \log^2(1-u) \log u - \log(1-u) Li_2(1-u) + Li_3(1-u) - \varsigma(3) - Li_3(u)$$

The second leg of (4.4.43) was demonstrated in (3.67e) in Volume I where we showed that

$$\sum_{n=1}^{\infty} t^n \sum_{k=1}^{n} \binom{n}{k} \frac{x^k}{(k+y)^s} = \frac{1}{1-t} \sum_{n=1}^{\infty} \frac{1}{(n+y)^s} \left[ \frac{xt}{(1-t)} \right]^n$$

$$= \frac{1}{1-t} \sum_{n=0}^{\infty} \frac{1}{(n+y)^s} \left[\frac{xt}{(1-t)}\right]^n - \frac{1}{(1-t)y^s}$$

and therefore

$$\sum_{n=1}^{\infty} t^n \sum_{k=1}^{n} \binom{n}{k} \frac{x^k}{(k+y)^s} = \frac{1}{1-t} \Phi\left(\frac{xt}{(1-t)}, s, y\right) - \frac{1}{(1-t)y^s}$$

The following identity [75aa] is easily derived

$$z\,\Phi(z,s,y+1) = \Phi(z,s,y) - \frac{1}{y^s}$$

and therefore we have

$$\frac{1}{1-t}\left[\Phi\left(\frac{xt}{(1-t)}, s, y\right) - \frac{1}{y^s}\right] = \frac{xt}{(1-t)^2} \Phi\left(\frac{xt}{(1-t)}, s, y+1\right)$$

Hence we obtain

$$\sum_{n=1}^{\infty} t^n \sum_{k=1}^{n} \binom{n}{k} \frac{x^k}{(k+y)^s} = \frac{xt}{(1-t)^2} \Phi\left(\frac{xt}{(1-t)}, s, y+1\right)$$

Letting $t = 1/2$ we get

(4.4.43zi) $\quad \displaystyle\sum_{n=1}^{\infty} \frac{1}{2^{n+1}} \sum_{k=1}^{n} \binom{n}{k} \frac{x^k}{(k+y)^s} = x\,\Phi(x, s, y+1)$

With $x = -1$ and $y = 0$ we get

$$\sum_{n=1}^{\infty} t^n \sum_{k=1}^{n} \binom{n}{k} \frac{(-1)^k}{k^s} = \frac{-t}{(1-t)^2} \Phi\left(\frac{-t}{1-t}, s, 1\right)$$

and since $Li_s(z) = z\,\Phi(z,s,1)$ we have another proof of (4.4.43ia)

$$\sum_{n=1}^{\infty} t^n \sum_{k=1}^{n} \binom{n}{k} \frac{(-1)^k}{k^s} = \frac{1}{1-t} Li_s\left(\frac{-t}{1-t}\right)$$

From (4.4.43sa) we see that with $y = 0$ we have

$$\int_0^1 \frac{\log^{s-1} p}{p} \log(1-p)\,dp = (-1)^s \Gamma(s) \sum_{n=1}^{\infty} \frac{1}{n 2^n} \sum_{k=1}^{n} \binom{n}{k} \frac{1}{k^s}$$

We note from (4.4.45) that

$$Li_{s+1}(x) = \sum_{n=1}^{\infty} \frac{1}{n 2^n} \sum_{k=1}^{n} \binom{n}{k} \frac{x^k}{k^s}$$

and hence we get

$$\int_0^1 \frac{\log^{s-1} p}{p} \log(1-p)\,dp = (-1)^s \Gamma(s) \varsigma(s+1)$$

which concurs with (4.4.91ga). Similarly, using (4.4.43s) we note that

(4.4.43zj) $\qquad \int_0^1 \frac{\log^{s-1} p}{p} \log(1-xp)\,dp = (-1)^s \Gamma(s) Li_{s+1}(x)$

in accordance with (4.4.38e).

Differentiating the above with respect to $x$ results in

(4.4.43zk) $\qquad \int_0^1 \frac{\log^{s-1} p}{1-xp}\,dp = (-1)^{s+1} \Gamma(s) \frac{Li_s(x)}{x}$

which agrees with (4.4.38b). A further differentiation results in

(4.4.43zl) $\qquad \int_0^1 \frac{p \log^{s-1} p}{(1-xp)^2}\,dp = (-1)^{s+1} \Gamma(s) \frac{Li_{s-1}(x) - Li_s(x)}{x^2}$

Integration of (4.4.43zj) gives us

$$\int_0^u \frac{dx}{x} \int_0^1 \frac{\log^{s-1} p}{p} \log(1-xp)\,dp = (-1)^s \Gamma(s) Li_{s+2}(u)$$

and we therefore have

(4.4.43zm) $$\int_0^1 \frac{\log^{s-1} p \, Li_2(pu)}{p} dp = (-1)^{s+1}\Gamma(s)Li_{s+2}(u)$$

in accordance with (4.4.42i).

$\square$

We easily determine that

$$\frac{\partial^2}{\partial y^2} \frac{\Gamma(y)\Gamma(z)}{\Gamma(y+z)} = B(y,z)\left\{[\psi'(y)-\psi'(y+z)]+[\psi(y)-\psi(y+z)]^2\right\}$$

$$\frac{\partial^2}{\partial y^2} \frac{\Gamma(y)\Gamma(z)}{\Gamma(y+z)}\bigg|_{y=1} = \frac{1}{z}\left\{[\psi'(1)-\psi'(1+z)]+[\psi(1)-\psi(1+z)]^2\right\}$$

$$= \frac{[\psi'(1)-\psi'(1+z)]}{z} + \frac{[\psi(1)-\psi(1+z)]}{z}[\psi(1)-\psi(1+z)]$$

As noted in [30a, p.34] we may consider this in the case where $z \to 0$: we have $B(y,1) = 1/y$.

We therefore have the limit

$$\lim_{z \to 0} \frac{\partial^2}{\partial y^2} \frac{\Gamma(y)\Gamma(z)}{\Gamma(y+z)}\bigg|_{y=1} = -\psi''(1) = 2\varsigma(3)$$

and we also see that

$$\frac{\partial^2}{\partial y^2} \frac{\Gamma(y)\Gamma(z)}{\Gamma(y+z)} = \int_0^1 p^{y-1} \log^2 p (1-p)^{z-1} dp$$

$$\lim_{z \to 0} \frac{\partial^2}{\partial y^2} \frac{\Gamma(y)\Gamma(z)}{\Gamma(y+z)}\bigg|_{y=1} = \int_0^1 \frac{\log^2 p}{(1-p)} dp$$

Therefore we see that

$$\int_0^1 \frac{\log^2 p}{(1-p)} dp = 2\varsigma(3)$$

Let us now consider another generalisation of Theorem 4.2.

**(ii) Theorem 4.2(b):**

We have for $s > 0$ and $|t| < 1$

(4.4.44) $$\sum_{n=0}^{\infty} t^n \sum_{k=0}^{n} \binom{n}{k} \frac{x^k}{(k+y)^s} = \frac{1}{\Gamma(s)} \int_{0}^{\infty} \frac{u^{s-1} e^{-yu}}{\left[1 - (1 + xe^{-u})t\right]} du$$

where the summation starts at $n = 0$ and $k = 0$ in this case.

**Proof:**

As before we have

(4.4.44a) $$\frac{1}{(k+y)^s} = \frac{1}{\Gamma(s)} \int_{0}^{\infty} u^{s-1} e^{-u(k+y)} du$$

We now consider the finite sum set out below (where the summation now starts at $k = 0$ and we specify that $y$ is neither zero nor a negative integer)

(4.4.44b) $$S_n^0(x, y) = \sum_{k=0}^{n} \binom{n}{k} \frac{x^k}{(k+y)^s}$$

(where the superscript highlights the fact that the summation starts this time at $k = 0$)

Now combine (4.4.44a) and (4.4.44b) to obtain

$$S_n^0(x, y) = \sum_{k=0}^{n} \binom{n}{k} \frac{x^k}{(k+y)^s} = \sum_{k=0}^{n} \binom{n}{k} x^k \cdot \frac{1}{\Gamma(s)} \int_{0}^{\infty} u^{s-1} e^{-u(k+y)} du$$

$$= \frac{1}{\Gamma(s)} \int_{0}^{\infty} u^{s-1} \sum_{k=0}^{n} \left\{ \binom{n}{k} e^{-u(k+y)} x^k \right\} du$$

$$= \frac{1}{\Gamma(s)} \int_{0}^{\infty} u^{s-1} \sum_{k=0}^{n} \left\{ \binom{n}{k} \left[ e^{-u} x \right]^k \right\} e^{-yu} du$$

Applying the binomial theorem we have

(4.4.44c) $$S_n^0(x, y) = \frac{1}{\Gamma(s)} \int_{0}^{\infty} u^{s-1} (1 + xe^{-u})^n e^{-yu} du$$

We now consider the following function $Z_s^0(x, y, t)$ which, for $|t| < 1$, is defined by the series

(4.4.44d) $\quad Z_s^0(x, y, t) = \dfrac{1}{2} \sum_{n=0}^{\infty} t^n \sum_{k=0}^{n} \binom{n}{k} \dfrac{x^k}{(k+y)^s}$

(where the superscript again highlights the fact that the summation starts this time at $n = 0$). We then substitute the relation obtained for $S_n^0(x, y)$ to deduce

$$Z_s^0(x, y, t) = \dfrac{1}{2} \sum_{n=0}^{\infty} t^n \dfrac{1}{\Gamma(s)} \int_0^{\infty} u^{s-1} \left(1 + xe^{-u}\right)^n e^{-yu} du$$

(4.4.44e) $\quad\quad\quad = \dfrac{1}{2} \dfrac{1}{\Gamma(s)} \int_0^{\infty} u^{s-1} \sum_{n=0}^{\infty} \left(1 + xe^{-u}\right)^n t^n e^{-yu} du$

(where we have again assumed that interchanging the order of summation and integration is permitted). Use of the geometric series simplifies (4.4.44e) to

(4.4.44f) $\quad Z_s^0(x, y, t) = \dfrac{1}{2} \sum_{n=0}^{\infty} t^n \sum_{k=0}^{n} \binom{n}{k} \dfrac{x^k}{(k+y)^s} = \dfrac{1}{2} \dfrac{1}{\Gamma(s)} \int_0^{\infty} \dfrac{u^{s-1} e^{-yu}}{\left[1 - \left(1 + xe^{-u}\right)t\right]} du$

We have for $t = 1/2$

$$Z_s^0(x, y, 1/2) = \sum_{n=0}^{\infty} \dfrac{1}{2^{n+1}} \sum_{k=0}^{n} \binom{n}{k} \dfrac{x^k}{(k+y)^s} = \dfrac{1}{\Gamma(s)} \int_0^{\infty} \dfrac{u^{s-1} e^{-(y-1)u}}{e^u - x} du$$

We then note from [126, p.121] that the Hurwitz-Lerch zeta function $\Phi(x, s, y)$ is defined by

(4.4.44fi) $\quad \Phi(x, s, y) = \sum_{n=0}^{\infty} \dfrac{x^n}{(n+y)^s} = \dfrac{1}{\Gamma(s)} \int_0^{\infty} \dfrac{u^{s-1} e^{-u(y-1)}}{e^u - x} du$

and hence we see that

(4.4.44fii) $\quad Z_s^0(x, y, 1/2) = \sum_{n=0}^{\infty} \dfrac{1}{2^{n+1}} \sum_{k=0}^{n} \binom{n}{k} \dfrac{x^k}{(k+y)^s} = \Phi(x, s, y) = \sum_{n=0}^{\infty} \dfrac{x^n}{(n+y)^s}$

With $y = 1$ we get

$$Z_s^0(x,1,1/2) = \sum_{n=0}^{\infty} \frac{1}{2^{n+1}} \sum_{k=0}^{n} \binom{n}{k} \frac{x^k}{(k+1)^s} = \frac{1}{\Gamma(s)} \int_0^{\infty} \frac{u^{s-1}}{e^u - x} du = Li_s(x)$$

With partial fractions we obtain

$$\int_0^{\infty} \frac{u^{s-1} e^{-2u}}{e^u - x} du = \frac{1}{x} \int_0^{\infty} \frac{u^{s-1}}{e^{2u}} du + \frac{1}{x^2} \int_0^{\infty} \frac{u^{s-1}}{e^u} du + \frac{1}{x^2} \int_0^{\infty} \frac{u^{s-1}}{e^u - x} du$$

$$= \frac{1}{x} \Gamma(s) \frac{1}{2^s} + \frac{1}{x^2} \Gamma(s) + \frac{1}{x^3} Li_s(x) \Gamma(s)$$

Therefore we get

$$\frac{1}{\Gamma(s)} \int_0^{\infty} \frac{u^{s-1} e^{-2u}}{e^u - x} du = \frac{1}{x} \frac{1}{2^s} + \frac{1}{x^2} + \frac{1}{x^3} Li_s(x)$$

We have

$$\Phi(x,s,a) = \sum_{n=0}^{\infty} \frac{x^n}{(n+a)^s} = \frac{1}{\Gamma(s)} \int_0^{\infty} \frac{u^{s-1} e^{-(a-1)u}}{e^u - x} du$$

and with $a = 3$ we get

$$\frac{1}{\Gamma(s)} \int_0^{\infty} \frac{u^{s-1} e^{-2u}}{e^u - x} du = \Phi(x,s,3) = \sum_{n=0}^{\infty} \frac{x^n}{(n+3)^s}$$

$$\sum_{n=0}^{\infty} \frac{x^n}{(n+3)^s} = \frac{1}{x} \frac{1}{2^s} + \frac{1}{x^2} + \frac{1}{x^3} Li_s(x)$$

$$\sum_{n=0}^{\infty} \frac{x^n}{(n+3)^s} = \frac{1}{3^s} + \frac{x}{4^s} + \frac{x^2}{5^s} + \frac{x^3}{6^s} + \ldots = \frac{1}{x^3} \left\{ \frac{x^3}{3^s} + \frac{x^4}{4^s} + \frac{x^5}{5^s} + \frac{x^6}{6^s} + \ldots \right\}$$

$$= \frac{1}{x^3} Li_3(x) - \frac{1}{x^3} \left\{ \frac{x}{1^s} + \frac{x^2}{2^s} \right\}$$

$$= \frac{1}{x^3} Li_3(x) - \frac{1}{x^2} - \frac{1}{x} \frac{1}{2^s}$$

We note from (4.4.44fi) (and we will also see this in (4.4.79)) that

$$\sum_{n=0}^{\infty}\frac{1}{2^{n+1}}\sum_{k=0}^{n}\binom{n}{k}\frac{(-1)^k}{(k+y)^s} = \sum_{n=0}^{\infty}\frac{(-1)^n}{(n+y)^s} = \Phi(-1,s,y)$$

where $\Phi(z,s,y)$ is the Hurwitz-Lerch zeta function. With $y=1$ this becomes the Hasse/Sondow identity (3.11)

$$\varsigma_a(s) = \sum_{n=0}^{\infty}\frac{1}{2^{n+1}}\sum_{k=0}^{n}\binom{n}{k}\frac{(-1)^k}{(k+1)^s} = \frac{1}{\Gamma(s)}\int_0^{\infty}\frac{u^{s-1}}{e^u+1}du$$

Differentiating (4.4.44f) with respect to $s$ we get

(4.4.44g)

$$\sum_{n=0}^{\infty}t^n\sum_{k=0}^{n}\binom{n}{k}\frac{x^k\log(k+y)}{(k+y)^s} = \frac{\Gamma'(s)}{[\Gamma(s)]^2}\int_0^{\infty}\frac{u^{s-1}e^{-yu}}{1-(1+xe^{-u})t}du - \frac{1}{\Gamma(s)}\int_0^{\infty}\frac{u^{s-1}e^{-yu}\log u}{1-(1+xe^{-u})t}du$$

Letting $t=1/2$, $s=1$ and $x=-1$ we have

$$\sum_{n=0}^{\infty}\frac{1}{2^{n+1}}\sum_{k=0}^{n}\binom{n}{k}\frac{(-1)^k\log(k+y)}{(k+y)} = -\gamma\int_0^{\infty}\frac{e^{-u(y-1)}}{e^u+1}du - \int_0^{\infty}\frac{e^{-u(y-1)}\log u}{e^u+1}du$$

Letting $y=1$ we get from the Hasse/Sondow identity

$$\varsigma_a'(1) = \sum_{n=0}^{\infty}\frac{1}{2^{n+1}}\sum_{k=0}^{n}\binom{n}{k}\frac{(-1)^k\log(k+1)}{(k+1)} = -\gamma\int_0^{\infty}\frac{1}{e^u+1}du - \int_0^{\infty}\frac{\log u}{e^u+1}du$$

We easily see that

$$\int_0^{\infty}\frac{1}{e^u+1}du = -\int_0^{\infty}\frac{-e^{-u}}{1+e^{-u}}du = -\log(1+e^{-u})\Big|_0^{\infty} = \log 2$$

Since $\varsigma_a(s) = \frac{1}{\Gamma(s)}\int_0^{\infty}\frac{u^{s-1}}{e^u+1}du$ we also see that

$$\varsigma_a(1) = \int_0^{\infty}\frac{1}{e^u+1}du = \log 2$$

We note that differentiating (3.11) results in (see also (4.4.112x))

$$\sum_{n=0}^{\infty}\frac{1}{2^{n+1}}\sum_{k=0}^{n}\binom{n}{k}\frac{(-1)^k \log(k+1)}{(k+1)^s} = -\varsigma_a'(s)$$

and hence

$$\sum_{n=0}^{\infty}\frac{1}{2^{n+1}}\sum_{k=0}^{n}\binom{n}{k}\frac{(-1)^k \log(k+1)}{k+1} = -\varsigma_a'(1)$$

We therefore see that

(4.4.44h) $\quad \displaystyle\int_0^{\infty}\frac{\log u}{e^u+1}du = -\gamma\log 2 + \varsigma_a'(1)$

We know from (C.61) of Volume VI that

$$\varsigma_a'(1) = \log 2\left[\gamma - \frac{\log 2}{2}\right]$$

and hence we have

(4.4.44i) $\quad \displaystyle\int_0^{\infty}\frac{\log u}{e^u+1}du = -\frac{\log^2 2}{2}$

It should be noted that the above result could have been obtained directly from (4.4.42c) by letting $s=1$.

Making the substitution $x = e^{-u}$ we get

$$\int_0^{\infty}\frac{\log u}{e^u+1}du = \int_0^1\frac{\log\log(1/x)}{1+x}dx$$

As shown in Appendix C of Volume VI, Adamchik [2a] has considered logarithmic integrals of this type and has shown that

$$\int_0^1\frac{x^{p-1}}{1+x^n}\log\log\left(\frac{1}{x}\right)dx = \frac{\gamma+\log(2n)}{2n}\left[\psi\left(\frac{p}{2n}\right)-\psi\left(\frac{n+p}{2n}\right)\right]+\frac{1}{2n}\left[\varsigma'\left(1,\frac{p}{2n}\right)-\varsigma'\left(1,\frac{n+p}{2n}\right)\right]$$

and in particular we have another proof of the above integral

$$\int_0^1 \frac{1}{1+x}\log\log\left(\frac{1}{x}\right)dx = -\frac{\log^2 2}{2}$$

A further differentiation of (4.4.44g) with respect to $s$ gives us

$$-\sum_{n=0}^{\infty} t^n \sum_{k=0}^{n} \binom{n}{k} \frac{x^k \log^2(k+y)}{(k+y)^s} =$$

$$\frac{[\Gamma(s)]^2 \Gamma''(s) - 2\Gamma(s)[\Gamma'(s)]^2}{[\Gamma(s)]^4} \int_0^{\infty} \frac{u^{s-1} e^{-yu}}{1-(1+xe^{-u})t} du + 2\frac{\Gamma'(s)}{[\Gamma(s)]^2} \int_0^{\infty} \frac{u^{s-1} e^{-yu} \log u}{1-(1+xe^{-u})t} du$$

$$-\frac{1}{\Gamma(s)} \int_0^{\infty} \frac{u^{s-1} e^{-yu} \log^2 u}{1-(1+xe^{-u})t} du$$

Letting $t = 1/2$, $s = 1$, $y = 1$ and $x = -1$ we have

$$-\sum_{n=0}^{\infty} \frac{1}{2^{n+1}} \sum_{k=0}^{n} \binom{n}{k} \frac{(-1)^k \log^2(k+1)}{(k+1)} = \{\Gamma''(1) - 2[\Gamma'(1)]^2\} \int_0^{\infty} \frac{du}{e^u+1} + 2\Gamma'(1) \int_0^{\infty} \frac{\log u}{e^u+1} du - \int_0^{\infty} \frac{\log^2 u}{e^u+1} du$$

From (E.16d) we have $\Gamma''(1) = \gamma^2 + \varsigma(2)$ and hence we obtain

$$\sum_{n=0}^{\infty} \frac{1}{2^{n+1}} \sum_{k=0}^{n} \binom{n}{k} \frac{(-1)^k \log^2(k+1)}{(k+1)} = [\gamma^2 - \varsigma(2)] \log 2 - \gamma \log^2 2 + \int_0^{\infty} \frac{\log^2 u}{e^u+1} du$$

Therefore we get

(4.4.44j) $$\int_0^{\infty} \frac{\log^2 u}{e^u+1} du = \varsigma_a''(1) - [-\gamma^2 + \varsigma(2) + \gamma \log 2] \log 2$$

Making the substitution $x = e^{-u}$ we get

$$\int_0^{\infty} \frac{\log^2 u}{e^u+1} du = \int_0^1 \frac{1}{1+x} \left[\log\log\left(\frac{1}{x}\right)\right]^2 dx$$

This type of integral is examined in Appendix C of Volume VI where we obtained the same result in a different way.

We note from [126, p.121] that the Hurwitz-Lerch zeta function $\Phi(x, s, y)$ is defined by

$$\Phi(x, s, y) = \sum_{n=0}^{\infty} \frac{x^n}{(n+y)^s} = \frac{1}{\Gamma(s)} \int_0^{\infty} \frac{u^{s-1} e^{-u(y-1)}}{e^u - x} \, du$$

and hence

$$\Phi(-1, 1, y) = \sum_{n=0}^{\infty} \frac{(-1)^n}{(n+y)} = \int_0^{\infty} \frac{e^{-u(y-1)}}{e^u + 1} \, du$$

Differentiating the Hurwitz-Lerch zeta function with respect to $s$ we get

$$\frac{\partial}{\partial s} \Phi(x, s, y) = -\sum_{n=0}^{\infty} \frac{x^n \log(n+y)}{(n+y)^s} = -\frac{\Gamma'(s)}{[\Gamma(s)]^2} \int_0^{\infty} \frac{u^{s-1} e^{-u(y-1)}}{e^u - x} \, du + \frac{1}{\Gamma(s)} \int_0^{\infty} \frac{u^{s-1} e^{-u(y-1)} \log u}{e^u - x} \, du$$

and with $s = 1$ and $x = -1$ we have

$$\left. \frac{\partial}{\partial s} \Phi(x, s, y) \right|_{s=1} = -\sum_{n=0}^{\infty} \frac{(-1)^n \log(n+y)}{(n+y)} = \gamma \int_0^{\infty} \frac{e^{-u(y-1)}}{e^u + 1} \, du - \int_0^{\infty} \frac{e^{-u(y-1)} \log u}{e^u + 1} \, du$$

Since

$$\sum_{n=0}^{\infty} \frac{1}{2^{n+1}} \sum_{k=0}^{n} \binom{n}{k} \frac{(-1)^k}{(k+y)^s} = \sum_{n=0}^{\infty} \frac{(-1)^n}{(n+y)^s} = \Phi(-1, s, y)$$

we deduce that

(4.4.44k) $$\sum_{n=0}^{\infty} \frac{1}{2^{n+1}} \sum_{k=0}^{n} \binom{n}{k} \frac{(-1)^k \log(k+y)}{k+y} = \sum_{n=0}^{\infty} \frac{(-1)^n \log(n+y)}{n+y}$$

We now integrate (4.4.44f) to obtain

$$\sum_{n=0}^{\infty} \frac{q^{n+1}}{n+1} \sum_{k=0}^{n} \binom{n}{k} \frac{x^k}{(k+y)^s} = \frac{1}{\Gamma(s)} \int_0^q \int_0^{\infty} \frac{u^{s-1} e^{-yu}}{\left[1 - (1 + xe^{-u})t\right]} \, du \, dt$$

Reversing the order of integration we get

$$\frac{1}{\Gamma(s)} \int_0^q \int_0^{\infty} \frac{u^{s-1} e^{-yu}}{\left[1 - (1 + xe^{-u})t\right]} \, du \, dt = \frac{1}{\Gamma(s)} \int_0^{\infty} u^{s-1} e^{-yu} \, du \int_0^q \frac{1}{\left[1 - (1 + xe^{-u})t\right]} \, dt$$

$$= -\frac{1}{\Gamma(s)} \int_0^{\infty} \frac{u^{s-1} e^{-yu} \log\left[1 - (1 + xe^{-u})q\right]}{1 + xe^{-u}} \, du$$

Therefore we have

(4.4.44l) $$-\frac{1}{\Gamma(s)}\int_0^\infty \frac{u^{s-1}e^{-yu}\log\left[1-\left(1+xe^{-u}\right)q\right]}{1+xe^{-u}}du = \sum_{n=0}^\infty \frac{q^{n+1}}{n+1}\sum_{k=0}^n \binom{n}{k}\frac{x^k}{(k+y)^s}$$

With the substitution $z = 1+xe^{-u}$ we obtain

$$\int_0^\infty \frac{u^{s-1}e^{-yu}\log\left[1-\left(1+xe^{-u}\right)q\right]}{1+xe^{-u}}du = (-1)^{s-1}\int_1^{1+x}\frac{\left(\frac{1-z}{x}\right)^y \log^{s-1}\left(\frac{1-z}{x}\right)\log[1-zq]}{z(z-1)}dz$$

$$= \sum_{n=0}^\infty \frac{q^{n+1}}{n+1}\sum_{k=0}^n \binom{n}{k}\frac{x^k}{(k+y)^s}$$

With $q=1$, $x=-1$ we get

$$\varsigma(s+1,u) = \int_0^\infty \frac{u^s e^{-yu}}{1-e^{-u}}du = \sum_{n=0}^\infty \frac{1}{n+1}\sum_{k=0}^n \binom{n}{k}\frac{(-1)^k}{(k+1)^s}$$

Differentiating (4.4.44l) with respect to $q$ results in (4.4.44f)

$$\frac{1}{\Gamma(s)}\int_0^\infty \frac{u^{s-1}e^{-yu}}{1-\left(1+xe^{-u}\right)q}du = \sum_{n=0}^\infty q^n \sum_{k=0}^n \binom{n}{k}\frac{x^k}{(k+y)^s}$$

We note from (4.4.44) that

$$\sum_{n=0}^\infty t^n \sum_{k=0}^n \binom{n}{k}\frac{x^k}{(k+y)^s} = \frac{1}{\Gamma(s)}\int_0^\infty \frac{u^{s-1}e^{-yu}}{\left[1-\left(1+xe^{-u}\right)t\right]}du$$

and integration gives us

$$\sum_{n=0}^\infty \frac{v^n}{n+1}\sum_{k=0}^n \binom{n}{k}\frac{x^k}{(k+y)^s} = \frac{1}{\Gamma(s)}\int_0^v dt \int_0^\infty \frac{u^{s-1}e^{-yu}}{\left[1-\left(1+xe^{-u}\right)t\right]}du$$

$$= \frac{1}{\Gamma(s)}\int_0^\infty u^{s-1}e^{-yu}du \int_0^v \frac{dt}{\left[1-\left(1+xe^{-u}\right)t\right]}$$

$$= -\frac{1}{\Gamma(s)} \int_0^\infty \frac{u^{s-1} e^{-yu}}{(1+xe^{-u})} du \int_0^v \frac{-(1+xe^{-u}) dt}{\left[1-(1+xe^{-u})t\right]}$$

$$= -\frac{1}{\Gamma(s)} \int_0^\infty \frac{u^{s-1} e^{-yu}}{(1+xe^{-u})} du \int_0^v \frac{-(1+xe^{-u}) dt}{\left[1-(1+xe^{-u})t\right]}$$

We have

$$\int_0^v \frac{-(1+xe^{-u}) dt}{\left[1-(1+xe^{-u})t\right]} = \log\left[1-(1+xe^{-u})v\right]$$

and therefore

$$\sum_{n=0}^\infty \frac{v^n}{n+1} \sum_{k=0}^n \binom{n}{k} \frac{x^k}{(k+y)^s} = -\frac{1}{\Gamma(s)} \int_0^\infty \frac{u^{s-1} e^{-yu} \log\left[1-(1+xe^{-u})v\right]}{1+xe^{-u}} du$$

Letting $x \to -x$ gives us

$$\sum_{n=0}^\infty \frac{1}{n+1} \sum_{k=0}^n \binom{n}{k} (-1)^k \frac{x^k}{(k+y)^s} = -\frac{1}{\Gamma(s)} \int_0^\infty \frac{u^{s-1} e^{-yu} \log[xe^{-u}]}{1-xe^{-u}} du$$

$$= \frac{1}{\Gamma(s)} \int_0^\infty \frac{u^s e^{-yu}}{1-xe^{-u}} du - \frac{\log x}{\Gamma(s)} \int_0^\infty \frac{u^{s-1} e^{-yu}}{1-xe^{-u}} du$$

$$= \frac{1}{\Gamma(s)} \int_0^\infty \frac{u^s e^{-(y-1)u}}{e^u - x} du - \frac{\log x}{\Gamma(s)} \int_0^\infty \frac{u^{s-1} e^{-(y-1)u}}{e^u - x} du$$

having regard to (4.4.44fi)

$$\Phi(x, s, y) = \frac{1}{\Gamma(s)} \int_0^\infty \frac{u^{s-1} e^{-(y-1)u}}{e^u - x} du$$

we then obtain

$$\sum_{n=0}^\infty \frac{1}{n+1} \sum_{k=0}^n \binom{n}{k} (-1)^k \frac{x^k}{(k+y)^s} = s\, \Phi(x, s+1, y) - \log x\, \Phi(x, s, y)$$

and with $x = 1$ we have

$$\sum_{n=0}^{\infty}\frac{1}{n+1}\sum_{k=0}^{n}\binom{n}{k}\frac{(-1)^k}{(k+y)^s}=s\,\Phi(1,s+1,y)$$

Letting $s\to s-1$ gives us

$$\sum_{n=0}^{\infty}\frac{1}{n+1}\sum_{k=0}^{n}\binom{n}{k}\frac{(-1)^k}{(k+y)^{s-1}}=(s-1)\Phi(1,s,y)$$

and since

$$\Phi(x,s,y)=\sum_{n=0}^{\infty}\frac{x^n}{(n+y)^s}$$

we obtain the Hasse identity

$$\frac{1}{s-1}\sum_{n=0}^{\infty}\frac{1}{n+1}\sum_{k=0}^{n}\binom{n}{k}\frac{(-1)^k}{(k+y)^{s-1}}=\sum_{n=0}^{\infty}\frac{1}{(n+y)^s}=\varsigma(s,y)$$

Differentiation results in

$$\sum_{n=0}^{\infty}\frac{v^n}{n+1}\sum_{k=0}^{n}\binom{n}{k}\frac{x^k\log(k+y)}{(k+y)^s}=$$

$$-\frac{1}{\Gamma(s)}\int_0^{\infty}\frac{u^{s-1}e^{-yu}\log u\,\log\!\left[1-(1+xe^{-u})v\right]}{1+xe^{-u}}du+\frac{\psi(s)}{\Gamma(s)}\int_0^{\infty}\frac{u^{s-1}e^{-yu}\log\!\left[1-(1+xe^{-u})v\right]}{1+xe^{-u}}du$$

$$\sum_{n=0}^{\infty}\frac{v^n}{n+1}\sum_{k=0}^{n}\binom{n}{k}\frac{(-1)^k\log(k+y)}{(k+y)^s}=\frac{1}{\Gamma(s)}\int_0^{\infty}\frac{u^s\log u\,e^{-(y-1)u}}{e^u-1}du-\frac{\psi(s)}{\Gamma(s)}\int_0^{\infty}\frac{u^s e^{-(y-1)u}}{e^u-1}du$$

We note from (4.4.13b) that

$$\int_0^1\frac{t^{y-1}\log^s t}{1-v(1-xt)}dt=(-1)^s s!\sum_{n=0}^{\infty}v^n\sum_{k=0}^{n}\binom{n}{k}(-1)^k\frac{x^k}{(k+y)^{s+1}}$$

Donal F. Connon
Elmhurst
Dundle Road
Matfield
Kent TN12 7HD
dconnon@btopenworld.com